\documentclass[a4paper,12pt,pdftex]{article}

\usepackage{stmaryrd}

\usepackage{amssymb,amsmath}
\usepackage{amsthm,comment}
\usepackage{bm}%
\usepackage{xcolor}
\usepackage {arydshln}
\usepackage{tikz}
  \usetikzlibrary{arrows.meta}
  \usetikzlibrary{decorations.markings}
  \usetikzlibrary{calc}

\def\bsqq{~$\blacksquare$} 
 \def\cb{{\cal B}} 
  
\def\cg{{\cal G}} \def\ch{{\cal H}} 
 \def\ck{{\cal K}} 
  
\def\cp{{\cal P}}   \def\cs{{\cal S}}

\def\bbc{\mathbb{C}}
\def\bbr{\mathbb{R}}
\def\bbn{\mathbb{N}}

\def\eqn#1{\begin{equation}\label{#1}}
\def\ee{\end{equation}}

\def\bea{\begin{eqnarray}}
\def\eea{\end{eqnarray}}

\def\eqnn#1{\begin{eqnarray}\label{#1}}
\newcommand{\eqna}[1]{\begin{subequations} \label{#1}
\begin{eqnarray}}
\def\eena{\end{eqnarray}
\end{subequations}}

\def\nn{\nonumber}
\def\ha{{\textstyle{1\over2}}}
\def\d{\delta}
\def\L{\Lambda} 

\def\rf#1#2{(\ref{#1}{#2})}

   \def\nt{\noindent}

\def\nn{\nonumber}

\def\ket#1{\left| #1\right\rangle}
\def\hf{\frac{1}{2}}
\def\Gau#1{\left\lfloor #1 \right\rfloor}

\begin{document}

\begin{center}
{\Large\bf
  Classification of the Reducible Verma Modules  over\\[4pt] the Jacobi Algebra $ {\cal G}_2$
}

\vspace{10mm}

{\large\bf  N. Aizawa$^1$, ~V.K. Dobrev$^2$, ~S. Doi$^1$}
	\\[10pt]
$^1$Department of Physical Science, Osaka Prefecture University, \\
Nakamozu Campus, Sakai, Osaka 599-8531, Japan
\\[10pt]
$^2$Institute of Nuclear Research and Nuclear Energy,\\
 Bulgarian Academy of Sciences, \\
72 Tsarigradsko Chaussee, 1784 Sofia, Bulgaria

\end{center}

\begin{abstract}
In the present paper we study the representations of the Jacobi algebra. More concretely, we
define, analogously to the case of semi-simple Lie algebras, the Verma modules over the Jacobi algebra
$\cg_2$. We study their reducibility and give explicit construction of the reducible Verma modules exhibiting
the corresponding singular vectors. Using this information we give a complete classification of the reducible Verma
modules. More than this we exhibit their interrelation of embeddings between these modules.
These embeddings are illustrated by diagrams of the embedding patterns so that each reducible Verma module appears in one
such diagram.
 \end{abstract}

\section{Introduction}

Invariant differential operators play a very important role in the description of physical
symmetries - recall, e.g., the examples of Dirac, Maxwell, Klein-Gordon, d'Almbert,  Schr\"odinger,
equations. Thus, when studying the applications of some symmetry object, always a very interesting question
is to find the invariant differential operators and equations with that symmetry.

The role of nonrelativistic symmetries in  theoretical physics was
always important.  Currently one of the  most popular fields in
theoretical physics - string theory, pretending to be a universal
theory - encompasses together relativistic quantum field theory,
classical gravity, and certainly, nonrelativistic quantum mechanics,
in such a way that it is not even necessary to separate these
components.

Since the cornerstone of quantum mechanics is the Schr\"o\-din\-ger
equation then
 it is not a surprise that the Schr\"o\-din\-ger group -
  the group that is the maximal group of
symmetry of the Schr\"o\-din\-ger equation - was the first to play
a prominent role in theoretical physics.
The latter is natural since originally the Schr\"odinger group,
actually the Schr\"odinger algebra, was introduced  in
\cite{Nie,Hag} as a nonrelativistic limit of
the vector-field realization of the conformal algebra.

Another interesting non-relativistic example is the Jacobi algebra \cite{EiZa,BeSc} which is the
semi-direct sum of the Heisenberg algebra and the $sp(n)$ algebra. Actually the lowest case of the
Jacobi algebra coincides with the lowest case of the Schr\"odinger algebra which makes it interesting
to apply to the Jacobi algebra the methods  applied to the Schr\"odinger algebra \cite{DDM,AiDo}.

For our approach, we recall that in the case of Lie algebras and groups there are several methods of finding the corresponding
invariant differential operators. We shall follow the method developed in \cite{Dob88,VKD4}.
In this method there is a correspondence between invariant differential operators and singular vectors of Verma modules
over the (complexified) Lie algebra in consideration.

Thus, in the present paper we develop the first stage of the project, namely, we give the complete classification of the
reducible Verma modules over the Jacoby algebra $\cg_2$. This is achieved by the study and explicit construction the singular vectors
of the Verma modules. We should note that this was started in \cite{Dob} where were given some
  examples the low level singular vectors.

The paper is organized as follows. In the next section we give the preliminaries needed.
Then we develop in detail the theory of lowest weight Verma modules over ${\cal G}_2$.
Then we study the singular vectors of Verma modules ${\cal G}_2$.
Further, we study the question of multiple reducibilities of Verma modules.
In the final section we give the classification of the reducible Verma modules. For this
we need to study the complete embedding pictures since some some reducible Verma modules
contain chains of embedded Verma modules. All results are illustrated by diagrams of the
embedding patterns.

\section{Preliminaries}

\subsection{Jacoby algebra $\cg_n$ : general case}

The Jacobi algebra is the semi-direct sum $\cg_n:= \ch_n\niplus 
sp(n,\bbr)_{\bbc}$ \cite{EiZa,BeSc}. The Heisenberg algebra
$\ch_n$ is generated by the boson creation (respectively,
annihilation) operators~${a}_i^{+}$~(${a}^-_i$),~$i,j
=1,\dots,n$, which verify the canonical commutation relations
\eqn{heis} \big[a^-_i,a^{+}_j\big]=\delta_{ij}, \qquad [a^-_i,a^-_j]
= \big[a_i^{+},a_j^{+}\big]= 0 . \ee
$\ch_n$ is an
ideal in $\cg_n$, i.e., $[\ch_n,\cg_n]=\ch_n$,
determined by the commutation relations (following the notation of \cite{Berc}):
\eqna{haspn}
&&\big[a^{+}_k,K^+_{ij}\big] = [a^-_k,K^-_{ij}]=0, \\
&&{} [a^-_i,K^+_{kj}] = \tfrac{1}{2}\delta_{ik}a^{+}_j+\tfrac{1}{2}\delta_{ij}a^{+}_k ,\qquad
 \big[K^-_{kj},a^{+}_i\big] = \tfrac{1}{2}\delta_{ik}a^-_j+\tfrac{1}{2}\delta_{ij}a^-_k , \\
&& \big[K^0_{ij},a^{+}_k\big] = \tfrac{1}{2}\delta_{jk}a^{+}_i,\qquad
\big[a^-_k,K^0_{ij}\big]= \tfrac{1}{2}\delta_{ik}a^-_{j} .
\eena
 $K^{\pm,0}_{ij}$ are the generators of the $\cs_n ~\equiv~ sp(n,\bbr)_{\bbc}$ algebra:
\eqna{baspn}
&& [K_{ij}^-,K_{kl}^-] = [K_{ij}^+,K_{kl}^+]=0 , \qquad 2\big[K^-_{ij},K^0_{kl}\big] = K_{il}^-\delta_{kj}+K^-_{jl}\delta_{ki}\label{baza23}, \\
&& 2[K_{ij}^-,K_{kl}^+] = K^0_{kj}\delta_{li}+
K^0_{lj}\delta_{ki}+K^0_{ki}\delta_{lj}+K^0_{li}\delta_{kj}\\
&& 2\big[K^+_{ij},K^0_{kl}\big] = -K^+_{ik}\delta_{jl}-K^+_{jk}\delta_{li},\quad
 2\big[K^0_{ji},K^0_{kl}\big] = K^0_{jl}\delta_{ki}-K^0_{ki}\delta_{lj} . 
\eena

In order to implement our approach we introduce a triangular decomposition of ~$\cg_n$~:
\eqn{deca} \cg_n\,=\, \cg_n^+\oplus \ck_n\oplus \cg_n^- \ ,\ee
using the triangular decomposition ~$\cs_n ~=~ \cs_n^+\oplus \ck_n\oplus \cs_n^- $,
where:
\eqnn{decpm} && \cg_n^\pm ~=~ \ch_n^\pm \oplus \cs_n^\pm \\
&& \ch^\pm_n ~=~ {\rm l.s.} \{\, {a}_i^{\pm} : i=1,\dots,n\}\ , \nn\\
&& \cs_n^+ ~=~ {\rm l.s.} \{\, K^+_{ij} ~:~ 1\leq i\leq j\leq n \}
\oplus {\rm l.s.} \{\, K^0_{ij} ~:~ 1\leq i<  j\leq n \} \nn\\
&& \cs_n^- ~=~ {\rm l.s.} \{\, K^-_{ij} ~:~ 1\leq i\leq j\leq n \}
\oplus {\rm l.s.} \{\, K^0_{ij} ~:~ 1\leq j < i \leq n \} \nn\\
&&\ck_n ~=~ {\rm l.s.} \{\, K^0_{ii} ~:~ 1\leq  i \leq n \}
\nn\eea
Note that the subalgebra ~$\ck_n$~ is abelian and
 is a Cartan subalgebra of ~$\cs_n$. Furthermore, not only ~$\cs_n^\pm$, but also
 ~$\cg_n^\pm$~ are its eigenspaces:
\eqn{ckcg} [ \ck_n, \cg_n^\pm] ~=~ \cg_n^\pm \ee
Thus, ~$\ck_n$~ plays for ~$\cg_n$~ the role that
Cartan subalgebras are playing for semi-simple Lie algebras.

Note that the algebra ~$\cg_1$~ is isomorphic to the (1+1)-dimensional Schr\"odinger algebra (without
central extension). The representations of the latter are well known, cf. \cite{DDM,AiDo,AiDo2,DLMZ}.
Thus, below in this paper we study the first new  case of the ~$\cg_n$~ series, namely, ~$\cg_2$.

\subsection{Jacoby algebra : case of ~$\cg_2$~}

First, for simplicity, we introduce the following notations for the basis of ~$\cs_2$~:
\eqna{bas2} \cs^+ ~:~&&~ b^+_i ~\equiv~ K^+_{ii}\ , ~~i=1,2; \quad
c^+ ~\equiv~ K^+_{12}\ ,\quad d^+ ~\equiv K^0_{12} \\
\cs^- ~:~&&~ b^-_i ~\equiv~ K^-_{ii}\ , ~~i=1,2; \quad
c^- ~\equiv~ K^-_{12}\ ,\quad d^- ~\equiv K^0_{21} \\
\ck ~:~&&~ h_i ~\equiv~ K^0_{ii} \ , ~~i=1,2.
\eena

Define $ \delta_1 := (\hf,0) $ and $ \delta_2 :=(0, \hf)$. Next, using \eqref{haspn} and \eqref{baspn} we give the eigenvalues of the basis of
~$\cg^+$~ w.r.t. ~$\ck$, namely, $ \mathrm{ad}h_1 $ and $ \mathrm{ad}h_2$   are expressed in terms of $\delta_1, \delta_2$ as follows:
\eqnn{eig2}
&& h_1 ~:~ (b^+_1, b^+_2, c^+, d^+ , a_1^+, a_2^+) ~:~ (1,0,\ha,\ha,\ha,0) \ , \\
&& h_2 ~:~ (b^+_1, b^+_2, c^+, d^+ , a_1^+, a_2^+) ~:~ (0,1,\ha,-\ha,0,\ha) \ , \nn\eea
(e.g., ~$[h_1,b^+_1] = b^+_1$, ~~$[h_2,d^+] = -\ha d^+$, etc).
Naturally, the  eigenvalues of the basis of ~$\cg^-$~ w.r.t. ~$\ck$~ are obtained from \eqref{eig2}
by multiplying every eigenvalue by (-1).
 
 Thus,  we can introduce the following grading of the basis of ~$\cg^+_2$:
\eqn{grad}
(b^+_1, b^+_2, c^+, d^+ , a_1^+, a_2^+) ~:~ (2\d_1,2\d_2,\d_1+\d_2,\d_1-\d_2,\d_1,\d_2)
\ee 
and the corresponding formula for $ {\cal G}_2^-$ is obtained by muliplying $-1$.

The grading of the   ~$\cs^+_2$~ part of the basis follows from the root system of  ~$\cs^+_2$,
while the grading of the ~$\ch^+_2$~ part of the basis is determined by consistency with commutation
relations \eqref{haspn}. It is consistent also with  formulae \eqref{eig2}.

 \medskip

Next we give the explicit commutation relations:

\noindent
$-$ Heisenberg algebra $ {\cal H}_2$
\begin{equation}
   [a_i^-, a_j^+ ] = \delta_{ij}, \qquad
   [a_i^{\pm}, a_j^{\pm}] = 0.
\end{equation}

\bigskip
\noindent
$-$ Symplectic algebra $ {\cal S}_2 = C_2$
\begin{alignat}{4}
  [h_1, b^{\pm}_1] &= \pm b^{\pm}_1, & \quad
  [h_1, b^{\pm}_2] &= 0, & \quad
  [h_1, c^{\pm}] &= \pm \hf c^{\pm}, & \quad
  [h_1, d^{\pm}] &= \pm \hf d^{\pm},
  \nn \\
  [h_2, b^{\pm}_1] &= 0, & \quad
  [h_2, b^{\pm}_2] &= \pm b^{\pm}_2, & \quad
  [h_2, c^{\pm}] &= \pm \hf c^{\pm}, & \quad
  [h_2, d^{\pm}] &= \mp \hf d^{\pm},
  \nn \\
  [h_1, h_2] &= 0, &
  [b_i^{\pm}, b_j^{\pm}] &= 0, &
  [b_i^-, b_j^+] &= 2 \delta_{ij} h_i, &
  [b_i^{\pm}, c^{\pm}] &= 0,
  \nn \\
  [b_1^{\pm}, c^{\mp}] &= \mp d^{\pm}, &
  [b_2^{\pm}, c^{\mp}] &= \mp d^{\mp}, &
  [b_1^{\pm}, d^{\pm}] &= 0, &
  [b_1^{\pm}, d^{\mp}] &= \mp c^{\pm},
  \nn \\
  [b_2^{\pm}, d^{\pm}] &= \mp c^{\pm}, &
  [b_2^{\pm}, d^{\mp}] &= 0, & 
  c^{\pm}, d^{\pm}] &= \mp \hf b_1^{\pm},   & 
  [c^{\pm}, d^{\mp}] &= \mp \hf b_2^{\pm}, 
  \nn  \\
  [c^-, c^+] &= \hf (h_1+h_2), &
  [d^-, d^+] &= \hf (h_2-h_1).
\end{alignat}

\bigskip
\noindent
$-$ Semidirect sum of $ {\cal H}_2$ and $ {\cal S}_2$
\begin{alignat}{4}
  [h_i, a_i^{\pm}] &= \pm \hf a_i^{\pm}, & \quad
  [a^{\pm}_i, b^{\pm}_j] &= 0, & \quad
  [a_i^{\pm}, b_i^{\mp}] &= \mp a_i^{\mp}, & \quad
  [a^{\pm}_i, c^{\pm}] &= 0,
  \nn \\
  [a_1^{\pm}, c^{\mp}] &= \mp \hf a_2^{\mp}, &
  [a_2^{\pm}, c^{\mp}] &= \mp \hf a_1^{\mp}, &
  [a_1^{\pm}, d^{\pm}] &=0, &
  [a_2^{\pm}, d^{\pm}] &= \mp \hf a_1^{\pm},
  \nn \\
  [a_1^{\pm}, d^{\mp}] &= \mp \hf a_2^{\pm}, &
  [a_2^{\pm}, d^{\mp}] &= 0
\end{alignat}
and
\begin{equation}
  [h_i, a_j^{\pm}] = [a_i^{\pm}, b_j^{\mp}] = 0 \quad \text{for} \quad i \neq j
\end{equation}

%
\section{Lowest weight Verma modules over ${\cal G}_2$} \label{SEC:Verma}
\setcounter{equation}{0}

As in \cite{Dob} we introduce Verma modules over the Jacobi algebra analogously to the case of
  of semi-simple algebras.
 Thus, we define a  lowest weight ~{\it Verma module} ~$V^\L$~ over $\cg_n$   as
the lowest  weight module over ~$\cg_n$~ with lowest  weight ~$\L \in \ck_n^*$~ and
lowest weight vector ~$v_0 \in V^\L$, induced from the
one-dimensional representation ~$V_0 \cong \bbc v_0$~ of
~$U(\cb_n)$~, (where ~$\cb_n  = \ck_n \oplus \cg_n^-$~ is a Borel
subalgebra of ~$\cg_n$), such that:
\eqnn{indb}
& &X ~v_0 ~~=~~ 0 , \quad \forall\, X\in \cg_n^- \cr
&&H ~v_0 ~~=~~ \L(H)~v_0\,, \quad \forall\, H \in \ck_n \eea

We introduce
\begin{equation}
   \hat{b}_k^+ := b_k^+ - \hf (a_k^+)^2, \qquad
   \hat{c}^+ := c^+ - \hf a_1^+ a_2^+.
\end{equation}
As a basis of $ V^{\Lambda}$  one may take ($\,\ket{0} = v_0 \,$)
\begin{equation}
  \ket{k_1, k_2, m_1, m_2, n_1, n_2} :=
  (\hat{b}_1^+)^{k_1} (\hat{b}_2^+)^{k_2} (a_1^+)^{m_1}  (a_2^+)^{m_2}   (\hat{c}^+)^{n_1} (d^+)^{n_2} \ket{0}
\end{equation}
where $ k_i, m_i, n_i \in \mathbb{Z}_{\geq 0}.$
This is an eigenvector of $ h_1, h_2 $ and its eigenvalue is given by
\begin{align}
    \Lambda &+ p_1 \delta_1 + p_2 \delta_2,  \label{weight-n}
     \\[3pt]
    & \Lambda := (\Lambda_1, \Lambda_2), \quad \Lambda_i := \Lambda(h_i)
    \nn \\
    & p_1 = K_1 + n_1 + n_2, \quad p_2 = K_2 + n_1 - n_2, \ (K_i := 2k_i+m_i).
    \label{p-def}
\end{align}
We also use the shorthand notation $ \ket{\underline{k}, \underline{m}, \underline{n}}$
with $ \underline{k} = (k_1, k_2), $ etc.

We also introduce
\begin{alignat}{3}
   \hat{d}^+ &:= d^+ - \hf a_1^+ a_2^-,
   &\qquad
   \hat{d}^- &:= d^- - \hf a_2^+ a_1^-,
   &\qquad
   N_k &:= a_k^+ a_k^-.
   \\
   \hat{b}_k^- &:= b_k^- - \frac{1}{2}(a_k^-)^2,
   &
   \hat{c}^- &:= c^- - \frac{1}{2} a_1^- a_2^-.
\end{alignat}
Next we compute the  action of $ {\cal G}_2^- $ on $\ket{\underline{k}, \underline{m}, \underline{n}}. $
The action of $ a_k^-$ is given by
\begin{align}
  a_1^- \ket{\underline{k}, \underline{m}, \underline{n}} &=
  m_1 \ket{\underline{k}, m_1-1, m_2, \underline{n}} ,
  \\
  a_2^- \ket{\underline{k}, \underline{m}, \underline{n}} &=
  m_2 \ket{\underline{k}, m_1, m_2-1, \underline{n}}.
\end{align}
The action of $ \hat{b}_k^-$ is given by
\begin{align}
  \hat{b}_1^- \ket{\underline{k}, \underline{m}, \underline{n}} &=
  k_1 \Big( 2\Lambda_1 + k_1 + n_1 + n_2 - \frac{3}{2} \Big) \ket{k_1-1,k_2, \underline{m}, \underline{n}}
  \nn \\
  & + \hf n_1 n_2 \Big(\Lambda_2-\Lambda_1-\hf(n_2-1) \Big)
  \ket{\underline{k}, \underline{m}, n_1-1, n_2-1}
  \nn \\
  &+ \frac{1}{4} n_1(n_1-1) \ket{k_1, k_2+1, \underline{m}, n_1-2, n_2},
  \\
  \hat{b}_2^- \ket{\underline{k}, \underline{m}, \underline{n}} &=
  k_2 \Big( 2\Lambda_2 + k_2 + n_1 - n_2 - \frac{3}{2} \Big) \ket{k_1,k_2-1, \underline{m}, \underline{n}}
  \nn \\
  & + n_1 \ket{\underline{k}, \underline{m}, n_1-1, n_2+1}
    + \frac{1}{4} n_1(n_1-1) \ket{k_1+1, k_2, \underline{m}, n_1-2, n_2}.
\end{align}
The action of $\hat{d}^-, \hat{c}^-$ is as follows:
\begin{align}
  \hat{d}^- \ket{\underline{k}, \underline{m}, \underline{n}}
  &= \frac{n_2}{2} \Big(\Lambda_2-\Lambda_1-\hf(n_2-1) \Big)
   \ket{\underline{k}, \underline{m}, n_1, n_2-1}
  \nn \\
  &+k_1 \ket{k_1-1, k_2, \underline{m}, n_1+1, n_2}
   + \frac{n_1}{2} \ket{k_1, k_2+1, \underline{m}, n_1-1, n_2},
  \\
  \hat{c}^- \ket{\underline{k}, \underline{m}, \underline{n}}
  &=
  \frac{n_1}{2} \Big(  \sum_{j=1}^2 (\Lambda_j + k_j) + \hf(n_1-2) \Big)
  \ket{\underline{k}, \underline{m}, n_1-1, n_2}
  \nn \\
  & + \hf k_2 n_2 \Big( \Lambda_2 - \Lambda_1 - \hf (n_2-1) \Big)
    \ket{k_1, k_2-1, \underline{m}, n_1, n_2-1}
  \nn \\
  &+ k_1 \ket{k_1-1, k_2, \underline{m}, n_1, n_2+1}
   + k_1 k_2 \ket{k_1-1, k_2-1, \underline{m}, n_1+1, n_2}.
\end{align}

The Verma module $ V^{\Lambda}$ has a weight space decomposition:
\begin{equation}\label{wsd}
  V^{\Lambda} = \bigoplus_{p_1,p_2} V^{\Lambda}_{p_1,p_2}
\end{equation}
where $ V^{\Lambda}_{p_1,p_2} $ is a subspace of $ V^{\Lambda}$ spanned by the vectors with the weight $ \Lambda + p_1 \delta_1 + p_2 \delta_2.$

\medskip
\noindent
\textbf{Remarks:}
\begin{enumerate}
  \renewcommand{\labelenumi}{(\roman{enumi})}
  \item $ p_1 \in \mathbb{Z}_{\geq 0}, \ p_2 \in \mathbb{Z}.$
  \item for a fixed value of $p_1,$ the smallest possible value of $p_2$ is $ p_2 = -p_1.$
  This is seen from the fact that the largest value of $n_2$ for a fixed $p_1$ is $n_2 = p_1.$
\end{enumerate}

%
%
%
%
\section{Singular vectors in Verma modules} \label{SEC:SVsearch}
\setcounter{equation}{0}

\subsection{Definitions and summary of singular vectors}\label{gendef}

We are interested  in the cases when the Verma modules are reducible. Namely, we are interested
in the cases when a Verma module $V^\L$  contains an invariant submodule
which is also a Verma module ~$V^{\L'}$, where ~$\L'\neq \L$, and holds the analog of
\eqna{indbb}
& &X ~\ket{0'} ~~=~~ 0 , \quad \forall\, X\in \cg_n^- \\
&&H ~\ket{0'} ~~=~~ \L'(H)~v'_0\,, \quad \forall\, H \in \ck_n 
\eena
Since ~$V^{\L'}$~ is an invariant submodule then there should be a mapping such that
~$\ket{0'}$~ is mapped to a ~{\it singular vector} ~$\ket{v_s} \in V^\L$~ fulfilling exactly \eqref{indbb}.
Thus, as in the semi-simple case there should be a polynomial ~$\cp$~ of ~$\cg_n^+$~ elements
which is eigenvector of $\ck_n$: ~$[H,\cp] = \L'(H)\cp$, ($\forall H\in\ck_n$), and then we would have: ~$\ket{v_s} ~=~\cp \ket{0}\,$.

The above situation we shall depict by the following diagram:

\begin{center}
  \begin{tikzpicture}
    \node[circle,fill=black,scale=0.3] (Lddd5) at (0,0) {};
    \node[circle,fill=black,scale=0.3] (Ldd4) at (3,0) {};
    \draw[decoration={markings, mark=at position 0.5 with {\arrow{Latex}}},postaction={decorate}] (Lddd5) to (Ldd4);
    \node[left] at (Lddd5) {$V^\Lambda$};
    \node[right] at (Ldd4) {$V^{\Lambda'}$};
  \end{tikzpicture}
\end{center}

Note that in the diagram the arrow points ~{\it to}~ the embedded Verma module.

\bigskip
Now we present the result of the complete search of singular vectors. 
We found five types of singular vectors and they exist in Verma modules with a particular value of the lowest weight. 
To specify the lowest weight, we introduced positive integers $ p^1, p^2, p^3, p^4, p^5 $ and $ q^3.$ 
\begin{enumerate}
 \renewcommand{\labelenumi}{(\roman{enumi})}
 \item $ \Lambda_1 - \Lambda_2 = \hf(1-p^1)$
 \begin{equation}
   \ket{v_s^{\Lambda'}} = (d^+)^{p^1} \ket{0}, \quad
   \Lambda' = \Lambda+p^1 (\delta_1-\delta_2). \label{SVcase1}
 \end{equation}
 \item $ {}^{\forall} \Lambda_1, \ \Lambda_2 = \frac{3}{4}- \frac{p^2}{2} $
 \begin{equation}
    \ket{v_s^{\Lambda'}} = (\hat{b}_2^+)^{p^2} \ket{0}, \quad
    \Lambda' = \Lambda + 2 p^2 \delta_2.
    \label{SVcase2}
 \end{equation}
 \item $ \Lambda_1 = \frac{5}{4}- \hf(p^3-q^3), \ \Lambda_2 = \frac{3}{4}-\hf q^3, \ (p^3 \neq q^3, p^3 \neq 2q^3)$
 
 There are three subcases depending on the range of $p^3, q^3.$ 
 The form of $ \ket{v_s^{\Lambda'}} $ and $\Lambda'$ are common for all the subcases:
 \begin{align}
     \ket{v_s^{\Lambda'}} &= \sum c(k,n) \ket{k, q^3-k-n, n, p^3-2k-n}, 
      \\
      \Lambda' &= \Lambda + p^3 \delta_1 + (2q^3-p^3) \delta_2.
 \end{align}
 They differ only the range of summation:
   \begin{enumerate}
     \item $ p^3 < q^3 \ \Rightarrow \  \displaystyle \sum = \sum_{k=0}^{\Gau{p^3/2}} \sum_{n=0}^{p^3-2k}$
     \item $ q^3 < p^3 < 2q^3 \ \Rightarrow \  \displaystyle \sum = \sum_{k=0}^{p^3-q^3} \sum_{n=0}^{q^3-k} + \sum_{k=p^3-q^3+1}^{\Gau{p^3/2}} \sum_{n=0}^{p^3-2k} $
     \item $ 2q^3 < p^3  \ \Rightarrow \  \displaystyle \sum = \sum_{k=0}^{q^3} \sum_{n=0}^{q^3-k} $
   \end{enumerate}
   \item $ \Lambda_1 + \Lambda_2 = 2 - \frac{p^4}{2}$
   \begin{align}
       \ket{v_s^{\Lambda'}} &=  \sum_{k=0}^{\Gau{p^4/2}} \sum_{n=0}^{p^4-2k}
       c(k,n) \ket{k,p^4-k-n,n,p^4-2k-n},
       \quad
       \Lambda' = \Lambda + p^4(\delta_1 + \delta_2),
   \end{align}
   \item $ \Lambda_1 = \frac{5}{4}- \frac{p^5}{2}, \ {}^{\forall} \Lambda_2 $
   \begin{align}
       \ket{v_s^{\Lambda'}} &=  \sum_{k=0}^{p^5} \sum_{n=0}^{p^5-k}
       c(k,n) \ket{k,p^5-k-n,n,2p^5-2k-n},
       \quad
       \Lambda' = \Lambda + 2p^5 \delta_1,
   \end{align}
\end{enumerate}
The coefficient $c(k,n)$ is given by
\begin{equation}
   c(k,n) = \frac{(-1)^{n}}{4^{k}k!n!}\frac{\Gamma(2\Lambda_{1}^*+p-\frac{3}{2})}{\Gamma(2\Lambda_{1}^*+p-\frac{3}{2}-k-n)} \frac{\Gamma(2\Lambda_{2}^*-p+q-\frac{3}{2}+2k+n)}{\Gamma(2\Lambda_{2}^*-p+q-\frac{3}{2})}
\end{equation}
where $ \Lambda_k^*,\ p $ and $q$ are given as follows:
\begin{itemize}
     \item type (iii): $ \Lambda_1^* = \frac{5}{4}-\frac{1}{2}(p^3-q^3),\  \Lambda_2^* = \frac{3}{4}-\frac{q^3}{2}, \ p=p^3,\ q=q^3$
     \item type (iv): $ \Lambda_1^* = \Lambda_1,\  \Lambda_2^* = \frac{3}{4}-\frac{p^4}{2}, \ p=q=p^4$
     \item type (iv): $ \Lambda_1^* = \frac{5}{4}-\frac{p^5}{2},\  \Lambda_2^* = \Lambda_2, \ p = 2p^5,\ q = p^5 $
\end{itemize}
We note that the weight $\Lambda'$ has a unified expression:
\begin{equation}
  \Lambda' = \Lambda + p\delta_1 + (2q-p) \delta_2
\end{equation}
where $p$ and $ q$ for type (i) (ii) are given by
\begin{itemize}
  \item type (i): $ p=p^1,\ q = 0$
  \item type (ii): $ p=0, \ q = p^2$
\end{itemize}

We prove the result in the rest of this section.

\subsection{General facts on the used singular vectors} \label{SEC:Generality}

To perform a complete search of singular vectors,
we need to specify the vectors in the weight space decomposition $ V^{\Lambda}_{p_1,p_2} $, cf. \eqref{wsd}.
To this end, observe first that
\begin{equation}
  p_1 + p_2 \ \text{and} \ K_1 + K_2 \ \text{have the same parity}  \label{p1p2Parity}
\end{equation}
This is due to the relation $ p_1+p_2 = K_1 + K_2 + 2n_1. $

The basis $ \ket{\underline{k}, \underline{m}, \underline{n}} $ of $ V^{\Lambda}$ is labelled by six non-negative integers.
However, if we fix the value of the pair  $(p_1,p_2), $
then two of the six non-negative integers become dependent parameters.
We take $K_1, k_1, k_2 $ and $ n_1$ as independent parameters.
Then we have
\begin{equation}
  m_1 = K_1 - 2k_1, \quad
  m_2 = p_1+p_2-K_1 - 2k_2-2n_1, \quad
  n_2 = p_1 - K_1-n_1.   \label{m1m2n2}
\end{equation}
The possible range of these parameters depends on the value of $(p_1,p_2).$
We will specify the range  later.
At this stage, one may write a singular vector as follows:
\begin{equation}
  \ket{v_s^{p_1 \delta_1 + p_2 \delta_2}}
  = \sum_{K_1,k_1,k_2,n_1} c(K_1, k_1, k_2, n_1) \ket{\underline{k}, \underline{m}, \underline{n}}
  \label{SV1stform}
\end{equation}
where $m_1, m_2, n_2$ are given by \eqref{m1m2n2} and the sum runs all possible values of $K_1, k_1, k_2 $ and $n_1.$
With this expression of $ \ket{v_s^{p_1 \delta_1 + p_2 \delta_2}} $ one may readily show the following:
\\ {\it 
  If $p_1+p_2$ is an odd integer, then there exist no singular vectors.
}\\  

 {\bf Proof:}\
 The condition
\begin{equation}
   a_1^- \ket{v_s^{p_1 \delta_1 + p_2 \delta_2}}
   =  \sum_{K_1,k_1,k_2,n_1} c(K_1, k_1, k_2, n_1) m_1 \ket{\underline{k}, \underline{m}, \underline{n}} = 0
   \label{conda1}
\end{equation}
means that $m_1 = K_1-2k_1$ must vanish to have $ c(K_1, k_1, k_2, n_1) \neq 0.$
Obviously, $ m_1 \neq 0 $ if $ K_1 $ is odd. Thus $K_1$ must be even.

The condition
\begin{equation}
   a_2^- \ket{v_s^{p_1 \delta_1 + p_2 \delta_2}}
   =  \sum_{K_1,k_1,k_2,n_1} c(K_1, k_1, k_2, n_1) m_2 \ket{\underline{k}, \underline{m}, \underline{n}} = 0
   \label{conda2}
\end{equation}
requires $ m_2 = p_1+p_2-K_1 - 2k_2-2n_1 = 0$, but $m_2$ never vanish if $ p_1+p_2$ is odd.
\bsqq

Even when $ p_1+p_2$ is even, the conditions \eqref{conda1} and \eqref{conda2} require $ m_1 = m_2 = 0$ which gives
\begin{equation}
  K_1 = 2k, \qquad k_2 = \rho -k -n, \qquad
  \rho = \frac{1}{2} (p_1+p_2)
\end{equation}
where we set $ k:=k_1, n:= n_1.$
Therefore, the singular vector has the form of
\begin{align}
  \ket{v_s^{p_1 \delta_1 + p_2 \delta_2}}
  = \sum_{k,n} c(k,n) \ket{k,\rho-k-n, n,p_1-2k-n},
  \label{SV2ndform}
\end{align}
where we introduced the new notation for the basis vector:
\begin{align}
  \ket{k,\ell, n,m} :=  (\hat{b}_1^+)^k (\hat{b}_2^+)^{\ell} (\hat{c}^+)^n (d^+)^{m} \ket{0}.
    \nn
\end{align}
By definition of $p_1$ and \eqref{SV2ndform} we see that
\begin{equation}
   0 \leq k \leq \Gau{\frac{p_1}{2}}, \quad
   0 \leq n \leq p_1, \quad
   0 \leq \rho -k -n, \quad
   0 \leq p_1 - 2k-n. \label{paramrange}
\end{equation}
These relations determine the possible range of summation over $k, n$ in \eqref{SV2ndform}.

\bigskip
In the sequel we show the followings:
\\ {\it 
 If $p_1+p_2$ is even, then there exists only one singular vectors for particular values of the lowest weights.
}\\  

First of all, we consider the two special cases for which we have $ k = n = 0,$ namely, there is no summation in \eqref{SV2ndform}:
\begin{enumerate}
  \renewcommand{\labelenumi}{\roman{enumi})}
  \item $ \rho = 0 \ \Leftrightarrow \ p_1 + p_2 = 0 $ \; i.e., \;
     $ \ket{v_s^{p_1 (\delta_1 - \delta_2)}} = \ket{0,0,0,p_1} = (d^+)^{p_1} \ket{0}.$
  \item $ p_1 = 0 $  \; i.e., \;
  $ \ket{v_s^{p_2 \delta_2}} = \ket{0,\frac{p_2}{2},0,0} = (\hat{b}_2^+)^{p_2/2} \ket{0}.$
\end{enumerate}
Then the following is immediate:
\\  
\begin{enumerate}
  \renewcommand{\labelenumi}{\roman{enumi}}
   \item if $ \rho = 0, $ then there exists only one singular vector for $ \Lambda_1-\Lambda_2 = \frac{1}{2}(1-p_1)$ which is given by
   \begin{equation}
      \ket{v_s^{p_1 (\delta_1 - \delta_2)}} = (d^+)^{p_1} \ket{0}. \label{Case1SV}
   \end{equation}
   \item if $p_1 = 0, $ then there exists only one singular vector for $ \Lambda_2 = \frac{1}{4}(3-p_2)$ which is given by
   \begin{equation}
      \ket{v_s^{p_2 \delta_2}}  = (\hat{b}_2^+)^{p_2/2} \ket{0}. \label{Case2SV}
   \end{equation}
\end{enumerate}

 {\bf Proof:}\
\begin{enumerate}
  \renewcommand{\labelenumi}{\roman{enumi})}
  \item It is immediate to see that
  $ \hat{b}_1^-, \hat{b_2}^- $ and $ \hat{c}^-$ annihilate  $\ket{v_s^{p_1 (\delta_1 - \delta_2)}}.$
  Also the following is readily seen
  \begin{equation}
     d^- \ket{v_s^{p_1 (\delta_1 - \delta_2)}} =
     \frac{p_1}{2} \big(\Lambda_2 - \Lambda_1 - \frac{1}{2}(p_1-1) \big) \ket{0,0,0,p_1-1}.
  \end{equation}
  \item One may immediately seen that $\hat{b}_1^-, \hat{c}^- $ and $ d^- $ annihilate $ \ket{v_s^{p_2 \delta_2}}.$
  One also see
  \begin{equation}
       \hat{b}_2^- \ket{v_s^{p_2 \delta_2}}
       = \frac{p_2}{2} \Big(2\Lambda_2+ \frac{p_2}{2} - \frac{3}{2} \Big) \ket{0,\frac{p_2}{2}-1,0,0}.
  \end{equation}
\end{enumerate}
\bsqq

Now we specify the possible range of the parameters $k, n$.
The range, which is determined by \eqref{paramrange}, is classified into three patterns:
\begin{enumerate}
  \renewcommand{\labelenumi}{\roman{enumi})}
  \setcounter{enumi}{2}
   \item $ 0 < p_1 \leq \rho \ \Rightarrow \ 0 \leq k \leq \Gau{\frac{p_1}{2}},\ 0 \leq n \leq p_1-2k $
   \item $ 0 < \rho < p_1$
    \begin{enumerate}
      \item $ p_2 \leq  0 \ \Rightarrow \ 0 \leq k \leq \rho, \ 0 \leq n \leq \rho-k$
      \item $ p_2 > 0 \ \Rightarrow \ 0 \leq k \leq \bar{\rho}, \ 0 \leq n \leq \rho-k $ and
      $ \bar{\rho}+1 \leq k \leq \Gau{\frac{p_1}{2}}, \ 0 \leq n \leq p_1 -2k $
    \end{enumerate}
\end{enumerate}
where $\bar{\rho} := \frac{1}{2}(p_1-p_2). $
This classification is easily seen by considering the intersection of two lines $ \ell_1: n = -k+\rho $ and $\ell_2: n=-2k+p_1$ in the $kn$-plane as depicted below:

\begin{center}
\begin{tikzpicture}
%
%
  \fill[black!20!white] (0,0)--(0,2)--(1,0)--cycle;
  \draw[thick,-latex] (-0.3,0)--(3,0); \node[right] at (3,0) {$k$};
  \draw[thick,-latex] (0,-0.3)--(0,3); \node[above] at (0,3) {$n$};
  \draw[domain=-0.3:2.8] plot (\x,{-\x+2.3}); \node at (3,-0.8) {$\ell_1$};
     \node[right] at (0,2.3) {$\rho$}; \node[right] at (2.3,0.2) {$\rho$};
  \draw[domain=-0.3:1.3] plot (\x,{-2*\x+2}); \node at (1.5,-0.8) {$\ell_2$};
     \node[left] at (0,2) {$p_1$}; \node[below] at (0.8,0) {$\frac{p_1}{2} $};

  \node at (1.5,-1.8) {iii) $p_1 \leq \rho$};
%
%
 \begin{scope}[xshift=5.5cm]
  \fill[black!20!white] (0,0)--(0,1)--(1,0)--cycle;
  \draw[thick,-latex] (-0.3,0)--(3,0); \node[right] at (3,0) {$k$};
  \draw[thick,-latex] (0,-0.3)--(0,3); \node[above] at (0,3) {$n$};
  \draw[domain=-0.3:1.4] plot (\x,{-\x+1}); \node at (1.3,-0.8) {$\ell_1$};
     \node[left] at (0,1) {$\rho$}; \node[below] at (0.8,0) {$\rho$};
  \draw[domain=-0.1:1.5] plot (\x,{-2*\x+2.6}); \node[right] at (0.8,1.5) {$\ell_2$};
     \node[left] at (0,2.5) {$p_1$}; \node[above] at (1.4,0) {$\frac{p_1}{2} $};

  \node at (1.5,-1.8) {iv) (a) $ \rho < p_1$ and $ p_2 \leq 0$};

 \end{scope}

%
%

 \begin{scope}[xshift=11cm]
  \fill[black!20!white] (0,0)--(0,1.8)--(0.8,1)--(1.3,0)--cycle;
  \draw[thick,-latex] (-0.3,0)--(3,0); \node[right] at (3,0) {$k$};
  \draw[thick,-latex] (0,-0.3)--(0,3); \node[above] at (0,3) {$n$};
  \draw[domain=-0.3:2.2] plot (\x,{-\x+1.8}); \node at (2.6,-0.7) {$\ell_1$};
     \node[left] at (0,1.8) {$\rho$}; \node[above] at (1.9,0) {$\rho$};
  \draw[domain=-0.1:1.6] plot (\x,{-2*\x+2.6}); \node[right] at (1.5,-0.9) {$\ell_2$};
     \node[left] at (0,2.5) {$p_1$}; \node[below] at (1.2,0) {$\frac{p_1}{2} $};
  \node[right] at (0.8,1.2) {$(\bar{\rho}, p_2)$};
  \node at (1.5,-1.8) {iv) (b) $\rho < p_1$ and $ p_2 > 0$};

 \end{scope}

\end{tikzpicture}
\end{center}

Before investigating the cases iii) and iv), we derive recurrence relations for $c(k,n)$ which are common for all the cases.
The singular vector \eqref{SV2ndform} must be annihilated by $\hat{b}_1^-, \hat{b}_2^-, \hat{c}^-, d^-.$

The condition $ \hat{b}_1^- \ket{v_s^{p_1 \delta_1 + p_2 \delta_2}} = 0$ gives
\begin{align}
  4(k+1) & \Big( 2\Lambda_1+p_1-k-\frac{5}{2} \Big) c(k+1,n-1)
  \nn \\
  & + 2n(p_1-2k-n) \tilde{\Lambda}\, c(k,n) + n(n+1) c(k,n+1) = 0.
  \label{RecRel1}
\end{align}
From the condition $ \hat{b}_2^- \ket{v_s^{p_1 \delta_1 + p_2 \delta_2}} = 0$ we have
\begin{align}
  4 (\rho-k-n) & \Big(2\Lambda_2-\bar{\rho}+k+n-\frac{3}{2} \Big) c(k,n)
  \nn \\
  & + (n+1) (n+2) c(k-1,n+2) + 4 (n+1) c(k,n+1) = 0.
  \label{RecRel2}
\end{align}
The condition $ \hat{c}^- \ket{v_s^{p_1 \delta_1 + p_2 \delta_2}} = 0$ gives
\begin{align}
  (n+1) & \Big( \Lambda_1 + \Lambda_2 + \rho - \frac{n}{2} -\frac{3}{2} \Big) c(k,n+1)
  + (\rho-k-n) (p_1-2k-n) \tilde{\Lambda} \, c(k,n)
  \nn \\
  & + 2 (k+1) c(k+1,n)
  + 2 (k+1) (\rho-k-n) c(k+1,n-1) = 0.
  \label{RecRel3}
\end{align}
From the condition $ d^- \ket{v_s^{p_1 \delta_1 + p_2 \delta_2}} = 0$ we have
\begin{align}
  2 (k+1) c(k+1,n-1) + (n+1) c(k,n+1)
   + (p_1-2k-n) \tilde{\Lambda}\, c(k,n) = 0.
   \label{RecRel4}
\end{align}
The possible range of $k$ and $n$ depends on the cases, so we will specify it later.
In addition to these relations, there exit more recurrence relations stems from the ``boundary values" (i.e., minimum and maximum values) of $k, n$ which also depend on the cases.  They will also be presented when we discuss the each cases.

We  find that the relations \eqref{RecRel1} and \eqref{RecRel4} are solved for arbitrarily $ \Lambda_1, \Lambda_2$ and the solution is given by
\begin{align}
  c(k,n) &=
  \frac{(-1)^n}{4^k k! n!}
  \frac{p_1 !}{(p_1-2k-n)!}
  \frac{\Gamma(2\Lambda_1+p_1-\frac{3}{2})}{\Gamma(2\Lambda_1+p_1-k-n-\frac{3}{2})}
  \nn \\
  & \times
  \frac{\Gamma(2\Lambda_2 - 2\Lambda_1 -p_1 +2k+n+ 1)}{\Gamma(2\Lambda_2 - 2\Lambda_1 -p_1 + 1)}
  \frac{\Gamma(4\Lambda_1+2p_1-2k-n-4)}{\Gamma(4\Lambda_1+2p_1-4)} \, c(0,0).
  \label{ckn}
\end{align}
However, \eqref{ckn} does not solve the relations \eqref{RecRel2} and \eqref{RecRel3} unless $ \Lambda_1, \Lambda_2$ take particular values.
As we will see later,  \eqref{RecRel2} and \eqref{RecRel3} together with other recurrence relations fix the value of the lowest weight.

In the next subsections we are finding more specific results.

%
\subsection{Case iii) \bm{$ 0 < p_1 \leq \rho$}} \label{SEC:Case3}

 {\it 
 There exists only one singular vector of the form
 \begin{equation}
      \ket{v_s^{p_1 \delta_1 + p_2 \delta_2}}
      = \sum_{k=0}^{\Gau{p_1/2}} \sum_{n=0}^{p_1-2k}
      c(k,n)
      \ket{k,\rho-k-n,n,p_1-2k-n}   \label{SVform1}
 \end{equation}
 for the following two cases
 \begin{itemize}
   \item $ 0 < p_1 < \rho$ and $ \Lambda_1 = \frac{5}{4}- \frac{1}{2}\bar{\rho}, \;
   \Lambda_2 = \frac{3}{4} - \frac{1}{2} \rho$
   \begin{equation}
      c(k,n) = \frac{p_1!\, \rho!}{4^k k! n! (p_1-2k-n)! (\rho-k-n)!}. \label{ckn3}
   \end{equation}
   \item $ 0 < p_1 = \rho $ and $ \Lambda_1 + \Lambda_2 = 2- \frac{p_1}{2}$
   \begin{equation}
      c(k,n) = \frac{p_1!}{4^k k! n! (p_1-2k-n)!}
      \frac{\Gamma(2\Lambda_1+p_1-\frac{3}{2})}{\Gamma(2\Lambda_1+p_1-k-n-\frac{3}{2})}.
      \label{ckn4}
   \end{equation}
 \end{itemize}
}

The~{\it Proof}~ of the above statement is straightforward by applying the conditions of annihilation of \eqref{SVform1} by the negative generators
and we omit it.

%
\subsection{Case iv)(a) \bm{$ 0 < \rho < p_1, \ p_2 \leq 0$}}

 {\it 
 There exists only one singular vector of the form
 \begin{equation}
      \ket{v_s^{p_1 \delta_1 + p_2 \delta_2}}
      = \sum_{k=0}^{\rho} \sum_{n=0}^{\rho-k}
      c(k,n)
      \ket{k,\rho-k-n,n,p_1-2k-n}   \label{SVform2}
 \end{equation}
 for the following two cases
 \begin{itemize}
   \item $ p_2 < 0 $ and $ \Lambda_1 = \frac{5}{4}- \frac{1}{2}\bar{\rho}, \;
   \Lambda_2 = \frac{3}{4} - \frac{1}{2} \rho.$
   \\[6pt]
   $c(k,n)$ is given by \eqref{ckn3}.
   \item $ p_2 = 0$ and $ \Lambda_1 = \frac{5}{4} - \frac{1}{4}p_1, \; {}^{\forall} \Lambda_2$
   \begin{equation}
      c(k,n) = \frac{(-1)^n}{4^k k! n!}
      \frac{(\frac{p_1}{2})!}{(\frac{p_1}{2}-k-n)!}
      \frac{\Gamma(2\Lambda_2-\frac{p_1}{2}-\frac{3}{2}+2k+n)}{\Gamma(2\Lambda_2-\frac{p_1}{2}-\frac{3}{2})}. \label{ckn5}
   \end{equation}
 \end{itemize}
}

As in the previous Case (iii) the~{\it Proof}~ is straightforward.

%
\subsection{Case iv)(b) \bm{$ 0 < \rho < p_1, \ p_2 > 0$}}

 {\it 
 There exists only one singular vector of the form
 \begin{equation}
      \ket{v_s^{p_1 \delta_1 + p_2 \delta_2}}
      =
      \left(
       \sum_{k=0}^{\bar{\rho}} \sum_{n=0}^{\rho-k}
       +
       \sum_{k=\bar{\rho}+1}^{\Gau{p_1/2}} \sum_{n=0}^{p_1-2k}
      \right)
      c(k,n)
      \ket{k,\rho-k-n,n,p_1-2k-n}   \label{SVform3}
 \end{equation}
 for $ \Lambda_1 = \frac{5}{4}- \frac{1}{2}\bar{\rho}, \; \Lambda_2 = \frac{3}{4} - \frac{1}{2} \rho$ and $ c(k,n) $ is given by \eqref{ckn3}.
}\\  

The~{\it Proof}~ is the same as the previous cases.

\bigskip
We  successfully carried out the search of  all singular vectors in $V^{\Lambda}$. 
For the sake of simplicity and for the convenience of later computation, we change the notations of parameters specifying the lowest weight.
By this change, all the new parameters $ p^1, p^2, p^3, p^4, p^5 $ and $ q^3 $ take a positive integer.

Case i) \ $ p_1 \ \to \ p^1$ in \eqref{Case1SV}

Case ii) \ $ p_2 \ \to \ 2p^2 $ in \eqref{Case2SV}

\noindent
These two cases corresponds to the type (i) (ii), respectively.
We have the common constraint on $ \Lambda_1, \Lambda_2$ and the common form of $ c(k,n)$ for Case iii) $(0 < p_1 < \rho)$, Case iv) (a) $ (p_2 < 0)$ and Case iv) (b) $(0 < \rho < p_1, p_2 > 0)$. By the change 
\begin{equation}
   p_1 \ \to \ p^3, \qquad \rho \ \to \ q^3
\end{equation}
which leads $ \bar{\rho} = p^3 - q^3, \ p_2 = 2q^3-p^3, $ 
these three cases correspond to the three subcases of the type (iii):
\begin{align}
   & \text{Case iii)} \ 0 < p_1 < \rho \quad \to \quad  0 < p^3 < q^3,
   \nonumber \\
   & \text{Case iv) (a)} \ p_2 < 0 \quad \to \quad 2q^3 < p^3
   \nonumber \\
   & \text{Case iv) (b)} \ 0 < \rho < p_1, 0 < p_2 \quad \to \quad q^3 < p^3 < 2q^3
\end{align}
For Case iii) with $ p_1 = \rho,$ we make the change $  p_1 \ \to \ p^4$ and this corresponds to the type (iv). 
Finally, for Case iv) (a) with  $ p_2 = 0 $
we make the change $ p_1 \ \to \ 2p^5 $ and this gives the type (v). 

\bigskip
In this way, we have completed the proof of the result in \S \ref{gendef}.

%
%
%

\section{Elementary reducibilities}
\label{SEC:EmbPat}
\setcounter{equation}{0}

In the precious section we have established that there are ~{\it five} types of singular vectors. Here we shall study the cases when some Verma modules have SV of various types.

\bigskip

\noindent
\textbf{A12345}

We start with Verma modules which are reducible under ~{\it all five} types. We denote these as  A12345. We shall see that these are related
to two more types: A234, which have SVs of types (ii), (iii), (iv); further: A2345 which have SVs of types (ii), (iii), (iv), (v).

We take the LW of type (iii) as  representative {and denote it $\Lambda^0$}:
 \eqna{A3}
  { \Lambda^{0}} &=& \Big(\frac{5}{4}-\frac{1}{2}(p^3-q^3), \frac{3}{4}-\frac{q^3}{2} \Big),
  \\
  && p^1 = p^3 -2q^3, \quad p^2 = q^3, \quad p^4 = p^3, \quad p^5 = p^3-q^3,
  \eena
where by $p^j$, $j=1,2,4,5$ we denote the parameters of the corresponding types of reducibility and we give how these are related with the parameters of type (iii) $p^3,q^3$.

 We establish the following:
\\ {\it
  The VM with the LW  {$ \Lambda^0$} $ = \Big(\frac{5}{4}-\frac{1}{2}(p^3-q^3), \frac{3}{4}-\frac{q^3}{2} \Big) $  has the SVs $\L'$ of the following weights:}
   \begin{itemize}
    \item $ p^3 > 2q^3 \quad \Rightarrow \quad $ \textbf{A12345}
     \eqn{emb12345}
       {\Lambda^{1}} = \Big( \frac{5}{4}- \frac{q^3}{2}, \frac{3}{4}-\hf (p^3-q^3) \Big).
     \end{equation}
     \item $  q^3 > p^3   \quad \Rightarrow \quad $ \textbf{A234}
        \eqna{emb234}
          {\Lambda^{2}} &= \Big( \frac{5}{4}-\frac{1}{2}(p^3-q^3), \frac{3}{4}+\frac{q^3}{2} \Big),
          \\
          {\Lambda^{3}} &= \Big( \frac{5}{4}+\frac{q^3}{2}, \frac{3}{4}-\hf(p^3-q^3) \Big),
          \\
          {\Lambda^{4}} &= \Big( \frac{5}{4}+\frac{q^3}{2}, \frac{3}{4}+ \hf (p^3-q^3) \Big).
       \eena
     \item $ 2q^3 > p^3 > q^3 \quad \Rightarrow \quad $ \textbf{A2345}
     \eqn{emb2345}
       {\Lambda^{5}} = \Big( \frac{5}{4}+ \hf (p^3-q^3), \frac{3}{4}-\frac{q^3}{2} \Big).
     \end{equation}
       \end{itemize}

Since  {$ \Lambda^{0} $} is the LW of Case (iii), this result covers all the LW coincide with $ \Lambda^{0},$ e.g., \textbf{A13}, \textbf{A123}.

\bigskip

\bigskip\noindent
\textbf{A1245}

Next we consider VMs which are reducible under conditions (i),(ii),(iv),(v) but not under case (iii).
Actually, we start with the coincidence of the LWs for cases (ii) and (v) which gives the relations:
 \eqnn{a25}
  \Lambda_1 &=& \frac{5}{4} - \frac{p^5}{2},
  \\
  \Lambda_1 + \hf (p^1-1) &=& \frac{3}{4}-\frac{p^2}{2} = -\Lambda_1 + 2 - \frac{p^4}{2} = \Lambda_2.
 \nn\eea
Take $ p^2, p^5$ as independent parameters (i.e., the LW of \textbf{A25}), then
\begin{equation}
   {\Lambda^{0}_{1245}= \Lambda^0_{25}} = \Big( \frac{5}{4}-\frac{p^5}{2}, \frac{3}{4}-\frac{p^2}{2} \Big),
   \qquad
   p^1 = p^5-p^2, \quad p^4 = p^5+p^2.
   \label{L1245}
\end{equation}

We have still to exclude case (iii) so we divide  into two cases:

\noindent
\textbf{a) $ p^5 \neq p^2$}

For this case, {$  \Lambda^0_{1245} $} takes a same value as {$ \Lambda^{0}. $}
This is seen by setting
\begin{equation}
   p^5 = p^3-q^3 > 0, \quad p^2 = q^3 \label{p52topq3}
\end{equation}
which brings {$ \Lambda^0_{1245}$} to the same form as {$ \Lambda^{0}. $}
The conditions $ p^1 > 0 $ and $ p^4 > 0 $ becomes
\begin{equation}
    p^1 = p^3 -2q^3 > 0, \qquad  p^4 = p^3 > 0
\end{equation}
which shows that \textbf{A1245} corresponds to $ p^3 > 2q^3$ case of {$ \Lambda^{0}, $} i.e. \textbf{A12345}.
Thus, \textbf{A1245} for $ p^5 \neq p^2$ is a subcase of \textbf{A12345}.

\medskip\noindent
\textbf{b) $ p^5 = p^2$}

For this case, {$  \Lambda^0_{1245} $} never coincide with {$ \Lambda^{0}. $}
Because the replacement \eqref{p52topq3} gives $ p^3 = 2q^3$ which is not allowed for $ \Lambda^{0}. $
Also, $p^1$ is not admissible since $ p^1 = 0$ for $p^5=p^2$ (Cf. $p^4 =2p^2$ is admissible).
Namely, the case \textbf{A1} is decoupled from \textbf{A1245}, since if we apply the \textbf{A1} conditions to the signature of \textbf{A245} we shall obtain \textbf{A1} with  $p^1=0$ which is not allowed.

Thus, we conclude:\\
 {\it
  The VM with the LW
\eqn{A245f} \Lambda^0_{245} = \Big( \frac{5}{4}-\frac{p^2}{2}, \frac{3}{4}-\frac{p^2}{2} \Big)\ee
 has the SVs  of the weights:}
   \eqna{a245}
     \Lambda^1_{245} &= \Big( \frac{5}{4}-\frac{p^2}{2}, \frac{3}{4}+\frac{p^2}{2} \Big),
      \\
      \Lambda^2_{245}&= \Big( \frac{5}{4}+\frac{p^2}{2}, \frac{3}{4}-\frac{p^2}{2} \Big),
      \\
      \Lambda^3_{245} &= \Big( \frac{5}{4}+\frac{p^2}{2}, \frac{3}{4}+\frac{p^2}{2} \Big)
   \eena
where the type of SVs are (ii), (v), (iv), respectively.
 Since  {$ \Lambda^0_{1245} = \Lambda^0_{25}, $}  \textbf{A25}, \textbf{A125} are the  subcases of \textbf{A12345} or  \textbf{A245}.

\bigskip

\bigskip\noindent
\textbf{A124}
Next we try to obtain VMs which are reducible under conditions (i),(ii),(iv) but not under cases (iii),(v).
Actually, we start with the coincidence of the LWs for cases (i) and (ii) which gives the relations:
 \begin{equation}
  \Lambda_1 + \hf (p^1-1) = \frac{3}{4}-\frac{p^2}{2} = -\Lambda_1 + 2 - \frac{p^4}{2}.
\end{equation}
Take $p^1, p^2$ as independent parameters, then
\begin{equation}
  \Lambda^0_{124} = \Lambda^0_{12} =  \Big( \frac{5}{4}-\hf (p^1+p^2), \frac{3}{4}-\frac{p^2}{2} \Big),
  \quad
  p^4 = p^1 + 2p^2.
\end{equation}
Obviously, this is a subcase of \textbf{A1245}.
Also \textbf{A12}  is a subcase of \textbf{A1245}.

\bigskip

\bigskip\noindent
\textbf{A145}
Next we try to obtain VMs which are reducible under conditions (i),(iv),(v) but not under cases (ii),(iii).
Actually, we start with the coincidence of the LWs for cases (iv) and (v) which gives the relations:
 \eqna{a145}
   \Lambda_1 &=& \frac{5}{4}-\frac{p^5}{2},
   \\
   \Lambda_1 + \hf (p^1-1)  &=&  -\Lambda_1 + 2 - \frac{p^4}{2} = \Lambda_2.
\eena
Take $p^4, p^5$ as independent parameters, then
\begin{equation}
  \Lambda^0_{145} = \Lambda^0_{45}= \Big(  \frac{5}{4}-\frac{p^5}{2}, \frac{3}{4}- \hf (p^4-p^5) \Big),
  \quad
  p^1 = 2p^5-p^4.  \label{A145}
\end{equation}
Compare this with $\Lambda^0_{1245},$ then one may see $ \Lambda^0_{145} = \Lambda^0_{1245}$ if $ p^4-p^5 > 0.$ Because by setting $ p^2 = p^4-p^5 > 0$
 we have:
 $$  \Lambda^0_{145} = \Big(  \frac{5}{4}-\frac{p^5}{2}, \frac{3}{4}- \frac{p^2}{2} \Big),
  \quad
  p^1 = p^5-p^2$$
 which is identical to \eqref{L1245}.

However, if $ p^4 \leq p^5$ then \eqref{A145} is never identical to $ \Lambda^0_{1245}.$
Note that $ p^1 > 0 $ for $ p^4 \leq p^5$ so that $ p^1$ takes an admissible value. Thus, we shall use:
\eqn{A145f} \Lambda^0_{145} =  \Big(  \frac{5}{4}-\frac{p^5}{2}, \frac{3}{4}- \hf (p^4-p^5) \Big), \quad p^4 \leq p^5 \ee
We find:
\\ {\it %
    The VM with the LW \eqref{A145f}   has the SVs  with  weights:}
     \eqna{L145}
       \Lambda^1_{145} &= \Big( \frac{5}{4}+\hf(p^5-p^4), \frac{3}{4}-\frac{p^5}{2} \Big),  
       \\
       \Lambda^2_{145}   &= \Big( \frac{5}{4}-\hf(p^5-p^4), \frac{3}{4}+\frac{p^5}{2} \Big), 
       \\
       \Lambda^5_{145}  &= \Big( \frac{5}{4}+\frac{p^5}{2}, \frac{3}{4}- \hf (p^4-p^5) \Big) 
    \eena
where the type of SVs are (i), (iv), (v), respectively.
The somewhat unnatural enumeration of $\L_{145}$'s is for later consitency with Section 6.
    Since $ \Lambda^0_{45} = \Lambda^0_{145}$, \textbf{A45} with $ p^4  \leq p^5$ always implies \textbf{A145}.

\bigskip

\bigskip\noindent
\textbf{A14}
Next we try to obtain VMs which are reducible under conditions (i),(iv)  but not under cases (ii),(iii),(v).
After some analysis we find that
   the VM with the LW
  \eqn{A14f}
  \Lambda^0_{14} = \Big( 1-\frac{q}{2}, \hf(1+p^1-q) \Big), \quad  p^1, q \in \mathbb{Z}_{\geq 1},  \quad  p^1 < 1+2q  \ee
  has the SVs   of   weights:
 \eqna{L14}
    \Lambda^1_{14}  &=& \Big( 1 + \hf (p^1-q),\hf (1-q) \Big),
     \\
    \Lambda^4_{14}  &=& \Big( \hf (3-p^1+q), 1+ \frac{q}{2} \Big).
  \eena

\bigskip

\bigskip\noindent
\textbf{A15}

Next we try to obtain VMs which are reducible under conditions (i),(v)  but not under cases (ii),(iii),(iv).
 The coincidence of the LWs gives the relations:
\begin{equation}
    \Lambda_1 = \frac{5}{4}-\frac{p^5}{2},
    \qquad
    \Lambda_1 + \hf (p^1-1) = \Lambda_2.
\end{equation}
Thus
\begin{equation}
   \Lambda^0_{15} = \Big( \frac{5}{4}-\frac{p^5}{2}, \frac{3}{4}-\hf (p^5-p^1) \Big).
\end{equation}
Comparing this with  $\Lambda^0_{145}$ we see that if   $ 2p^5 > p^1, $ then $ \Lambda^0_{15} = \Lambda^0_{145}. $

On the other hand, if $ 2p^5 \leq p^1, $  then $ \Lambda^0_{15} \neq \Lambda^0_{145} $ and $ \Lambda^0_{15}$ is not identical to any LW considered so far.
 Thus, we  shall use:
 \eqn{A15f}
 \Lambda^0_{15}= \Big( \frac{5}{4}-\frac{p^5}{2}, \frac{3}{4}-\hf (p^5-p^1) \Big), \quad  2p^5 \leq p^1 \ee
We find: \\
  {\it
  The VM with  LW \eqref{A15f}  has   SVs  of   weights:}
  \eqna{L15}
  \Lambda^1_{15}   &=& \Big( \frac{5}{4} + \hf (p^1-p^5), \frac{3}{4}-\frac{p^5}{2} \Big), 
     \\
  \Lambda^2_{15}    &=& \Big(  \frac{5}{4}+\frac{p^5}{2}, \frac{3}{4}-\hf (p^5-p^1) \Big).
  \eena

\bigskip

\bigskip\noindent
\textbf{A24}

Next we look for  VMs which are reducible under conditions (ii),(iv)  but not under cases (i),(iii),(v).
The coincidence of the LWs gives the relation:
\begin{equation}
  \frac{3}{4}-\frac{p^2}{2} = -\Lambda_1 + 2 - \frac{p^4}{2}.
\end{equation}
Thus
\begin{equation}
  \Lambda^0_{24} = \Big( \frac{5}{4}-\hf (p^4-p^2), \frac{3}{4}-\frac{p^2}{2} \Big).
\end{equation}
Comparing this with $ \Lambda^{0}, $ we  see that $ \Lambda^0_{24} = \Lambda^{0} $ if $ p^4 \neq p^2 $ or $ p^4 \neq 2p^2$ by the correspondence
\begin{equation}
  p^4 \ \leftrightarrow \ p^3, \qquad p^2 \ \leftrightarrow \ q^3
\end{equation}
Thus if $ p^4 \neq p^2 $ or $ p^4 \neq 2p^2,$ then \textbf{A24} is a subcase of \textbf{A12345}.

 On the other hand for $ p^4 = 2p^2 $ we have:
\begin{equation}
   \Lambda^0_{24} = \Big( \frac{5}{4}-\frac{p^2}{2}, \frac{3}{4}-\frac{p^2}{2} \Big)
\end{equation}
which is identical to the LW of \textbf{A245}.

Finally, if $ p^4 = p^2, $ then the  LW:
\eqn{A24f} \Lambda^0_{24} =  \Big( \frac{5}{4}, \frac{3}{4}-\frac{p^2}{2} \Big)\ee
 is not identical to any LW considered so far. Thus,  we have:

 {\it
 The VM with the LW \eqref{A24f}   has the SVs  of  weights:}
  \eqna{L24}
   \Lambda^1_{24}   &=& \Big( \frac{5}{4}, \frac{3}{4}+\frac{p^2}{2} \Big), 
     \\
   \Lambda^2_{24}  &=& \Big( \frac{5}{4}+\frac{p^2}{2}, \frac{3}{4} \Big). 
 \eena

\bigskip

\bigskip\noindent
\textbf{A1}

Next we try to obtain VMs which are reducible only under condition (i).
We need to exclude all the coincidence with the type (i) LW. After some tedious analysis we obtain that it is enough to
exclude coincidence with case A14. The latter is given by:
\begin{equation}
   \Big( \Lambda_1, \Lambda_1+ \hf (p^1-1) \Big), 
   \quad   \Lambda_1 = \frac{1}{4}(5-p^1-p^4).
      \label{A1}
\end{equation}
 Thus, we have:

  {\it 
  The VM of   LW
   \eqna{L1}
     \Lambda^0_1 = \Big( \Lambda_1, \Lambda_1+ \hf (p^1-1) \Big),
     \\
     -4\Lambda_1 +5 -p^1 \notin \bbn 
    \eena
    has only the type (i) SV of   weight:}
  \eqn{L11} \Lambda^1_1  = \Big( \Lambda_1 + \frac{p^1}{2},  \Lambda_1 -\hf \Big).  \ee
  (The condition in \rf{L1}b excludes the SVs of type (iv) - the excluded positive integer would be identified with $p^4$.)

\bigskip

\bigskip\noindent
\textbf{A2}

Next we try to obtain VMs which are reducible only under conditions (ii).
We need to exclude all the coincidence with the type (i) LW. After some tedious analysis we obtain that it is enough to
exclude coincidence with case A24. The latter is given by:
\eqn{A2f}
  \Big( \Lambda_1, \frac{3}{4}-\frac{p^2}{2} \Big).
\end{equation}
Thus, we have:

  {\it
  The VM of   LW
  \eqna{L2}
     \Lambda^0_2  = \Big( \Lambda_1, \frac{3}{4}-\frac{p^2}{2} \Big),
     \\
     -2\Lambda_1 + \frac{5}{2} + p^2 \notin \bbn
  \eena
  has only the type (ii) SV of   weight:}
  \eqn{L21} \Lambda^1_2  = \Big( \Lambda_1, \frac{3}{4}+\frac{p^2}{2} \Big).   \ee
  (The condition in \rf{L2}b excludes the SVs of type (iv).)

\bigskip\noindent
\textbf{A3}

  We have seen that the LW $ \Lambda^{0} $ has at least \textbf{A234} (see the part \textbf{A12345}), so there exists no VM having only type (iii) SV.

\bigskip\noindent
\textbf{A4}

Next we try to obtain VMs which are reducible only under conditions (iv).
We need to exclude all the coincidence with the type (iv) LW. After some tedious analysis we obtain that it is enough to
exclude coincidence with case A14. The latter is given by:

\eqn{A4f}
   \Big(\Lambda_1, -\Lambda_1 + 2 -\frac{p^4}{2} \Big), 
   \qquad  \Lambda_1 =  \frac{1}{4}(5-p^1-p^4)
\end{equation}
Thus, we have

 {\it
  The VM of the LW
  \eqna{L4}
     \Lambda^0_4  = \Big( \Lambda_1, -\Lambda_1 + 2 -\frac{p^4}{2}  \Big),
     \\
     -4\Lambda_1 +5 -p^4 \notin\bbn
  \eena
  has only the type (iv) SV of  weight:}
  \eqn{L41} \Lambda^1_4  = \Big( \Lambda_1+ \frac{p^4}{2}, -\Lambda_1 + 2 \Big).   \ee
  (The condition in \rf{L4}b excludes the SVs of type (i).)

\bigskip\noindent
\textbf{A5}

Next we try to obtain VMs which are reducible only under conditions (v).
We need to exclude all the coincidence with the type (v) LW. After some tedious analysis we obtain that it is enough to
exclude coincidence with cases A15 and A45 since they   contain all other cases.
With an abuse of notation we represent them as:
\eqn{A5f}
   \Lambda^0_{15} \cup \Lambda^0_{45}= \Big( \frac{5}{4}-\frac{p^5}{2}, \frac{3}{4}-\frac{r}{2} \Big), \quad
   r \in \mathbb{Z}
\end{equation}
The relation which gives the coincidence with the type (v) LW
\begin{equation}
   \Big( \frac{5}{4}-\frac{p^5}{2}, \Lambda_2 \Big)
   = \Lambda^0_{15} \cup \Lambda^0_{45}
\end{equation}
is reduced to $ \Lambda_2 = \frac{3}{4}-\frac{r}{2}.$\\
Thus, we have:

  {\it
  The VM of   LW
   \eqn{L5}
     \Lambda^0_5  = \Big(  \frac{5}{4}-\frac{p^5}{2}, \Lambda_2  \Big),
     \quad
     \Lambda_2 \neq \frac{3}{4}-\frac{r}{2}, \ r \in \mathbb{Z}
  \ee
  has only the type (v) SV of   weight:}
  \eqn{L51} \Lambda^1_5  = \Big(   \frac{5}{4}+\frac{p^5}{2}, \Lambda_2   \Big).   \ee

\bigskip\noindent

In {\bf conclusion}, we have the following types of elementary embeddings: A12345, A234, A2345, A245, A145, A14, A15, A24, A1, A2, A4, A5, which are distinguished
from one another by conditions given above.  In the next section we shall find their complete  embedding pictures.


\section{Complete embedding pictures}
\setcounter{equation}{0}

 We shall find the complete  embedding pictures of the reducible Verma modules. These are obtained using the
elementary picture from \S \ref{gendef}  not only to the initial modules but also to their invariant submodules.

\subsection{Case A12345}

As we have seen the VM which is reducible under case (iii) may have all five reducibilities (some under some conditions).
The initial embedding picture is as follows:

\bigskip\bigskip

%
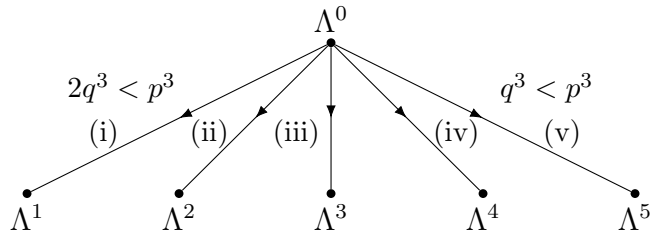
\begin{center}
  \begin{tikzpicture}
    \node[circle,fill=black,scale=0.3] (L3) at (0,10) {};
    \node[circle,fill=black,scale=0.3] (Ld1) at (-4,8) {};
    \node[circle,fill=black,scale=0.3] (Ld2) at (-2,8) {};
    \node[circle,fill=black,scale=0.3] (Ld3) at (0,8) {};
    \node[circle,fill=black,scale=0.3] (Ld4) at (2,8) {};
    \node[circle,fill=black,scale=0.3] (Ld5) at (4,8) {};
    \draw[decoration={markings, mark=at position 0.5 with {\arrow{Latex}}},postaction={decorate}] (L3) to (Ld1);
    \node[left,xshift=-7] at ($(L3)!0.6!(Ld1)$) {\small (i)};
    \draw[decoration={markings, mark=at position 0.5 with {\arrow{Latex}}},postaction={decorate}] (L3) to (Ld2);
    \node[left] at ($(L3)!0.6!(Ld2)$) {\small (ii)};
    \draw[decoration={markings, mark=at position 0.5 with {\arrow{Latex}}},postaction={decorate}] (L3) to (Ld3);
    \node[left] at ($(L3)!0.6!(Ld3)$) {\small (iii)};
    \draw[decoration={markings, mark=at position 0.5 with {\arrow{Latex}}},postaction={decorate}] (L3) to (Ld4);
    \node[right] at ($(L3)!0.6!(Ld4)$) {\small (iv)};
    \draw[decoration={markings, mark=at position 0.5 with {\arrow{Latex}}},postaction={decorate}] (L3) to (Ld5);
    \node[right,xshift=7] at ($(L3)!0.6!(Ld5)$) {\small (v)};
    \node[above] at (L3) {$\Lambda^0$}; 
    \node[below] at (Ld1) {$\Lambda^1$}; 
    \node[below] at (Ld2) {$\Lambda^2$}; 
    \node[below] at (Ld3) {$\Lambda^3$}; 
    \node[below] at (Ld4) {$\Lambda^4$}; 
    \node[below] at (Ld5) {$\Lambda^5$}; 
    \node[left,xshift=-20] at ($(L3)!0.3!(Ld1)$) {\small $ 2q^3 < p^3$};
    \node[right,xshift=25] at ($(L3)!0.3!(Ld5)$) {\small $ q^3 < p^3$};
  \end{tikzpicture}
   \end{center}

where the LWs $ \L^{a}, \ a = 0, 1, \dots, 5 $ are given in  \eqref{A3} - \eqref{emb2345}.
%
%
%
We now obtain the embeddings of each of these five cases.

\subsubsection{Invariant modules in $ V^{\Lambda^1} $}

 After some analysis we find that $ V^{\Lambda^1}$ has the type (ii), (iii), (iv) and (v) SVs for all possible values of $ p^3, q^3.$ Their weights are given by
~$\L^{12}$,  ~$\L^{5}$,  ~$\widetilde{\L}^{0}$, ~$\L^{3}$, respectively, where:
\eqna{Lamvv} 
   \Lambda^{12} &=&
   \Big( \frac{5}{4}-\frac{q^3}{2}, \frac{3}{4}+\hf(p^3-q^3) \Big),\\
      \widetilde{\L}^{0} &=&                
   \Big( \frac{5}{4}+\hf(p^3-q^3), \frac{3}{4}+\frac{q^3}{2} \Big),
   \eena
 and  use notation $\widetilde{\L}^{0}$ for the LW which differs from $\L^0$ by the signs of the additions to $5/4,3/4$. 

Next, we seek SVs in $V^{\L^{12}}. $ 
We find that it has  types (i) and (v) SVs for all possible values of $p^3, q^3$ and the weights of SVs are given by respectively by  ~$\L^5$, $\L^4$.

Finally, we seek SVs in $ V^{\widetilde{\L}^{0}}. $ 
We find that it has  no SVs. 

The results are summarized in the diagram:


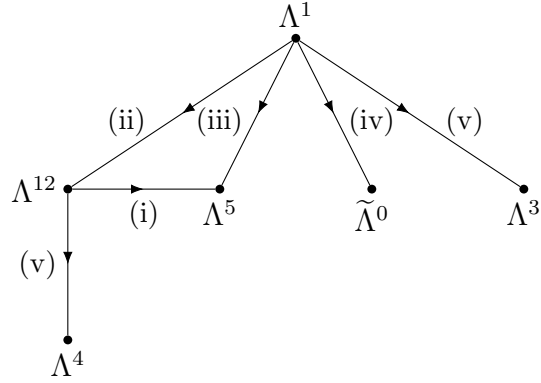
\begin{center}
  \begin{tikzpicture}
    \node[circle,fill=black,scale=0.3] (L1) at (0,10) {};
    \node[circle,fill=black,scale=0.3] (Ld2) at (-3,8) {};
    \node[circle,fill=black,scale=0.3] (Ld3) at (-1,8) {};
    \node[circle,fill=black,scale=0.3] (Ld4) at (1,8) {};
    \node[circle,fill=black,scale=0.3] (Ld5) at (3,8) {};
    \node[circle,fill=black,scale=0.3] (Ldd2) at (-3,6) {};
    \draw[decoration={markings, mark=at position 0.5 with {\arrow{Latex}}},postaction={decorate}] (L1) to (Ld2);
    \node[left,xshift=-5] at ($(L1)!0.55!(Ld2)$) {\small (ii)};    %
    \draw[decoration={markings, mark=at position 0.5 with {\arrow{Latex}}},postaction={decorate}] (L1) to (Ld3);
    \node[left] at ($(L1)!0.55!(Ld3)$) {\small (iii)};    %
    \draw[decoration={markings, mark=at position 0.5 with {\arrow{Latex}}},postaction={decorate}] (L1) to (Ld4);
    \node[right] at ($(L1)!0.55!(Ld4)$) {\small (iv)};    %
    \draw[decoration={markings, mark=at position 0.5 with {\arrow{Latex}}},postaction={decorate}] (L1) to (Ld5);
    \node[right,xshift=5] at ($(L1)!0.55!(Ld5)$) {\small (v)};    %
    \draw[decoration={markings, mark=at position 0.5 with {\arrow{Latex}}},postaction={decorate}] (Ld2) to (Ldd2);
    \node[left] at ($(Ld2)!0.5!(Ldd2)$) {\small (v)};    %
    \draw[decoration={markings, mark=at position 0.5 with {\arrow{Latex}}},postaction={decorate}] (Ld2) to (Ld3);
    \node[below] at ($(Ld2)!0.5!(Ld3)$) {\small (i)};    %
    \node[above] at (L1) {$\Lambda^1$}; 
    \node[left] at (Ld2) {$\Lambda^{12}$}; 
    \node[below] at (Ld3) {$\Lambda^5$}; 
    \node[below] at (Ld4) {$\widetilde{\L}^{0}$}; 
    \node[below] at (Ld5) {$\Lambda^3$}; 
    \node[below] at (Ldd2) {$\Lambda^4$}; 
  \end{tikzpicture}
\end{center}


\subsubsection{Invariant modules in $V^{\L^2} $}

\medskip\noindent After some analysis we find that  $V^{\L^2}$ has the following SVs:
\begin{itemize}
   \item $ 2q^3 <  p^3  \ \Rightarrow \ $ type (i) (iv) (v) SVs of weights ~$\L^3$, $\L^{12}$, $\widetilde{\L}^{0}$;
   \item $ q^3 < p^3 < 2q^3 \ \Rightarrow \ $ type (i) (v) SVs of weights $ \L^3  $ and $ \widetilde{\L}^{0}, $ respectively.
    \item $ p^3 < q^3 \ \Rightarrow \ $type (i) SV of weights $ \L^{3}. $ 
\end{itemize}

Further embedding in $ V^{\L^{12}} $ when $ 2q^3 < p^3$ are given
by type (i) and (v) SVs of weight $ \L^5 $ and $ \L^4 , $ respectively.


Further embedding in $ V^{\widetilde{\L}^{0}} $ where $ q^3 < p^3$ is given  by   type (i) SV of the weight $\L^4$.


There are no further embeddings  and thus,
the embedding picture of $ V^{\Lambda^2} $ 
is given as follows:

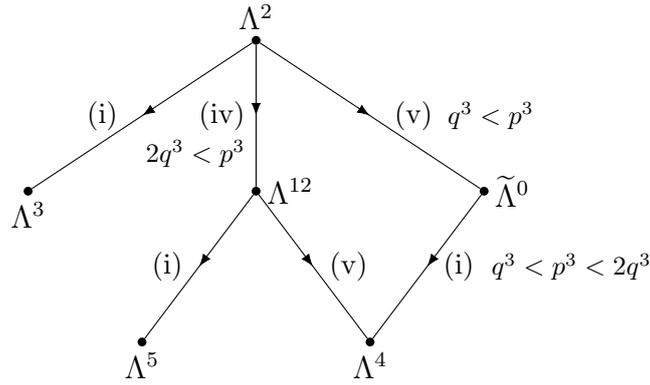
\begin{center}
  \begin{tikzpicture}
    \node[circle,fill=black,scale=0.3] (L2) at (0,10) {};
    \node[circle,fill=black,scale=0.3] (Ld1) at (-3,8) {};
    \node[circle,fill=black,scale=0.3] (Ld4) at (0,8) {};
    \node[circle,fill=black,scale=0.3] (Ld5) at (3,8) {};
    \node[circle,fill=black,scale=0.3] (Ldd1) at (-1.5,6) {};
    \node[circle,fill=black,scale=0.3] (Ldd5) at (1.5,6) {};
    \draw[decoration={markings, mark=at position 0.5 with {\arrow{Latex}}},postaction={decorate}] (L2) to (Ld1);
    \node[left,xshift=-5] at ($(L2)!0.5!(Ld1)$) {\small (i)};
    \draw[decoration={markings, mark=at position 0.5 with {\arrow{Latex}}},postaction={decorate}] (L2) to (Ld4);
    \node[left] at ($(L2)!0.5!(Ld4)$) {\small (iv)};
    \draw[decoration={markings, mark=at position 0.5 with {\arrow{Latex}}},postaction={decorate}] (L2) to (Ld5);
    \node[right,xshift=5] at ($(L2)!0.5!(Ld5)$) {\small (v)};
    \draw[decoration={markings, mark=at position 0.5 with {\arrow{Latex}}},postaction={decorate}] (Ld4) to (Ldd1);
    \node[left,xshift=-2] at ($(Ld4)!0.5!(Ldd1)$) {\small (i)};
    \draw[decoration={markings, mark=at position 0.5 with {\arrow{Latex}}},postaction={decorate}] (Ld4) to (Ldd5);
    \node[right,xshift=2] at ($(Ld4)!0.5!(Ldd5)$) {\small (v)};
    \draw[decoration={markings, mark=at position 0.5 with {\arrow{Latex}}},postaction={decorate}] (Ld5) to (Ldd5);
    \node[right,xshift=2] at ($(Ld5)!0.5!(Ldd5)$) {\small (i)};
    \node[above] at (L2) {$\L^2$}; 
    \node[below] at (Ld1) {$\L^3$} ; 
    \node[right] at (Ld4) {$\L^{12}$} ; 
    \node[right] at (Ld5) {$\widetilde{\L}^{0}$} ; 
    \node[below] at (Ldd1) {$\L^5$} ; 
    \node[below] at (Ldd5) {$\L^4$} ; 
    \node[left] at ($(L2)!0.75!(Ld4)$) {{\footnotesize $  2q^3 < p^3$}};
    \node[right,xshift=25] at ($(L2)!0.5!(Ld5)$) {{\footnotesize $  q^3 < p^3$}};
    \node[right,xshift=20] at ($(Ld5)!0.5!(Ldd5)$){{\footnotesize $  q^3 < p^3 < 2q^3$}};
  \end{tikzpicture}
\end{center}
%


\subsubsection{Invariant modules in $  V^{\L^3} $}


\medskip\noindent
After some analysis we find that  $ V^{\L^3} $ 
has the following SVs:

\begin{itemize}
   \item $ 2q^3 <  p^3  \ \Rightarrow \ $ type (ii) (iii) (iv) SVs of weights ~$ \L^4$, $\widetilde{\L}^{0}$, $\L^5$, resp.
    \item $ q^3 < p^3 < 2q^3 \ \Rightarrow \ $ type (ii) SV of weights $ \L^4. $
\end{itemize}

There are no further embeddings.  Thus, the embedding picture of $ V^{\L^3} $ 
is given as follows:

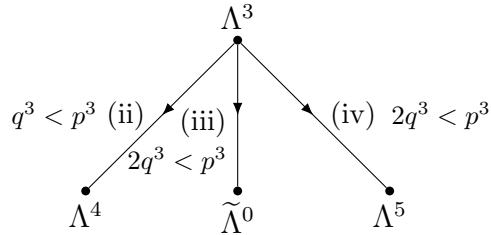
\begin{center}
  \begin{tikzpicture}
    \node[circle,fill=black,scale=0.3] (L3) at (0,10) {};
    \node[circle,fill=black,scale=0.3] (Ld2) at (-2,8) {};
    \node[circle,fill=black,scale=0.3] (Ld3) at (0,8) {};
    \node[circle,fill=black,scale=0.3] (Ld4) at (2,8) {};
    \draw[decoration={markings, mark=at position 0.5 with {\arrow{Latex}}},postaction={decorate}] (L3) to (Ld2);
    \node[left,xshift=-2] at ($(L3)!0.5!(Ld2)$) {\small (ii)};
    \draw[decoration={markings, mark=at position 0.5 with {\arrow{Latex}}},postaction={decorate}] (L3) to (Ld3);
    \node[left] at ($(L3)!0.55!(Ld3)$) {\small (iii)};
    \draw[decoration={markings, mark=at position 0.5 with {\arrow{Latex}}},postaction={decorate}] (L3) to (Ld4);
    \node[right,xshift=2] at ($(L3)!0.5!(Ld4)$) {\small (iv)};
    \node[above] at (L3)  {$\Lambda^3 $}; 
    \node[below] at (Ld2) {$\Lambda^4 $}; 
    \node[below] at (Ld3) {$\widetilde{\L}^{0} $}; 
    \node[below] at (Ld4){$\Lambda^5 $}; 
    \node[left,xshift=-20] at ($(L3)!0.5!(Ld2)$) {\footnotesize $ q^3 < p^3$};
    \node[left] at ($(L3)!0.8!(Ld3)$) {\footnotesize $ 2q^3 < p^3$};
    \node[right,xshift=25] at ($(L3)!0.5!(Ld4)$) {\footnotesize $ 2q^3 < p^3$};
  \end{tikzpicture}
\end{center}
%


\subsubsection{Invariant submodules in $ V^{\L^4} $}


\medskip\noindent
After some analysis we find that  $ V^{\L^4} $ has the following SVs:
\begin{itemize}
   \item $ 2q^3 <  p^3  \ \Rightarrow \ $ type (i) SV of weight $\widetilde{\L}^{0}$;
   \item $ p^3 < q^3   \ \Rightarrow \ $ type (ii) SV of weight ~$\L^3$;
\end{itemize}

There are no further embeddings.  Thus, the embedding picture of $ V^{\L^4} $ is given as follows:

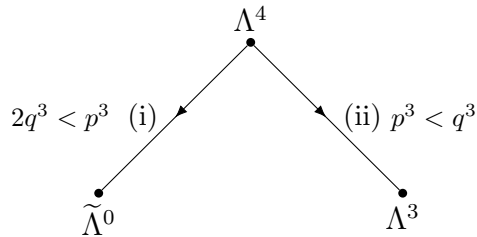
\begin{center}
  \begin{tikzpicture}
    \node[circle,fill=black,scale=0.3] (L4) at (0,10) {};
    \node[circle,fill=black,scale=0.3] (Ld1) at (-2,8) {};
    \node[circle,fill=black,scale=0.3] (Ld2) at (2,8) {};
    \draw[decoration={markings, mark=at position 0.5 with {\arrow{Latex}}},postaction={decorate}] (L4) to (Ld1);
    \node[left,xshift=-2] at ($(L4)!0.5!(Ld1)$) {\small (i)};
    \draw[decoration={markings, mark=at position 0.5 with {\arrow{Latex}}},postaction={decorate}] (L4) to (Ld2);
    \node[right,xshift=2] at ($(L4)!0.5!(Ld2)$) {\small (ii)};
    \node[above] at (L4){$\Lambda^4 $}; 
    \node[below] at (Ld1) {$\widetilde{\L}^{0} $}; 
    \node[below] at (Ld2) {$\Lambda^3 $}; 
    \node[left,xshift=-20] at ($(L4)!0.5!(Ld1)$) {\footnotesize{$ 2q^3 < p^3$}};
    \node[right,xshift=20] at ($(L4)!0.5!(Ld2)$) {\footnotesize{$p^3 < q^3 $}};
  \end{tikzpicture}
\end{center}
%


\subsubsection{Invariant modules in $ V^{\L^5} $}


\medskip\noindent
After some analysis we find that  $ V^{\L^5} $ has the following SVs:

\begin{itemize}
  \item $ q^3 < p^3 < 2q^3 \ \Rightarrow \ $ type (ii) (iii) (iv) SVs of weights ~$\widetilde{\L}^{0}$, $\L^4$, $\L^3$, resp.
  \item $  2q^3 < p^3 \ \Rightarrow \ $ type (ii)  SV of weight $ \widetilde{\L}^{0}.$ 
\end{itemize}

Further embedding in $ V^{\widetilde{\L}^{0}} $ 
where $ q^3 < p^3 $ shows that  it has the type (i) SV of weight ~$\L^4$
for $ q^3 < p^3 < 2q^3.$

There are no further embeddings.  Thus, the embedding picture of  $ V^{\L^5} $ 
is given as follows:

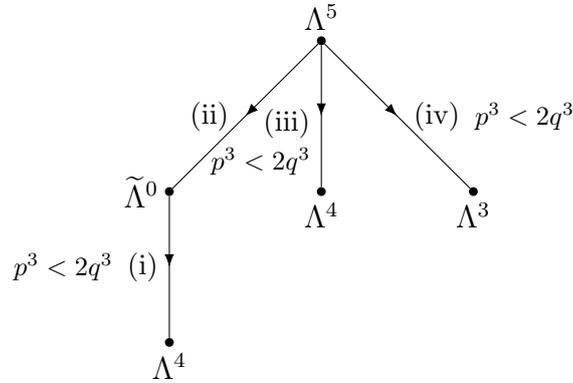
\begin{center}
  \begin{tikzpicture}
    \node[circle,fill=black,scale=0.3] (L5) at (0,10) {};
    \node[circle,fill=black,scale=0.3] (Ld2) at (-2,8) {};
    \node[circle,fill=black,scale=0.3] (Ld3) at (0,8) {};
    \node[circle,fill=black,scale=0.3] (Ld4) at (2,8) {};
    \node[circle,fill=black,scale=0.3] (Ldd1) at (-2,6) {};
    \draw[decoration={markings, mark=at position 0.5 with {\arrow{Latex}}},postaction={decorate}] (L5) to (Ld2);
    \node[left,xshift=-2] at ($(L5)!0.5!(Ld2)$) {\small (ii)};
    \draw[decoration={markings, mark=at position 0.5 with {\arrow{Latex}}},postaction={decorate}] (L5) to (Ld3);
    \node[left] at ($(L5)!0.55!(Ld3)$) {\small (iii)};
    \draw[decoration={markings, mark=at position 0.5 with {\arrow{Latex}}},postaction={decorate}] (L5) to (Ld4);
    \node[right,xshift=2] at ($(L5)!0.5!(Ld4)$) {\small (iv)};
    \draw[decoration={markings, mark=at position 0.5 with {\arrow{Latex}}},postaction={decorate}] (Ld2) to (Ldd1);
    \node[left] at ($(Ld2)!0.5!(Ldd1)$) {\small (i)};
    %
    \node[above] at (L5) {$\Lambda^5 $}; 
    \node[left] at (Ld2) {$\widetilde{\L}^{0} $}; 
    \node[below] at (Ld3) {$\Lambda^4 $}; 
    \node[below] at (Ld4){$\Lambda^3 $}; 
    \node[below] at (Ldd1) {$\Lambda^4 $}; 
    \node[left,xshift=-17] at ($(Ld2)!0.5!(Ldd1)$)  {\footnotesize $ p^3 < 2q^3$};
    \node[left] at ($(L5)!0.8!(Ld3)$) {\footnotesize $ p^3 < 2q^3$};
    \node[right,xshift=25] at ($(L5)!0.5!(Ld4)$) {\footnotesize $ p^3 < 2q^3$};
  \end{tikzpicture}
\end{center}
%
This completes the embeddings in  \textbf{A12345}.

\subsubsection{Complete diagram for \textbf{A12345} and subcases}

Doing the complete diagram it may be better to distinguish the cases, since pairs of modules have submodules of opposite directions depending on the values of the parameters $p^3,q^3$. The distinction of cases was already done in \eqref{emb234},\eqref{emb2345},\eqref{emb12345} and so we proceed:

 \bigskip\noindent
\textbf{A234}

The LW is $ \L^0 $ with the constraint $ p^3 < q^3 $ which we denote by $ \L^0_{234}.$
For the constraint $ p^3 < q^3 $ we have the following embedding patterns.
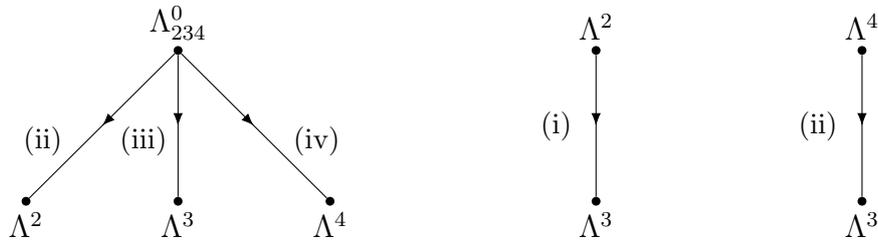
\begin{center}
  \begin{tikzpicture}
    \node[circle,fill=black,scale=0.3] (L234) at (0,10) {};
    \node[circle,fill=black,scale=0.3] (Ld2) at (-2,8) {};
    \node[circle,fill=black,scale=0.3] (Ld3) at (0,8) {};
    \node[circle,fill=black,scale=0.3] (Ld4) at (2,8) {};
    \draw[decoration={markings, mark=at position 0.5 with {\arrow{Latex}}},postaction={decorate}] (L234) to (Ld2);
    \node[left,xshift=-5] at ($(L234)!0.6!(Ld2)$) {\small (ii)};
    \draw[decoration={markings, mark=at position 0.5 with {\arrow{Latex}}},postaction={decorate}] (L234) to (Ld3);
    \node[left] at ($(L234)!0.6!(Ld3)$) {\small (iii)};
    \draw[decoration={markings, mark=at position 0.5 with {\arrow{Latex}}},postaction={decorate}] (L234) to (Ld4);
    \node[right,xshift=5] at ($(L234)!0.6!(Ld4)$) {\small (iv)};
    \node[above] at (L234) {$\Lambda^0_{234}$};
    \node[below] at (Ld2) {$\Lambda^2$};
    \node[below] at (Ld3) {$\Lambda^3$};
    \node[below] at (Ld4) {$\Lambda^4$};
    \begin{scope}[xshift=55mm]
      \node[circle,fill=black,scale=0.3] (L2) at (0,10) {};
      \node[circle,fill=black,scale=0.3] (L23) at (0,8) {};
      \draw[decoration={markings, mark=at position 0.5 with {\arrow{Latex}}},postaction={decorate}] (L2) to (L23);
      \node[left,xshift=-5] at ($(L2)!0.5!(L23)$) {\small (i)};
      \node[above] at (L2) {$\Lambda^2$};
      \node[below] at (L23) {$\Lambda^3$};
    \end{scope}
    \begin{scope}[xshift=90mm]
      \node[circle,fill=black,scale=0.3] (L2) at (0,10) {};
      \node[circle,fill=black,scale=0.3] (L23) at (0,8) {};
      \draw[decoration={markings, mark=at position 0.5 with {\arrow{Latex}}},postaction={decorate}] (L2) to (L23);
      \node[left,xshift=-5] at ($(L2)!0.5!(L23)$) {\small (ii)};
      \node[above] at (L2) {$\Lambda^4$};
      \node[below] at (L23) {$\Lambda^3$};
    \end{scope}
   \end{tikzpicture}
\end{center}
$ V^{\Lambda^3}$ does not have SV.
Combining these diagrams, we obtain the following diagram.

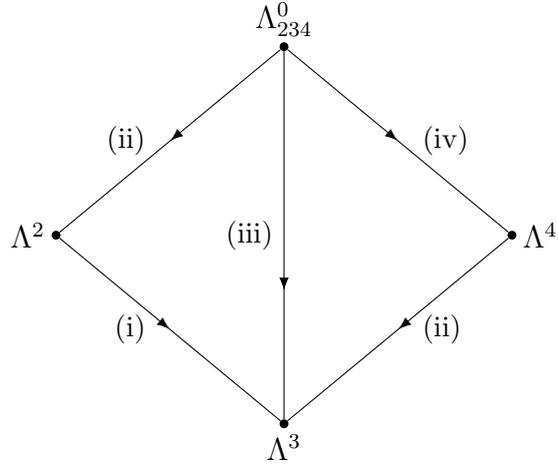
\eqn{A234d}
  \begin{tikzpicture}
    \node[circle,fill=black,scale=0.3] (L245) at (0,10) {};
    \node[circle,fill=black,scale=0.3] (Ld2) at (-3,7.5) {};
    \node[circle,fill=black,scale=0.3] (Ld4) at (0,5) {};
    \node[circle,fill=black,scale=0.3] (Ld5) at (3,7.5) {};
    \draw[decoration={markings, mark=at position 0.5 with {\arrow{Latex}}},postaction={decorate}] (L245) to (Ld2);
    \node[left,xshift=-5] at ($(L245)!0.5!(Ld2)$) {\small (ii)};
    \draw[decoration={markings, mark=at position 0.65 with {\arrow{Latex}}},postaction={decorate}] (L245) to (Ld4);
    \node[left] at ($(L245)!0.5!(Ld4)$) {\small (iii)};
    \draw[decoration={markings, mark=at position 0.5 with {\arrow{Latex}}},postaction={decorate}] (L245) to (Ld5);
    \node[right,xshift=5] at ($(L245)!0.5!(Ld5)$) {\small (iv)};
    \draw[decoration={markings, mark=at position 0.5 with {\arrow{Latex}}},postaction={decorate}] (Ld2) to (Ld4);
    \node[left,xshift=-5] at ($(Ld2)!0.5!(Ld4)$) {\small (i)};
    \draw[decoration={markings, mark=at position 0.5 with {\arrow{Latex}}},postaction={decorate}] (Ld5) to (Ld4);
    \node[right,xshift=5] at ($(Ld5)!0.5!(Ld4)$) {\small (ii)};
    \node[above] at (L245) {$\Lambda^0_{234}$};
    \node[left] at (Ld2) {$\Lambda^2$};
    \node[below] at (Ld4) {$\Lambda^3$};
    \node[right] at (Ld5) {$\Lambda^4$};
  \end{tikzpicture}
\ee

\bigskip\bigskip\noindent
\textbf{A2345}

The LW is $ \L^0 $ with the constraint $ q^3 < p^3 < 2q^3 $ which we denote by $ \L^0_{2345}.$
For the constraint $ q^3 < p^3 < 2q^3 $ we have the following embedding patterns.

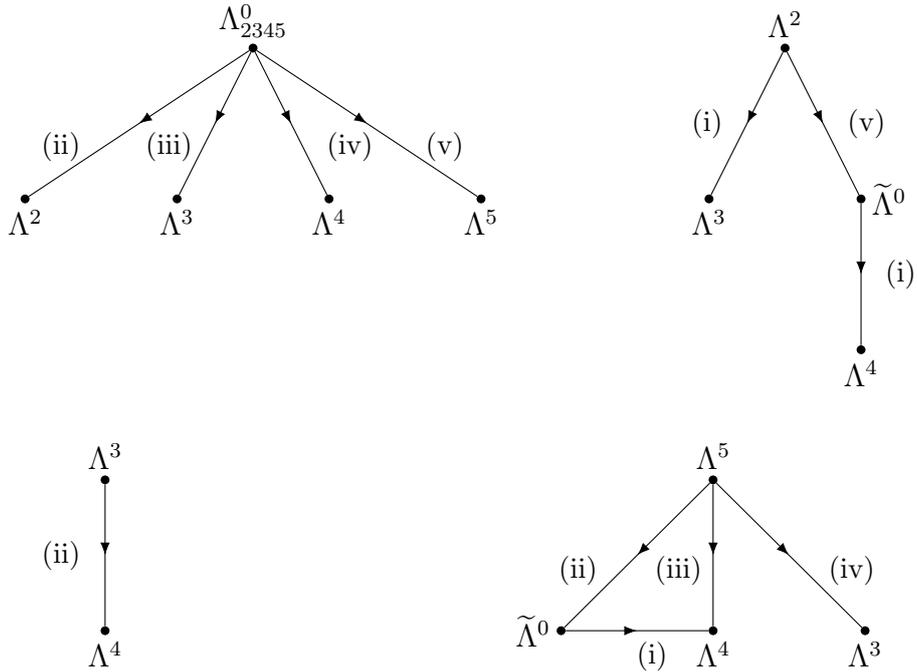
\begin{center}
  \begin{tikzpicture}
    \node[circle,fill=black,scale=0.3] (L234) at (0,10) {};
    \node[circle,fill=black,scale=0.3] (Ld2) at (-3,8) {};
    \node[circle,fill=black,scale=0.3] (Ld3) at (-1,8) {};
    \node[circle,fill=black,scale=0.3] (Ld4) at (1,8) {};
    \node[circle,fill=black,scale=0.3] (Ld5) at (3,8) {};
    \draw[decoration={markings, mark=at position 0.5 with {\arrow{Latex}}},postaction={decorate}] (L234) to (Ld2);
    \node[left,xshift=-5] at ($(L234)!0.65!(Ld2)$) {\small (ii)};
    \draw[decoration={markings, mark=at position 0.5 with {\arrow{Latex}}},postaction={decorate}] (L234) to (Ld3);
    \node[left] at ($(L234)!0.65!(Ld3)$) {\small (iii)};
    \draw[decoration={markings, mark=at position 0.5 with {\arrow{Latex}}},postaction={decorate}] (L234) to (Ld4);
    \node[right,xshift=5] at ($(L234)!0.65!(Ld4)$) {\small (iv)};
    \draw[decoration={markings, mark=at position 0.5 with {\arrow{Latex}}},postaction={decorate}] (L234) to (Ld5);
    \node[right,xshift=5] at ($(L234)!0.65!(Ld5)$) {\small (v)};
    \node[above] at (L234) {$\Lambda^0_{2345}$};
    \node[below] at (Ld2) {$\Lambda^2$};
    \node[below] at (Ld3) {$\Lambda^3$};
    \node[below] at (Ld4) {$\Lambda^4$};
    \node[below] at (Ld5) {$\Lambda^5$};
    \begin{scope}[xshift=70mm]
      \node[circle,fill=black,scale=0.3] (L2) at (0,10) {};
      \node[circle,fill=black,scale=0.3] (L23) at (-1,8) {};
      \node[circle,fill=black,scale=0.3] (L214) at (1,8) {};
      \node[circle,fill=black,scale=0.3] (L24) at (1,6) {};
      \draw[decoration={markings, mark=at position 0.5 with {\arrow{Latex}}},postaction={decorate}] (L2) to (L23);
      \node[left,xshift=-5] at ($(L2)!0.5!(L23)$) {\small (i)};
      \draw[decoration={markings, mark=at position 0.5 with {\arrow{Latex}}},postaction={decorate}] (L2) to (L214);
      \node[right,xshift=5] at ($(L2)!0.5!(L214)$) {\small (v)};
      \draw[decoration={markings, mark=at position 0.5 with {\arrow{Latex}}},postaction={decorate}] (L214) to (L24);
      \node[right,xshift=5] at ($(L214)!0.5!(L24)$) {\small (i)};
      \node[above] at (L2) {$\Lambda^2$};
      \node[below] at (L23) {$\Lambda^3$};
      \node[right] at (L214) {$\widetilde{\L}^{0}$};
      \node[below] at (L24) {$\Lambda^4$};
    \end{scope}
   \end{tikzpicture}
\end{center}

\begin{center}
  \begin{tikzpicture}
    \node[circle,fill=black,scale=0.3] (L3) at (0,10) {};
    \node[circle,fill=black,scale=0.3] (L34) at (0,8) {};
    \draw[decoration={markings, mark=at position 0.5 with {\arrow{Latex}}},postaction={decorate}] (L3) to (L34);
    \node[left,xshift=-5] at ($(L3)!0.5!(L34)$) {\small (ii)};
    \node[above] at (L3) {$\Lambda^3$};
    \node[below] at (L34) {$\Lambda^4$};
    \begin{scope}[xshift=80mm]
      \node[circle,fill=black,scale=0.3] (L5) at (0,10) {};
      \node[circle,fill=black,scale=0.3] (L514) at (-2,8) {};
      \node[circle,fill=black,scale=0.3] (L54) at (0,8) {};
      \node[circle,fill=black,scale=0.3] (L53) at (2,8) {};
      \draw[decoration={markings, mark=at position 0.5 with {\arrow{Latex}}},postaction={decorate}] (L5) to (L514);
      \node[left,xshift=-5] at ($(L5)!0.6!(L514)$) {\small (ii)};
      \draw[decoration={markings, mark=at position 0.5 with {\arrow{Latex}}},postaction={decorate}] (L5) to (L54);
      \node[left] at ($(L5)!0.6!(L54)$) {\small (iii)};
      \draw[decoration={markings, mark=at position 0.5 with {\arrow{Latex}}},postaction={decorate}] (L5) to (L53);
      \node[right,xshift=5] at ($(L5)!0.6!(L53)$) {\small (iv)};
      \draw[decoration={markings, mark=at position 0.5 with {\arrow{Latex}}},postaction={decorate}] (L514) to (L54);
      \node[below] at ($(L514)!0.6!(L54)$) {\small (i)};
      \node[above] at (L5) {$\Lambda^5$};
      \node[left] at (L514) {$\widetilde{\L}^{0}$};
      \node[below] at (L54) {$\Lambda^4$};
      \node[below] at (L53) {$\Lambda^3$};
    \end{scope}
   \end{tikzpicture}
\end{center}
$ V^{\Lambda^4}$ does not have SV.

Combining these diagrams, we obtain the following embedding diagram:

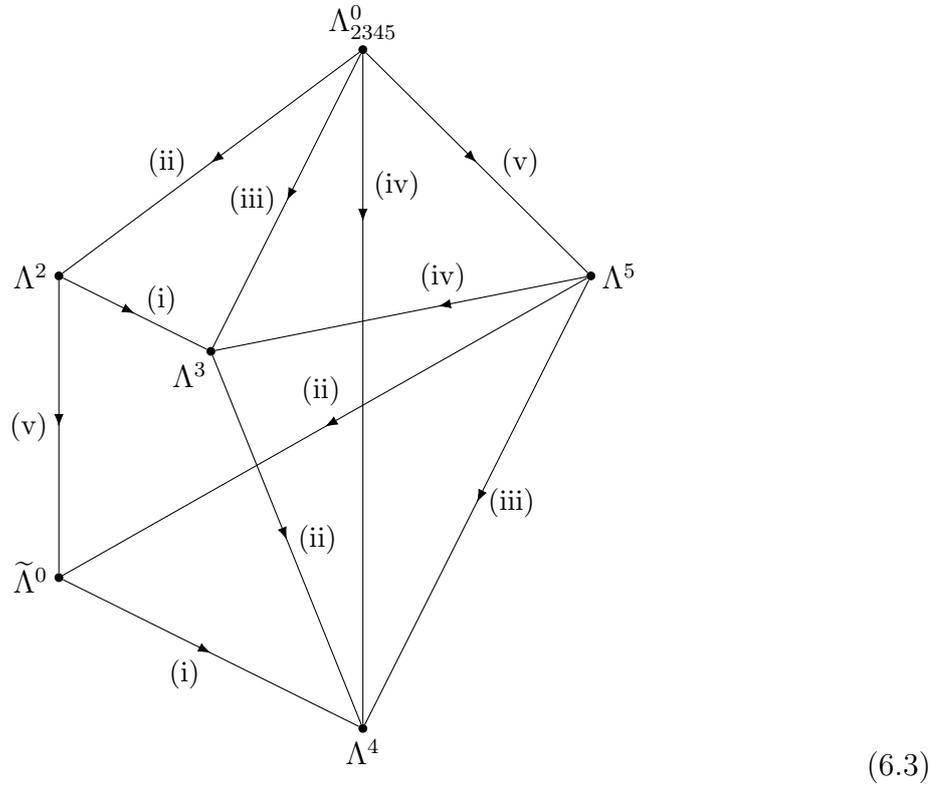
\eqn{A2345d}
  \begin{tikzpicture}
    \node[circle,fill=black,scale=0.3] (L234) at (0,10) {};
    \node[above] at (L234) {$\Lambda^0_{2345}$};
    \node[circle,fill=black,scale=0.3] (Ld2) at (-4,7) {};
    \node[left] at (Ld2) {$\Lambda^2$};
    \node[circle,fill=black,scale=0.3] (Ld3) at (-2,6) {};
    \node[below,xshift=-8] at (Ld3) {$\Lambda^3$};
    \node[circle,fill=black,scale=0.3] (Ld4) at (0,1) {};
    \node[below] at (Ld4) {$\Lambda^4$};
    \node[circle,fill=black,scale=0.3] (Ld5) at (3,7) {};
    \node[right] at (Ld5) {$\Lambda^5$};
    \node[circle,fill=black,scale=0.3] (Ld14) at (-4,3) {};
    \node[left] at (Ld14) {$\widetilde{\L}^{0}$};
    \draw[decoration={markings, mark=at position 0.5 with {\arrow{Latex}}},postaction={decorate}] (L234) to (Ld2);
    \node[left,xshift=-5] at ($(L234)!0.5!(Ld2)$) {\small (ii)};
    \draw[decoration={markings, mark=at position 0.5 with {\arrow{Latex}}},postaction={decorate}] (L234) to (Ld3);
    \node[left] at ($(L234)!0.5!(Ld3)$) {\small (iii)};
    \draw[decoration={markings, mark=at position 0.25 with {\arrow{Latex}}},postaction={decorate}] (L234) to (Ld4);
    \node[right] at ($(L234)!0.2!(Ld4)$) {\small (iv)};
    \draw[decoration={markings, mark=at position 0.5 with {\arrow{Latex}}},postaction={decorate}] (L234) to (Ld5);
    \node[right,xshift=5] at ($(L234)!0.5!(Ld5)$) {\small (v)};
    \draw[decoration={markings, mark=at position 0.5 with {\arrow{Latex}}},postaction={decorate}] (Ld2) to (Ld14);
    \node[left] at ($(Ld2)!0.5!(Ld14)$) {\small (v)};
    \draw[decoration={markings, mark=at position 0.5 with {\arrow{Latex}}},postaction={decorate}] (Ld2) to (Ld3);
    \node[right,yshift=5] at ($(Ld2)!0.5!(Ld3)$) {\small (i)};
    \draw[decoration={markings, mark=at position 0.5 with {\arrow{Latex}}},postaction={decorate}] (Ld3) to (Ld4);
    \node[right] at ($(Ld3)!0.5!(Ld4)$) {\small (ii)};
    \draw[decoration={markings, mark=at position 0.5 with {\arrow{Latex}}},postaction={decorate}] (Ld14) to (Ld4);
    \node[left,yshift=-8] at ($(Ld14)!0.5!(Ld4)$) {\small (i)};
    \draw[decoration={markings, mark=at position 0.4 with {\arrow{Latex}}},postaction={decorate}] (Ld5) to (Ld3);
    \node[left,yshift=8] at ($(Ld5)!0.3!(Ld3)$) {\small (iv)};
    \draw[decoration={markings, mark=at position 0.5 with {\arrow{Latex}}},postaction={decorate}] (Ld5) to (Ld14);
    \node[left,yshift=8] at ($(Ld5)!0.45!(Ld14)$) {\small (ii)};
    \draw[decoration={markings, mark=at position 0.5 with {\arrow{Latex}}},postaction={decorate}] (Ld5) to (Ld4);
    \node[right] at ($(Ld5)!0.5!(Ld4)$) {\small (iii)};
   \end{tikzpicture}
\ee

\bigskip\bigskip
\noindent
\textbf{A12345}

The LW is $ \L^0 $ with the constraint $ 2q^3 < p^3 $ for which we keep using the notation $ \L^0. $
%
%
For the constraint $ 2q^3 < p^3 $ we have the following embedding patterns:

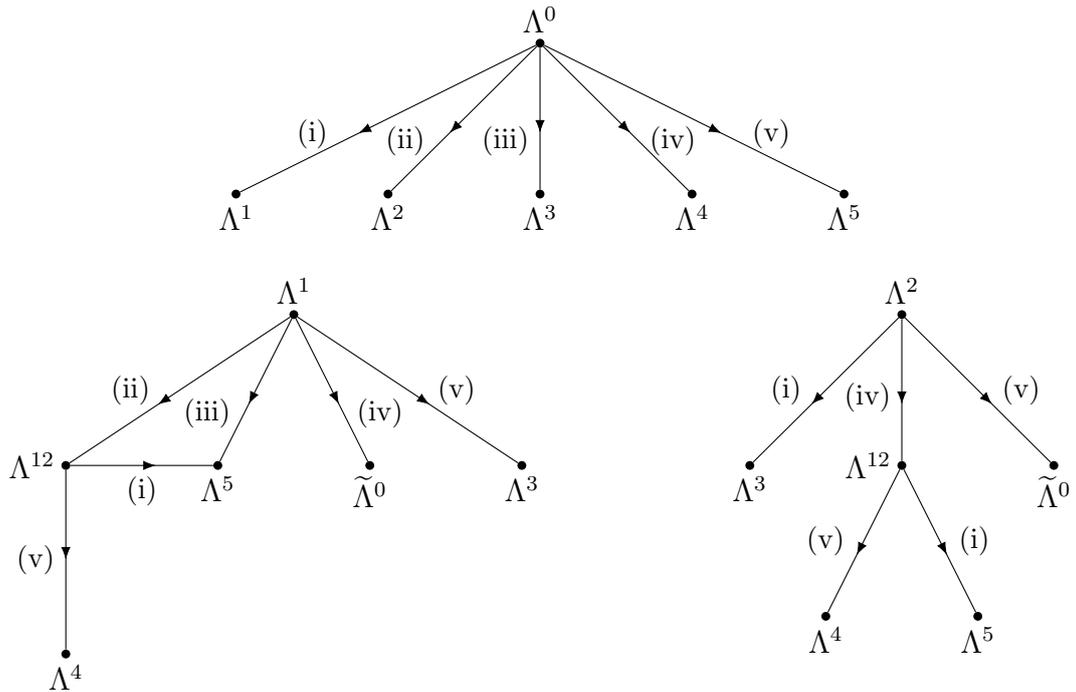
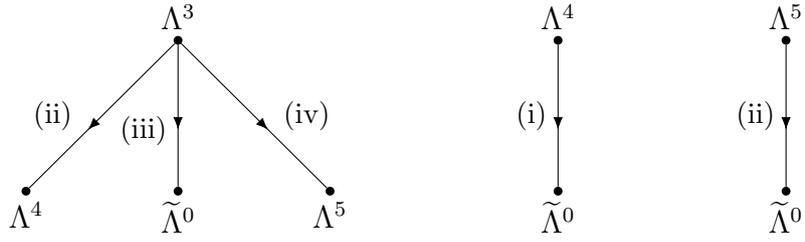
\begin{center}
  \begin{tikzpicture}
    \node[circle,fill=black,scale=0.3] (L3) at (0,10) {};
    \node[above] at (L3) {$\Lambda^0$};
    \node[circle,fill=black,scale=0.3] (Ld1) at (-4,8) {};
    \node[below] at (Ld1) {$\Lambda^1$};
    \node[circle,fill=black,scale=0.3] (Ld2) at (-2,8) {};
    \node[below] at (Ld2) {$\Lambda^2$};
    \node[circle,fill=black,scale=0.3] (Ld3) at (0,8) {};
    \node[below] at (Ld3) {$\Lambda^3$};
    \node[circle,fill=black,scale=0.3] (Ld4) at (2,8) {};
    \node[below] at (Ld4) {$\Lambda^4$};
    \node[circle,fill=black,scale=0.3] (Ld5) at (4,8) {};
    \node[below] at (Ld5) {$\Lambda^5$};
    \draw[decoration={markings, mark=at position 0.6 with {\arrow{Latex}}},postaction={decorate}] (L3) to (Ld1);
    \node[left,xshift=-7] at ($(L3)!0.6!(Ld1)$) {\small (i)};
    \draw[decoration={markings, mark=at position 0.6 with {\arrow{Latex}}},postaction={decorate}] (L3) to (Ld2);
    \node[left,xshift=-2] at ($(L3)!0.65!(Ld2)$) {\small (ii)};
    \draw[decoration={markings, mark=at position 0.6 with {\arrow{Latex}}},postaction={decorate}] (L3) to (Ld3);
    \node[left] at ($(L3)!0.65!(Ld3)$) {\small (iii)};
    \draw[decoration={markings, mark=at position 0.6 with {\arrow{Latex}}},postaction={decorate}] (L3) to (Ld4);
    \node[right] at ($(L3)!0.65!(Ld4)$) {\small (iv)};
    \draw[decoration={markings, mark=at position 0.6 with {\arrow{Latex}}},postaction={decorate}] (L3) to (Ld5);
    \node[right,xshift=7] at ($(L3)!0.6!(Ld5)$) {\small (v)};

  \end{tikzpicture}
\end{center}

\begin{center}
  \begin{tikzpicture}
    \node[circle,fill=black,scale=0.3] (L1) at (0,10) {};
    \node[above] at (L1) {$\Lambda^1$};
    \node[circle,fill=black,scale=0.3] (L12) at (-3,8) {};
    \node[left] at (L12) {$\Lambda^{12}$};
    \node[circle,fill=black,scale=0.3] (L5) at (-1,8) {};
    \node[below] at (L5) {$\Lambda^5$};
    \node[circle,fill=black,scale=0.3] (L14) at (1,8) {};
    \node[below] at (L14) {$\widetilde{\L}^{0}$};
    \node[circle,fill=black,scale=0.3] (L3) at (3,8) {};
    \node[below] at (L3) {$\Lambda^3$};
    \node[circle,fill=black,scale=0.3] (L4) at (-3,5.5) {};
    \node[below] at (L4) {$\Lambda^4$};
    \draw[decoration={markings, mark=at position 0.6 with {\arrow{Latex}}},postaction={decorate}] (L1) to (L12);
    \node[left,xshift=-7] at ($(L1)!0.5!(L12)$) {\small (ii)};
    \draw[decoration={markings, mark=at position 0.6 with {\arrow{Latex}}},postaction={decorate}] (L1) to (L5);
    \node[left,xshift=-1] at ($(L1)!0.65!(L5)$) {\small (iii)};
    \draw[decoration={markings, mark=at position 0.6 with {\arrow{Latex}}},postaction={decorate}] (L1) to (L14);
    \node[right,xshift=1] at ($(L1)!0.65!(L14)$) {\small (iv)};
    \draw[decoration={markings, mark=at position 0.6 with {\arrow{Latex}}},postaction={decorate}] (L1) to (L3);
    \node[right,xshift=7] at ($(L1)!0.5!(L3)$) {\small (v)};
    \draw[decoration={markings, mark=at position 0.5 with {\arrow{Latex}}},postaction={decorate}] (L12) to (L4);
    \node[left] at ($(L12)!0.5!(L4)$) {\small (v)};
    \draw[decoration={markings, mark=at position 0.6 with {\arrow{Latex}}},postaction={decorate}] (L12) to (L5);
    \node[below] at ($(L12)!0.5!(L5)$) {\small (i)};
    \begin{scope}[xshift=80mm]
      \node[circle,fill=black,scale=0.3] (L2) at (0,10) {};
      \node[above] at (L2) {$\Lambda^2$};
      \node[circle,fill=black,scale=0.3] (Ld3) at (-2,8) {};
      \node[below] at (Ld3) {$\Lambda^3$};
      \node[circle,fill=black,scale=0.3] (Ld12) at (0,8) {};
      \node[left] at (Ld12) {$\Lambda^{12}$};
      \node[circle,fill=black,scale=0.3] (Ld14) at (2,8) {};
      \node[below] at (Ld14) {$\widetilde{\L}^{0}$};
      \node[circle,fill=black,scale=0.3] (Ld4) at (-1,6) {};
      \node[below] at (Ld4) {$\Lambda^4$};
      \node[circle,fill=black,scale=0.3] (Ld5) at (1,6) {};
      \node[below] at (Ld5) {$\Lambda^5$};
      \draw[decoration={markings, mark=at position 0.6 with {\arrow{Latex}}},postaction={decorate}] (L2) to (Ld3);
      \node[left,xshift=-5] at ($(L2)!0.5!(Ld3)$) {\small (i)};
      \draw[decoration={markings, mark=at position 0.6 with {\arrow{Latex}}},postaction={decorate}] (L2) to (Ld12);
      \node[left] at ($(L2)!0.55!(Ld12)$) {\small (iv)};
      \draw[decoration={markings, mark=at position 0.6 with {\arrow{Latex}}},postaction={decorate}] (L2) to (Ld14);
      \node[right,xshift=5] at ($(L2)!0.5!(Ld14)$) {\small (v)};
      \draw[decoration={markings, mark=at position 0.6 with {\arrow{Latex}}},postaction={decorate}] (Ld12) to (Ld4);
      \node[left,xshift=-3] at ($(Ld12)!0.5!(Ld4)$) {\small (v)};
      \draw[decoration={markings, mark=at position 0.6 with {\arrow{Latex}}},postaction={decorate}] (Ld12) to (Ld5);
      \node[right,xshift=3] at ($(Ld12)!0.5!(Ld5)$) {\small (i)};
    \end{scope}
  \end{tikzpicture}
\end{center}

\begin{center}
  \begin{tikzpicture}
    \node[circle,fill=black,scale=0.3] (L3) at (0,10) {};
    \node[above] at (L3) {$\Lambda^3$};
    \node[circle,fill=black,scale=0.3] (L4) at (-2,8) {};
    \node[below] at (L4) {$\Lambda^4$};
    \node[circle,fill=black,scale=0.3] (L14) at (0,8) {};
    \node[below] at (L14) {$\widetilde{\L}^{0}$};
    \node[circle,fill=black,scale=0.3] (L5) at (2,8) {};
    \node[below] at (L5) {$\Lambda^{5}$};
    \draw[decoration={markings, mark=at position 0.6 with {\arrow{Latex}}},postaction={decorate}] (L3) to (L4);
    \node[left,xshift=-7] at ($(L3)!0.5!(L4)$) {\small (ii)};
    \draw[decoration={markings, mark=at position 0.6 with {\arrow{Latex}}},postaction={decorate}] (L3) to (L14);
    \node[left] at ($(L3)!0.6!(L14)$) {\small (iii)};
    \draw[decoration={markings, mark=at position 0.6 with {\arrow{Latex}}},postaction={decorate}] (L3) to (L5);
    \node[right,xshift=7] at ($(L3)!0.5!(L5)$) {\small (iv)};
    \begin{scope}[xshift=50mm]
       \node[circle,fill=black,scale=0.3] (Ld4) at (0,10) {};
       \node[above] at (Ld4) {$\Lambda^4$};
       \node[circle,fill=black,scale=0.3] (Ld14) at (0,8) {};
       \node[below] at (Ld14) {$\widetilde{\L}^{0}$};
       \draw[decoration={markings, mark=at position 0.6 with {\arrow{Latex}}},postaction={decorate}] (Ld4) to (Ld14);
       \node[left] at ($(Ld4)!0.5!(Ld14)$) {\small (i)};
    \end{scope}
    \begin{scope}[xshift=80mm]
       \node[circle,fill=black,scale=0.3] (Ld4) at (0,10) {};
       \node[above] at (Ld4) {$\Lambda^5$};
       \node[circle,fill=black,scale=0.3] (Ld14) at (0,8) {};
       \node[below] at (Ld14) {$\widetilde{\L}^{0}$};
       \draw[decoration={markings, mark=at position 0.6 with {\arrow{Latex}}},postaction={decorate}] (Ld4) to (Ld14);
       \node[left] at ($(Ld4)!0.5!(Ld14)$) {\small (ii)};
    \end{scope}

  \end{tikzpicture}
\end{center}
$ V^{\widetilde{\L}^{0}}$ does not have SV.

Combining these diagrams, we obtain the complete embedding diagram of this case:

\eqn{A12345d}
  \begin{tikzpicture}
    \node[circle,fill=black,scale=0.3] (L0) at (0,9) {};
    \node[above] at (L0) {$\Lambda^0$};
    \node[circle,fill=black,scale=0.3] (L1) at (-4,7) {};
    \node[left] at (L1) {$\Lambda^1$};
    \node[circle,fill=black,scale=0.3] (L2) at (4,7) {};
    \node[right] at (L2) {$\Lambda^2$};
    \node[circle,fill=black,scale=0.3] (L3) at (0,2) {};
    \node[right] at (L3) {$\Lambda^3$};
    \node[circle,fill=black,scale=0.3] (L4) at (4,-3) {};
    \node[right] at (L4) {$\Lambda^4$};
    \node[circle,fill=black,scale=0.3] (L5) at (-4,-3) {};
    \node[left] at (L5) {$\Lambda^5$};
    \node[circle,fill=black,scale=0.3] (L12) at (4,2) {};
    \node[right] at (L12) {$\Lambda^{12}$};
    \node[circle,fill=black,scale=0.3] (L14) at (0,-5) {};
    \node[below] at (L14) {$\widetilde{\L}^{0}$};
    \draw[decoration={markings, mark=at position 0.5 with {\arrow{Latex}}},postaction={decorate}] (L0) to (L1);
    \node[left,xshift=-7] at ($(L0)!0.4!(L1)$) {\small (i)};
    \draw[decoration={markings, mark=at position 0.5 with {\arrow{Latex}}},postaction={decorate}] (L0) to (L2);
    \node[right,xshift=7] at ($(L0)!0.4!(L2)$) {\small (ii)};
    \draw[decoration={markings, mark=at position 0.5 with {\arrow{Latex}}},postaction={decorate}] (L0) to (L3);
    \node[right] at ($(L0)!0.5!(L3)$) {\small (iii)};
    \draw[decoration={markings, mark=at position 0.2 with {\arrow{Latex}}},postaction={decorate}] (L0) to (L4);
    \node[right] at ($(L0)!0.2!(L4)$) {\small (iv)};
    \draw[decoration={markings, mark=at position 0.2 with {\arrow{Latex}}},postaction={decorate}] (L0) to (L5);
    \node[left] at ($(L0)!0.2!(L5)$) {\small (v)};
    \draw[decoration={markings, mark=at position 0.75 with {\arrow{Latex}}},postaction={decorate}] (L1) to (L3);
    \node[right] at ($(L1)!0.7!(L3)$) {\small (v)};
    \draw[decoration={markings, mark=at position 0.2 with {\arrow{Latex}}},postaction={decorate}] (L1) to (L12);
    \node[right,yshift=2] at ($(L1)!0.2!(L12)$) {\small (ii)};
    \draw[decoration={markings, mark=at position 0.25 with {\arrow{Latex}}},postaction={decorate}] (L1) to (L14);
    \node[left] at ($(L1)!0.25!(L14)$) {\small (iv)};
    \draw[decoration={markings, mark=at position 0.5 with {\arrow{Latex}}},postaction={decorate}] (L1) to (L5);
    \node[left] at ($(L1)!0.5!(L5)$) {\small (iii)};
    \draw[decoration={markings, mark=at position 0.75 with {\arrow{Latex}}},postaction={decorate}] (L2) to (L3);
    \node[left,xshift=-2] at ($(L2)!0.7!(L3)$) {\small (i)};
    \draw[decoration={markings, mark=at position 0.25 with {\arrow{Latex}}},postaction={decorate}] (L2) to (L14);
    \node[right] at ($(L2)!0.25!(L14)$) {\small (v)};
    \draw[decoration={markings, mark=at position 0.5 with {\arrow{Latex}}},postaction={decorate}] (L2) to (L12);
    \node[right] at ($(L2)!0.5!(L12)$) {\small (iv)};
    \draw[decoration={markings, mark=at position 0.55 with {\arrow{Latex}}},postaction={decorate}] (L3) to (L4);
    \node[right] at ($(L3)!0.55!(L4)$) {\small (ii)};
    \draw[decoration={markings, mark=at position 0.6 with {\arrow{Latex}}},postaction={decorate}] (L3) to (L14);
    \node[left] at ($(L3)!0.6!(L14)$) {\small (iii)};
    \draw[decoration={markings, mark=at position 0.55 with {\arrow{Latex}}},postaction={decorate}] (L3) to (L5);
    \node[left] at ($(L3)!0.55!(L5)$) {\small (iv)};
    \draw[decoration={markings, mark=at position 0.5 with {\arrow{Latex}}},postaction={decorate}] (L12) to (L4);
    \node[right] at ($(L12)!0.5!(L4)$) {\small (v)};
    \draw[decoration={markings, mark=at position 0.8 with {\arrow{Latex}}},postaction={decorate}] (L12) to (L5);
    \node[right,yshift=-5] at ($(L12)!0.8!(L5)$) {\small (i)};
    \draw[decoration={markings, mark=at position 0.5 with {\arrow{Latex}}},postaction={decorate}] (L4) to (L14);
    \node[right,yshift=-7] at ($(L4)!0.5!(L14)$) {\small (i)};
    \draw[decoration={markings, mark=at position 0.5 with {\arrow{Latex}}},postaction={decorate}] (L5) to (L14);
    \node[left,yshift=-7] at ($(L5)!0.5!(L14)$) {\small (ii)};

  \end{tikzpicture}
\ee



\subsection{Case A245}

In \eqref{a245} we have found the initial embeddings of this case which are given by the diagram:

  \bigskip\bigskip

%
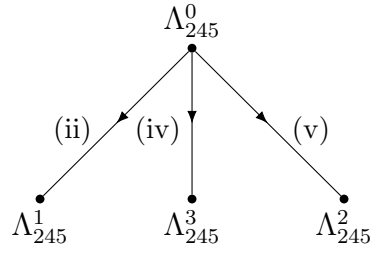
\begin{center}
  \begin{tikzpicture}
    \node[circle,fill=black,scale=0.3] (L245) at (0,10) {};
    \node[circle,fill=black,scale=0.3] (Ld2) at (-2,8) {};
    \node[circle,fill=black,scale=0.3] (Ld4) at (0,8) {};
    \node[circle,fill=black,scale=0.3] (Ld5) at (2,8) {};
    \draw[decoration={markings, mark=at position 0.5 with {\arrow{Latex}}},postaction={decorate}] (L245) to (Ld2);
    \node[left,xshift=-2] at ($(L245)!0.55!(Ld2)$) {\small (ii)};
    \draw[decoration={markings, mark=at position 0.5 with {\arrow{Latex}}},postaction={decorate}] (L245) to (Ld4);
    \node[left] at ($(L245)!0.55!(Ld4)$) {\small (iv)};
    \draw[decoration={markings, mark=at position 0.5 with {\arrow{Latex}}},postaction={decorate}] (L245) to (Ld5);
    \node[right,xshift=2] at ($(L245)!0.55!(Ld5)$) {\small (v)};
    \node[above] at (L245) {$\Lambda^0_{245}$};
    \node[below] at (Ld2) {$\Lambda^1_{245} $};
    \node[below] at (Ld4) {$\Lambda^3_{245} $};
    \node[below] at (Ld5) {$\Lambda^2_{245} $};
  \end{tikzpicture}
\end{center}

where ~$\Lambda^1_{245} $, $\Lambda^2_{245} $, $\Lambda^3_{245} $ are given in  \eqref{a245}.



\bigskip

 Next we  look for further embeddings in the above invariant submodules.
After some analysis we find that the  VM
  $ V^{\Lambda^1_{245} }$ has the type (i) and (v) SVs of weights $\Lambda^2_{245}$, $\Lambda^3_{245}$, 
  thus we have:

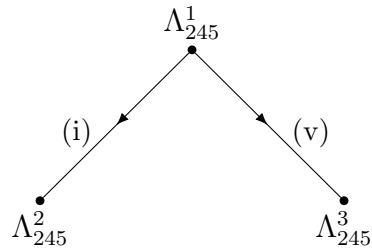
\begin{center}
  \begin{tikzpicture}
    \node[circle,fill=black,scale=0.3] (L2) at (0,10) {};
    \node[circle,fill=black,scale=0.3] (Ld1) at (-2,8) {};
    \node[circle,fill=black,scale=0.3] (Ld5) at (2,8) {};
    \draw[decoration={markings, mark=at position 0.5 with {\arrow{Latex}}},postaction={decorate}] (L2) to (Ld1);
    \node[left,xshift=-2] at ($(L2)!0.55!(Ld1)$) {\small (i)};
    \draw[decoration={markings, mark=at position 0.5 with {\arrow{Latex}}},postaction={decorate}] (L2) to (Ld5);
    \node[right,xshift=2] at ($(L2)!0.55!(Ld5)$) {\small (v)};
    \node[above] at (L2) {$\Lambda^1_{245}  $};
    \node[below] at (Ld1) {$\Lambda^2_{245}  $};
    \node[below] at (Ld5) {$\Lambda^3_{245}  $};
  \end{tikzpicture}
\end{center}


\bigskip

 Further we find that the  VM
  $ V^{\Lambda^3_{245}}$ has no SV.

Next  we find that the  VM
 $ V^{\Lambda^2_{245}} $ has the type (ii) SV of weight $\Lambda^3_{245}$. 

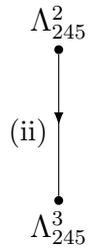
\begin{center}
  \begin{tikzpicture}
    \node[circle,fill=black,scale=0.3] (L5) at (0,10) {};
    \node[circle,fill=black,scale=0.3] (Ld2) at (0,8) {};
    \draw[decoration={markings, mark=at position 0.5 with {\arrow{Latex}}},postaction={decorate}] (L5) to (Ld2);
    \node[left] at ($(L5)!0.55!(Ld2)$) {\small (ii)};
    \node[above] at (L5) {$\Lambda^2_{245} $};
    \node[below] at (Ld2) {$\Lambda^3_{245} $};
  \end{tikzpicture}
\end{center}


Combining the last three  subcases we give the complete diagram for \textbf{A245} :

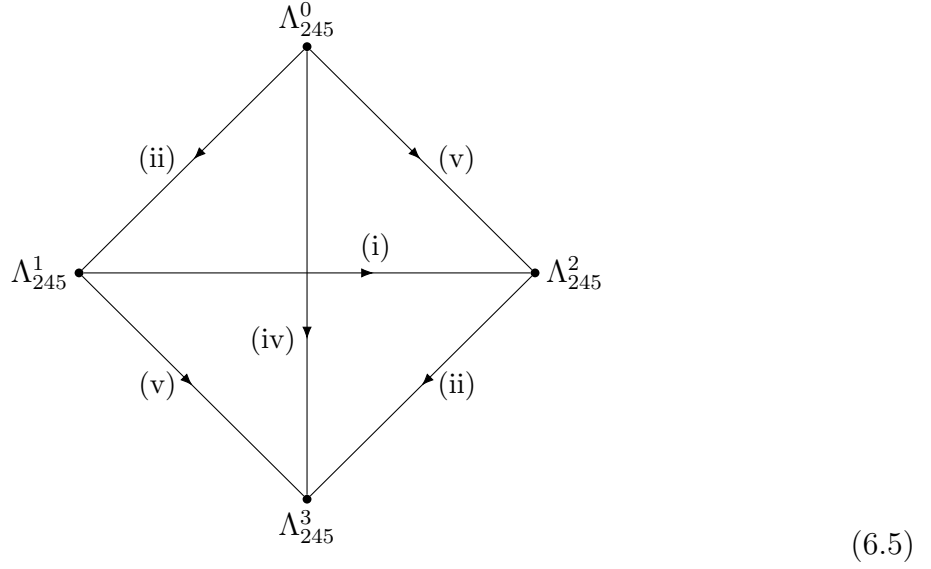
\eqn{A245}
  \begin{tikzpicture}
    \node[circle,fill=black,scale=0.3] (L245) at (0,10) {};
    \node[circle,fill=black,scale=0.3] (Ld2) at (-3,7) {};
    \node[circle,fill=black,scale=0.3] (Ld4) at (0,4) {};
    \node[circle,fill=black,scale=0.3] (Ld5) at (3,7) {};
    \draw[decoration={markings, mark=at position 0.5 with {\arrow{Latex}}},postaction={decorate}] (L245) to (Ld2);
    \node[left,xshift=-2] at ($(L245)!0.5!(Ld2)$) {\small (ii)};
    \draw[decoration={markings, mark=at position 0.65 with {\arrow{Latex}}},postaction={decorate}] (L245) to (Ld4);
    \node[left] at ($(L245)!0.65!(Ld4)$) {\small (iv)};
    \draw[decoration={markings, mark=at position 0.5 with {\arrow{Latex}}},postaction={decorate}] (L245) to (Ld5);
    \node[right,xshift=2] at ($(L245)!0.5!(Ld5)$) {\small (v)};
    \draw[decoration={markings, mark=at position 0.65 with {\arrow{Latex}}},postaction={decorate}] (Ld2) to (Ld5);
    \node[above] at ($(Ld2)!0.65!(Ld5)$) {\small (i)};
    \draw[decoration={markings, mark=at position 0.5 with {\arrow{Latex}}},postaction={decorate}] (Ld2) to (Ld4);
    \node[left,xshift=-2] at ($(Ld2)!0.5!(Ld4)$) {\small (v)};
    \draw[decoration={markings, mark=at position 0.5 with {\arrow{Latex}}},postaction={decorate}] (Ld5) to (Ld4);
    \node[right,xshift=2] at ($(Ld5)!0.5!(Ld4)$) {\small (ii)};
    \node[above] at (L245) {$\Lambda^0_{245}$};
    \node[left] at (Ld2) {$\Lambda^1_{245} $};
    \node[below] at (Ld4) {$\Lambda^3_{245} $};
    \node[right] at (Ld5) {$\Lambda^2_{245}  $};
  \end{tikzpicture}
\ee

\bigskip


\subsection{Case  \textbf{A145}}

In \eqref{L145} we have found the initial embeddings of this case which are given by the diagram:

\bigskip


%
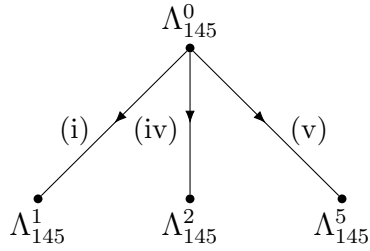
\begin{center}
  \begin{tikzpicture}
    \node[circle,fill=black,scale=0.3] (L145) at (0,10) {};
    \node[circle,fill=black,scale=0.3] (Ld1) at (-2,8) {};
    \node[circle,fill=black,scale=0.3] (Ld4) at (0,8) {};
    \node[circle,fill=black,scale=0.3] (Ld5) at (2,8) {};
    \draw[decoration={markings, mark=at position 0.5 with {\arrow{Latex}}},postaction={decorate}] (L145) to (Ld1);
    \node[left,xshift=-2] at ($(L145)!0.55!(Ld1)$) {\small (i)};
    \draw[decoration={markings, mark=at position 0.5 with {\arrow{Latex}}},postaction={decorate}] (L145) to (Ld4);
    \node[left] at ($(L145)!0.55!(Ld4)$) {\small (iv)};
    \draw[decoration={markings, mark=at position 0.5 with {\arrow{Latex}}},postaction={decorate}] (L145) to (Ld5);
    \node[right,xshift=2] at ($(L145)!0.55!(Ld5)$) {\small (v)};
    \node[above] at (L145) {$\Lambda^0_{145}$};
    \node[below] at (Ld1) {$\Lambda^1_{145} $};
    \node[below] at (Ld4) {$\Lambda^2_{145} $};
    \node[below] at (Ld5) {$\Lambda^5_{145} $};
  \end{tikzpicture}
\end{center}


\bigskip

where $ \L^1_{145}, \L^2_{145}, \L^5_{145} $ are given in \eqref{L145}.
%
%
Next we find the further embeddings in the above invariant submodules.

\medskip
 First we find that the VM
  $\Lambda^1_{145}   $ has the type (iii) (ii) (iv) SVs respectively of weights $\L^5_{145} $ and:
 \eqna{a145na}
 \Lambda^3_{145}   &=&
  \Big( \frac{5}{4} + \hf(p^5-p^4), \frac{3}{4}+\frac{p^5}{2} \Big),
  \\
   \Lambda^4_{145}   &=&
  \Big( \frac{5}{4}+\frac{p^5}{2}, \frac{3}{4} - \hf(p^5-p^4) \Big).
\eena

\medskip  Further we find that the VMs $ \Lambda^3_{145}    $
and $ \Lambda^4_{145}   $ have the type (i) SV of weight $ \L^5_{145}. $

Thus, the embedding picture for $\Lambda^1_{145} $ 
is:

\begin{center}
  \begin{tikzpicture}
    \node[circle,fill=black,scale=0.3] (L1) at (0,10) {};
    \node[circle,fill=black,scale=0.3] (Ld2) at (-2,8) {};
    \node[circle,fill=black,scale=0.3] (Ld3) at (0,6) {};
    \node[circle,fill=black,scale=0.3] (Ld4) at (2,8) {};
    \draw[decoration={markings, mark=at position 0.5 with {\arrow{Latex}}},postaction={decorate}] (L1) to (Ld2);
    \node[left,xshift=-2] at ($(L1)!0.5!(Ld2)$) {\small (ii)};
    \draw[decoration={markings, mark=at position 0.5 with {\arrow{Latex}}},postaction={decorate}] (L1) to (Ld3);
    \node[left] at ($(L1)!0.5!(Ld3)$) {\small (iii)};
    \draw[decoration={markings, mark=at position 0.5 with {\arrow{Latex}}},postaction={decorate}] (L1) to (Ld4);
    \node[right,xshift=2] at ($(L1)!0.5!(Ld4)$) {\small (iv)};
    \draw[decoration={markings, mark=at position 0.5 with {\arrow{Latex}}},postaction={decorate}] (Ld2) to (Ld3);
    \node[left,xshift=-2] at ($(Ld2)!0.5!(Ld3)$) {\small (i)};
    \draw[decoration={markings, mark=at position 0.5 with {\arrow{Latex}}},postaction={decorate}] (Ld4) to (Ld3);
    \node[right,xshift=2] at ($(Ld4)!0.5!(Ld3)$) {\small (ii)};
    \node[above] at (L1) {$\Lambda^1_{145}$};
   \node[left] at (Ld2) {$\Lambda^3_{145}  $};
   \node[below] at (Ld3) {$\Lambda^5_{145} $};
   \node[right] at (Ld4) {$\Lambda^4_{145} $};
  \end{tikzpicture}
\end{center}


\medskip
Next we find that the VM
  $ \Lambda^2_{145} $ has the types (i) and (v) SVs of weights \rf{a145na}{b},   \rf{a145na}{a} (for $p^4 < p^5$), respectively.
 Thus, we have the diagram:

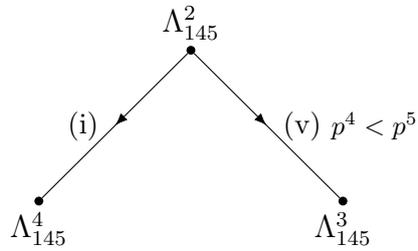
\begin{center}
  \begin{tikzpicture}
    \node[circle,fill=black,scale=0.3] (L4) at (0,10) {};
    \node[circle,fill=black,scale=0.3] (Ld1) at (-2,8) {};
    \node[circle,fill=black,scale=0.3] (Ld5) at (2,8) {};
    \draw[decoration={markings, mark=at position 0.5 with {\arrow{Latex}}},postaction={decorate}] (L4) to (Ld1);
    \node[left,xshift=-2] at ($(L4)!0.5!(Ld1)$) {\small (i)};
    \draw[decoration={markings, mark=at position 0.5 with {\arrow{Latex}}},postaction={decorate}] (L4) to (Ld5);
    \node[right,xshift=2] at ($(L4)!0.5!(Ld5)$) {\small (v)};
    \node[above] at (L4) {$\Lambda^2_{145} $};
    \node[below] at (Ld1) {$\Lambda^4_{145} $};
    \node[below] at (Ld5) {$\Lambda^3_{145} $};
    \node[right,xshift=20] at ($(L4)!0.5!(Ld5)$) {\footnotesize $ p^4 < p^5$};
  \end{tikzpicture}
\end{center}

\medskip
 Finally  we find that the VM
  $\Lambda^5_{145} $ has no SV.

\bigskip


Thus, the complete diagram for \textbf{A145} is:

\eqn{A145dg}
  \begin{tikzpicture}
    \node[circle,fill=black,scale=0.3] (L145) at (0,10) {};
    \node[circle,fill=black,scale=0.3] (Ld1) at (-2,8) {};
    \node[circle,fill=black,scale=0.3] (Ld4) at (2,8) {};
    \node[circle,fill=black,scale=0.3] (Ld5) at (0,3) {};
    \node[circle,fill=black,scale=0.3] (Ldd2) at (-2,5) {};
    \node[circle,fill=black,scale=0.3] (Ldd4) at (2,5) {};
    \draw[decoration={markings, mark=at position 0.5 with {\arrow{Latex}}},postaction={decorate}] (L145) to (Ld1);
    \node[left,xshift=-2] at ($(L145)!0.5!(Ld1)$) {\small (i)};
    \draw[decoration={markings, mark=at position 0.5 with {\arrow{Latex}}},postaction={decorate}] (L145) to (Ld4);
    \node[right,xshift=2] at ($(L145)!0.5!(Ld4)$) {\small (iv)};
    \draw[decoration={markings, mark=at position 0.2 with {\arrow{Latex}}},postaction={decorate}] (L145) to (Ld5);
    \node[left] at ($(L145)!0.2!(Ld5)$) {\small (v)};
    \draw[decoration={markings, mark=at position 0.5 with {\arrow{Latex}}},postaction={decorate}] (Ld1) to (Ldd2);
    \node[left] at ($(Ld1)!0.5!(Ldd2)$) {\small (ii)};
    \draw[decoration={markings, mark=at position 0.5 with {\arrow{Latex}}},postaction={decorate}] (Ld4) to (Ldd4);
    \node[right] at ($(Ld4)!0.5!(Ldd4)$) {\small (i)};
    \draw[decoration={markings, mark=at position 0.5 with {\arrow{Latex}}},postaction={decorate}] (Ldd2) to (Ld5);
    \node[left,xshift=-2] at ($(Ldd2)!0.5!(Ld5)$) {\small (i)};
    \draw[decoration={markings, mark=at position 0.5 with {\arrow{Latex}}},postaction={decorate}] (Ldd4) to (Ld5);
    \node[right,xshift=2] at ($(Ldd4)!0.5!(Ld5)$) {\small (ii)};
    \draw[decoration={markings, mark=at position 0.7 with {\arrow{Latex}}},postaction={decorate}] (Ld1) to (Ld5);
    \node[right] at ($(Ld1)!0.6!(Ld5)$) {\small (iii)};
    \draw[decoration={markings, mark=at position 0.7 with {\arrow{Latex}}},postaction={decorate}] (Ld1) to (Ldd4);
    \node[right,xshift=2] at ($(Ld1)!0.65!(Ldd4)$) {\small (iv)};
    \draw[decoration={markings, mark=at position 0.3 with {\arrow{Latex}}},postaction={decorate}] (Ld4) to (Ldd2);
    \node[left,xshift=-2,yshift=2] at ($(Ld4)!0.2!(Ldd2)$) {\small (v)};
    \node[above] at (L145) {$\Lambda^0_{145}$};
    \node[left] at (Ld1) {$\Lambda^1_{145}  $};
    \node[right] at (Ld4) {$\Lambda^2_{145}  $};
    \node[below] at (Ld5) {$\Lambda^5_{145}  $};
    \node[left] at (Ldd2) {$\Lambda^3_{145}  $};
    \node[right] at (Ldd4) {$\Lambda^4_{145} $};
    \node[above,xshift=-5,yshift=5] at ($(Ld4)!0.2!(Ldd2)$) {\footnotesize $ p^4<p^5$};
  \end{tikzpicture}
\ee

\bigskip


\subsection{Case  \textbf{A14}}

 In \eqref{L14} we have found the initial embeddings of this case which are given by the diagram:

 \bigskip

\begin{center}
  \begin{tikzpicture}
    \node[circle,fill=black,scale=0.3] (L14) at (0,10) {};
    \node[circle,fill=black,scale=0.3] (Ld1) at (-2,8) {};
    \node[circle,fill=black,scale=0.3] (Ld4) at (2,8) {};
    \draw[decoration={markings, mark=at position 0.5 with {\arrow{Latex}}},postaction={decorate}] (L14) to (Ld1);
    \node[left,xshift=-2] at ($(L14)!0.5!(Ld1)$) {\small (i)};
    \draw[decoration={markings, mark=at position 0.5 with {\arrow{Latex}}},postaction={decorate}] (L14) to (Ld4);
    \node[right,xshift=2] at ($(L14)!0.5!(Ld4)$) {\small (iv)};
    \node[above] at (L14) {$ \L^{0}_{14}$};
    \node[below] at (Ld1) {$\Lambda^1_{14} $};
    \node[below] at (Ld4) {$\Lambda^4_{14} $};
    %
  \end{tikzpicture}
\end{center}
see \eqref{A14f}, \eqref{L14} for $\L^a_{14}.$
%
%
%
  Further we distinguish two cases depending whether $ p^1 $ is bigger or smaller than $ \frac{1}{2}+q$ since the embedding diagrams are different.

  \bigskip\noindent
{\bf Case 1:}~~~ $ p^1 < \frac{1}{2}+q$ \quad we have the following embeddings:
\begin{center}
  \begin{tikzpicture}
    \node[circle,fill=black,scale=0.3] (L234) at (0,10) {};
    \node[circle,fill=black,scale=0.3] (Ld2) at (-3,8) {};
    \node[circle,fill=black,scale=0.3] (Ld3) at (-1,8) {};
    \node[circle,fill=black,scale=0.3] (Ld4) at (1,8) {};
    \node[circle,fill=black,scale=0.3] (Ld5) at (3,8) {};
    \draw[decoration={markings, mark=at position 0.5 with {\arrow{Latex}}},postaction={decorate}] (L234) to (Ld2);
    \node[left,xshift=-5] at ($(L234)!0.65!(Ld2)$) {\small (ii)};
    \draw[decoration={markings, mark=at position 0.5 with {\arrow{Latex}}},postaction={decorate}] (L234) to (Ld3);
    \node[left] at ($(L234)!0.65!(Ld3)$) {\small (iii)};
    \draw[decoration={markings, mark=at position 0.5 with {\arrow{Latex}}},postaction={decorate}] (L234) to (Ld4);
    \node[right,xshift=5] at ($(L234)!0.65!(Ld4)$) {\small (iv)};
    \draw[decoration={markings, mark=at position 0.5 with {\arrow{Latex}}},postaction={decorate}] (L234) to (Ld5);
    \node[right,xshift=5] at ($(L234)!0.65!(Ld5)$) {\small (v)};
    \node[above] at (L234) {$\Lambda^1_{14}$};
    \node[below] at (Ld2) {$\Lambda^{2}_{14}$};
    \node[below] at (Ld3) {$\Lambda^{3}_{14}$};
    \node[below] at (Ld4) {$\Lambda^{5}_{14}$};
    \node[below] at (Ld5) {$\Lambda^{6}_{14}$};
    \begin{scope}[xshift=75mm]
        \node[circle,fill=black,scale=0.3] (L14) at (0,10) {};
        \node[circle,fill=black,scale=0.3] (Ld1) at (-2,8) {};
        \node[circle,fill=black,scale=0.3] (Ld4) at (2,8) {};
        \draw[decoration={markings, mark=at position 0.5 with {\arrow{Latex}}},postaction={decorate}] (L14) to (Ld1);
        \node[left,xshift=-5] at ($(L14)!0.6!(Ld1)$) {\small (v)};
        \draw[decoration={markings, mark=at position 0.5 with {\arrow{Latex}}},postaction={decorate}] (L14) to (Ld4);
        \node[right,xshift=5] at ($(L14)!0.6!(Ld4)$) {\small (i)};
        \node[above] at (L14) {$\Lambda^{2}_{14}$};
        \node[below] at (Ld1) {$\Lambda^{4}_{14}$};
        \node[below] at (Ld4) {$\Lambda^{3}_{14}$};
    \end{scope}
  \end{tikzpicture}
\end{center}

\begin{center}
  \begin{tikzpicture}
      \node[circle,fill=black,scale=0.3] (L2) at (0,10) {};
      \node[circle,fill=black,scale=0.3] (L23) at (0,8) {};
      \draw[decoration={markings, mark=at position 0.5 with {\arrow{Latex}}},postaction={decorate}] (L2) to (L23);
      \node[left] at ($(L2)!0.5!(L23)$) {\small (iv)};
      \node[above] at (L2) {$\Lambda^{3}_{14}$};
      \node[below] at (L23) {$\Lambda^5_{14}$};
    \begin{scope}[xshift=50mm]
    \node[circle,fill=black,scale=0.3] (L234) at (0,10) {};
    \node[circle,fill=black,scale=0.3] (Ld2) at (-2,8) {};
    \node[circle,fill=black,scale=0.3] (Ld3) at (0,8) {};
    \node[circle,fill=black,scale=0.3] (Ld4) at (2,8) {};
    \draw[decoration={markings, mark=at position 0.5 with {\arrow{Latex}}},postaction={decorate}] (L234) to (Ld2);
    \node[left,xshift=-5] at ($(L234)!0.6!(Ld2)$) {\small (ii)};
    \draw[decoration={markings, mark=at position 0.5 with {\arrow{Latex}}},postaction={decorate}] (L234) to (Ld3);
    \node[left] at ($(L234)!0.6!(Ld3)$) {\small (iv)};
    \draw[decoration={markings, mark=at position 0.5 with {\arrow{Latex}}},postaction={decorate}] (L234) to (Ld4);
    \node[right,xshift=5] at ($(L234)!0.6!(Ld4)$) {\small (iii)};
    \node[above] at (L234) {$\Lambda^{6}_{14}$};
    \node[below] at (Ld2) {$\Lambda^{4}_{14}$};
    \node[below] at (Ld3) {$\Lambda^{3}_{14}$};
    \node[below] at (Ld4) {$\Lambda^{5}_{14}$};
    \end{scope}
    \begin{scope}[xshift=100mm]
       \node[circle,fill=black,scale=0.3] (L4) at (0,10) {};
      \node[circle,fill=black,scale=0.3] (L14) at (0,8) {};
      \draw[decoration={markings, mark=at position 0.5 with {\arrow{Latex}}},postaction={decorate}] (L4) to (L14);
      \node[left] at ($(L4)!0.5!(L14)$) {\small (i)};
      \node[above] at (L4) {$\Lambda^{4}_{14}$};
      \node[below] at (L14) {$\Lambda^{5}_{14}$};
    \end{scope}
   \end{tikzpicture}
\end{center}
where
\eqna{a14c1}
   \Lambda^{2}_{14} &= \Big( 1+ \hf(p^1-q), 1+\frac{q}{2} \Big),
  \\
  \Lambda^{3}_{14} &= \Big( \frac{3}{2}+\frac{q}{2}, \hf(1+p^1-q) \Big),
  \\
  \Lambda^{5}_{14} &= \Big( \frac{3}{2}+\frac{q}{2},1-\hf(p^1-q) \Big),
  \\
  \Lambda^{6}_{14} &= \Big( \frac{3}{2}- \hf(p^1-q), \hf(1-q) \Big).
 \eena

The complete embedding pattern of Case 1 is:

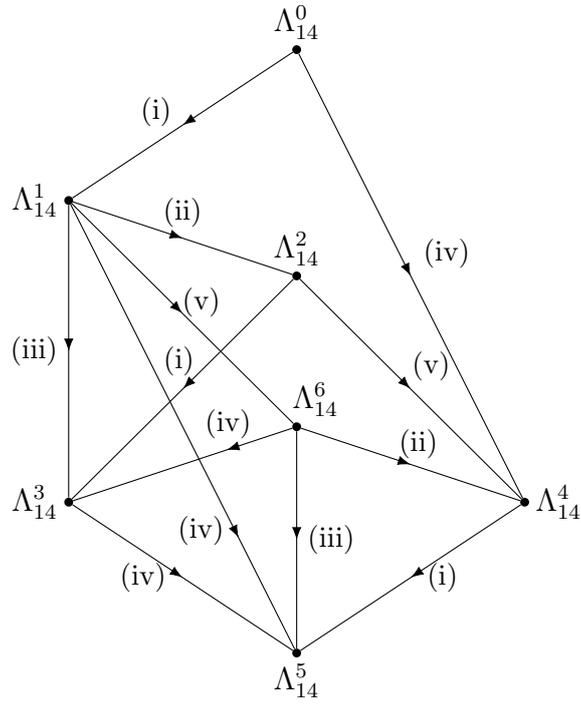
\eqn{A141}
  \begin{tikzpicture}
    \node[circle,fill=black,scale=0.3] (L0) at (0,10) {};
    \node[circle,fill=black,scale=0.3] (L1) at (-3,8) {};
    \node[circle,fill=black,scale=0.3] (L4) at (3,4) {};
    \node[circle,fill=black,scale=0.3] (L13) at (-3,4) {};
    \node[circle,fill=black,scale=0.3] (L12) at (0,7) {};
    \node[circle,fill=black,scale=0.3] (L15) at (0,5) {};
    \node[circle,fill=black,scale=0.3] (L14) at (0,2) {};
    \draw[decoration={markings, mark=at position 0.5 with {\arrow{Latex}}},postaction={decorate}] (L0) to (L1);
    \node[left,yshift=5] at ($(L0)!0.5!(L1)$) {\small (i)};
    \draw[decoration={markings, mark=at position 0.5 with {\arrow{Latex}}},postaction={decorate}] (L0) to (L4);
    \node[right,xshift=5] at ($(L0)!0.45!(L4)$) {\small (iv)};
    \draw[decoration={markings, mark=at position 0.5 with {\arrow{Latex}}},postaction={decorate}] (L1) to (L12);
    \node[above] at ($(L1)!0.5!(L12)$) {\small (ii)};
    \draw[decoration={markings, mark=at position 0.5 with {\arrow{Latex}}},postaction={decorate}] (L1) to (L13);
    \node[left] at ($(L1)!0.5!(L13)$) {\small (iii)};
    \draw[decoration={markings, mark=at position 0.75 with {\arrow{Latex}}},postaction={decorate}] (L1) to (L14);
    \node[left] at ($(L1)!0.73!(L14)$) {\small (iv)};
    \draw[decoration={markings, mark=at position 0.5 with {\arrow{Latex}}},postaction={decorate}] (L1) to (L15);
    \node[right] at ($(L1)!0.45!(L15)$) {\small (v)};
    \draw[decoration={markings, mark=at position 0.5 with {\arrow{Latex}}},postaction={decorate}] (L12) to (L13);
    \node[left,yshift=2] at ($(L12)!0.4!(L13)$) {\small (i)};
    \draw[decoration={markings, mark=at position 0.5 with {\arrow{Latex}}},postaction={decorate}] (L12) to (L4);
    \node[right,xshift=5] at ($(L12)!0.4!(L4)$) {\small (v)};
    \draw[decoration={markings, mark=at position 0.5 with {\arrow{Latex}}},postaction={decorate}] (L13) to (L14);
    \node[left,xshift=-2] at ($(L13)!0.5!(L14)$) {\small (iv)};
    \draw[decoration={markings, mark=at position 0.3 with {\arrow{Latex}}},postaction={decorate}] (L15) to (L13);
    \node[above] at ($(L15)!0.3!(L13)$) {\small (iv)};
    \draw[decoration={markings, mark=at position 0.5 with {\arrow{Latex}}},postaction={decorate}] (L15) to (L14);
    \node[right] at ($(L15)!0.5!(L14)$) {\small (iii)};
    \draw[decoration={markings, mark=at position 0.5 with {\arrow{Latex}}},postaction={decorate}] (L15) to (L4);
    \node[right,yshift=5] at ($(L15)!0.4!(L4)$) {\small (ii)};
    \draw[decoration={markings, mark=at position 0.5 with {\arrow{Latex}}},postaction={decorate}] (L4) to (L14);
    \node[right,xshift=2] at ($(L4)!0.5!(L14)$) {\small (i)};
    \node[above] at (L0) {$\Lambda^0_{14}$};
    \node[left] at (L1) {$\Lambda^{1}_{14}$};
    \node[right] at (L4) {$\Lambda^{4}_{14}$};
    \node[left] at (L13) {$\Lambda^{3}_{14}$};
    \node[above] at (L12) {$\Lambda^{2}_{14}$};
    \node[above,xshift=7] at (L15) {$\Lambda^{6}_{14}$};
    \node[below] at (L14) {$\Lambda^{5}_{14}$};
  \end{tikzpicture}
\ee

\bigskip

{\bf Case 2:}~~~ $  \frac{1}{2}+q < p^1$ \quad we have the following embeddings.
The LWs are the same as Case 1).
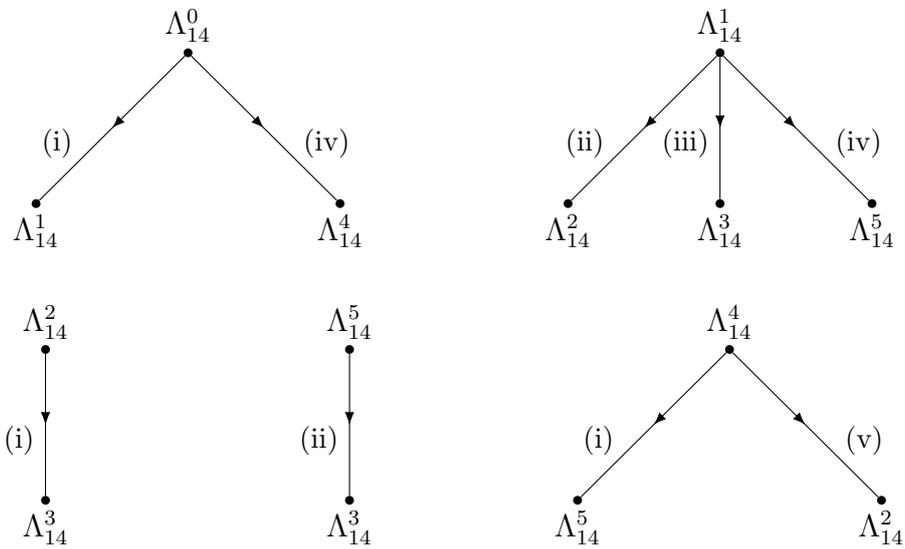
\begin{center}
 \begin{tikzpicture}
    \node[circle,fill=black,scale=0.3] (L14) at (0,10) {};
    \node[circle,fill=black,scale=0.3] (Ld1) at (-2,8) {};
    \node[circle,fill=black,scale=0.3] (Ld4) at (2,8) {};
    \draw[decoration={markings, mark=at position 0.5 with {\arrow{Latex}}},postaction={decorate}] (L14) to (Ld1);
    \node[left,xshift=-5] at ($(L14)!0.6!(Ld1)$) {\small (i)};
    \draw[decoration={markings, mark=at position 0.5 with {\arrow{Latex}}},postaction={decorate}] (L14) to (Ld4);
    \node[right,xshift=5] at ($(L14)!0.6!(Ld4)$) {\small (iv)};
    \node[above] at (L14) {$\Lambda^{0}_{14}$};
    \node[below] at (Ld1) {$\Lambda^{1}_{14}$};
    \node[below] at (Ld4) {$\Lambda^{4}_{14}$};
    \begin{scope}[xshift=70mm]
     \node[circle,fill=black,scale=0.3] (L234) at (0,10) {};
     \node[circle,fill=black,scale=0.3] (Ld2) at (-2,8) {};
     \node[circle,fill=black,scale=0.3] (Ld3) at (0,8) {};
     \node[circle,fill=black,scale=0.3] (Ld4) at (2,8) {};
     \draw[decoration={markings, mark=at position 0.5 with {\arrow{Latex}}},postaction={decorate}] (L234) to (Ld2);
     \node[left,xshift=-5] at ($(L234)!0.6!(Ld2)$) {\small (ii)};
     \draw[decoration={markings, mark=at position 0.5 with {\arrow{Latex}}},postaction={decorate}] (L234) to (Ld3);
     \node[left] at ($(L234)!0.6!(Ld3)$) {\small (iii)};
     \draw[decoration={markings, mark=at position 0.5 with {\arrow{Latex}}},postaction={decorate}] (L234) to (Ld4);
     \node[right,xshift=5] at ($(L234)!0.6!(Ld4)$) {\small (iv)};
     \node[above] at (L234) {$\Lambda^{1}_{14}$};
     \node[below] at (Ld2) {$\Lambda^{2}_{14}$};
     \node[below] at (Ld3) {$\Lambda^{3}_{14}$};
     \node[below] at (Ld4) {$\Lambda^{5}_{14}$};
    \end{scope}
 \end{tikzpicture}
\end{center}

\begin{center}
  \begin{tikzpicture}
    \node[circle,fill=black,scale=0.3] (L12) at (0,10) {};
    \node[circle,fill=black,scale=0.3] (L13) at (0,8) {};
     \draw[decoration={markings, mark=at position 0.5 with {\arrow{Latex}}},postaction={decorate}] (L12) to (L13);
     \node[left] at ($(L12)!0.6!(L13)$) {\small (i)};
     \node[above] at (L12) {$\Lambda^{2}_{14}$};
     \node[below] at (L13) {$\Lambda^{3}_{14}$};
    \begin{scope}[xshift=40mm]
    \node[circle,fill=black,scale=0.3] (L12) at (0,10) {};
    \node[circle,fill=black,scale=0.3] (L13) at (0,8) {};
     \draw[decoration={markings, mark=at position 0.5 with {\arrow{Latex}}},postaction={decorate}] (L12) to (L13);
     \node[left] at ($(L12)!0.6!(L13)$) {\small (ii)};
     \node[above] at (L12) {$\Lambda^{5}_{14}$};
     \node[below] at (L13) {$\Lambda^{3}_{14}$};
     \end{scope}
    \begin{scope}[xshift=90mm]
    \node[circle,fill=black,scale=0.3] (L14) at (0,10) {};
    \node[circle,fill=black,scale=0.3] (Ld1) at (-2,8) {};
    \node[circle,fill=black,scale=0.3] (Ld4) at (2,8) {};
    \draw[decoration={markings, mark=at position 0.5 with {\arrow{Latex}}},postaction={decorate}] (L14) to (Ld1);
    \node[left,xshift=-5] at ($(L14)!0.6!(Ld1)$) {\small (i)};
    \draw[decoration={markings, mark=at position 0.5 with {\arrow{Latex}}},postaction={decorate}] (L14) to (Ld4);
    \node[right,xshift=5] at ($(L14)!0.6!(Ld4)$) {\small (v)};
    \node[above] at (L14) {$\Lambda^{4}_{14}$};
    \node[below] at (Ld1) {$\Lambda^{5}_{14}$};
    \node[below] at (Ld4) {$\Lambda^{2}_{14}$};
    \end{scope}
 \end{tikzpicture}
\end{center}

\noindent
The complete embedding pattern of Case 2 is:

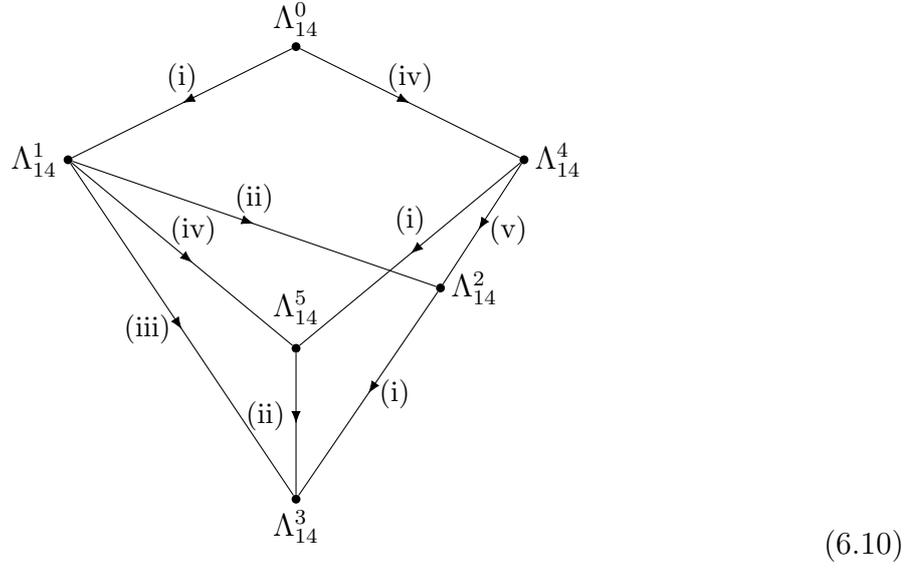
\eqn{A142}
  \begin{tikzpicture}
    \node[circle,fill=black,scale=0.3] (L0) at (0,10) {};
    \node[circle,fill=black,scale=0.3] (L1) at (-3,8.5) {};
    \node[circle,fill=black,scale=0.3] (L4) at (3,8.5) {};
    \node[circle,fill=black,scale=0.3] (L14) at (0,6) {};
    \node[circle,fill=black,scale=0.3] (L12) at (1.9,6.8) {};
    \node[circle,fill=black,scale=0.3] (L13) at (0,4) {};
    \draw[decoration={markings, mark=at position 0.5 with {\arrow{Latex}}},postaction={decorate}] (L0) to (L1);
    \node[above] at ($(L0)!0.5!(L1)$) {\small (i)};
    \draw[decoration={markings, mark=at position 0.5 with {\arrow{Latex}}},postaction={decorate}] (L0) to (L4);
    \node[above] at ($(L0)!0.5!(L4)$) {\small (iv)};
    \draw[decoration={markings, mark=at position 0.5 with {\arrow{Latex}}},postaction={decorate}] (L1) to (L13);
    \node[left] at ($(L1)!0.5!(L13)$) {\small (iii)};
    \draw[decoration={markings, mark=at position 0.55 with {\arrow{Latex}}},postaction={decorate}] (L1) to (L14);
    \node[above,yshift=3] at ($(L1)!0.55!(L14)$) {\small (iv)};
    \draw[decoration={markings, mark=at position 0.5 with {\arrow{Latex}}},postaction={decorate}] (L1) to (L12);
    \node[above] at ($(L1)!0.5!(L12)$) {\small (ii)};
    \draw[decoration={markings, mark=at position 0.5 with {\arrow{Latex}}},postaction={decorate}] (L4) to (L14);
    \node[above,yshift=3] at ($(L4)!0.5!(L14)$) {\small (i)};
    \draw[decoration={markings, mark=at position 0.55 with {\arrow{Latex}}},postaction={decorate}] (L4) to (L12);
    \node[right] at ($(L4)!0.55!(L12)$) {\small (v)};
    \draw[decoration={markings, mark=at position 0.5 with {\arrow{Latex}}},postaction={decorate}] (L14) to (L13);
    \node[left,yshift=3] at ($(L14)!0.5!(L13)$) {\small (ii)};
    \draw[decoration={markings, mark=at position 0.5 with {\arrow{Latex}}},postaction={decorate}] (L12) to (L13);
    \node[right] at ($(L12)!0.5!(L13)$) {\small (i)};
    \node[above] at (L0) {$\Lambda^{0}_{14}$};
    \node[left] at (L1) {$\Lambda^{1}_{14}$};
    \node[right] at (L4) {$\Lambda^{4}_{14}$};
    \node[above,yshift=5] at (L14) {$\Lambda^{5}_{14}$};
    \node[right] at (L12) {$\Lambda^{2}_{14}$};
    \node[below] at (L13) {$\Lambda^{3}_{14}$};
  \end{tikzpicture}
\ee

  \bigskip


\subsection{Case  \textbf{A15}}

 In \eqref{L15} we have found the initial embeddings of this case which are given by the diagram:

 \bigskip

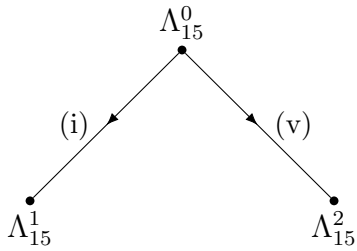
\begin{center}
  \begin{tikzpicture}
    \node[circle,fill=black,scale=0.3] (L15) at (0,10) {};
    \node[circle,fill=black,scale=0.3] (Ld1) at (-2,8) {};
    \node[circle,fill=black,scale=0.3] (Ld4) at (2,8) {};
    \draw[decoration={markings, mark=at position 0.5 with {\arrow{Latex}}},postaction={decorate}] (L15) to (Ld1);
    \node[left,xshift=-2] at ($(L15)!0.5!(Ld1)$) {\small (i)};
    \draw[decoration={markings, mark=at position 0.5 with {\arrow{Latex}}},postaction={decorate}] (L15) to (Ld4);
    \node[right,xshift=2] at ($(L15)!0.5!(Ld4)$) {\small (v)};
    \node[above] at (L15) {$\Lambda^0_{15}$};
    \node[below] at (Ld1) {$\Lambda^1_{15}$};
    \node[below] at (Ld4) {$\Lambda^2_{15}$};
    %
  \end{tikzpicture}
\end{center}
%
%

\bigskip

Next we find the further embeddings in the above invariant submodules.

\medskip
We find that  $\Lambda^1_{15}   $ has the type (ii) SV of weight:
\eqn{a15i}
  \Lambda^3_{15}  = \Big( \frac{5}{4}+\hf(p^1-p^5), \frac{3}{4}+\frac{p^5}{2} \Big).
\end{equation}

Next we find that
  $ \Lambda^2_{15}  $ has the type (i) SV if $2p^5 < p^1$ and its weight is given by \eqref{a15i}.

  We also find that $\Lambda^3_{15}   $ has no SV.


Finally, we find the complete diagram for \textbf{A15}:

\eqn{A15}
  \begin{tikzpicture}
    \node[circle,fill=black,scale=0.3] (L15) at (0,10) {};
    \node[circle,fill=black,scale=0.3] (Ld1) at (-2,8) {};
    \node[circle,fill=black,scale=0.3] (Ldd2) at (0,6) {};
    \node[circle,fill=black,scale=0.3] (Ld5) at (2,8) {};
    \draw[decoration={markings, mark=at position 0.5 with {\arrow{Latex}}},postaction={decorate}] (L15) to (Ld1);
    \node[left,xshift=-2] at ($(L15)!0.5!(Ld1)$) {\small (i)};
    \draw[decoration={markings, mark=at position 0.5 with {\arrow{Latex}}},postaction={decorate}] (L15) to (Ld5);
    \node[right,xshift=2] at ($(L15)!0.5!(Ld5)$) {\small (v)};
    \draw[decoration={markings, mark=at position 0.5 with {\arrow{Latex}}},postaction={decorate}] (Ld1) to (Ldd2);
    \node[left,xshift=-2] at ($(Ld1)!0.5!(Ldd2)$) {\small (ii)};
    \draw[decoration={markings, mark=at position 0.5 with {\arrow{Latex}}},postaction={decorate}] (Ld5) to (Ldd2);
    \node[right,xshift=2] at ($(Ld5)!0.5!(Ldd2)$) {\small (i)};
    \node[above] at (L15) {$\Lambda^0_{15}$};
    \node[left] at (Ld1) {$\Lambda^1_{15}$} ;
    \node[right] at (Ld5) {$\Lambda^2_{15}$} ;
    \node[below] at (Ldd2) {$\Lambda^3_{15}$} ;
  \end{tikzpicture}
\ee

\bigskip


\subsection{Case  \textbf{A24}}

In \eqref{L24} we have found the initial embeddings of this case which are given by the diagram:

 \bigskip

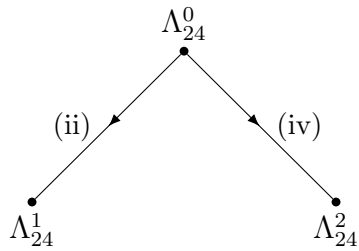
\begin{center}
  \begin{tikzpicture}
    \node[circle,fill=black,scale=0.3] (L24) at (0,10) {};
    \node[circle,fill=black,scale=0.3] (Ld1) at (-2,8) {};
    \node[circle,fill=black,scale=0.3] (Ld4) at (2,8) {};
    \draw[decoration={markings, mark=at position 0.5 with {\arrow{Latex}}},postaction={decorate}] (L24) to (Ld1);
    \node[left,xshift=-2] at ($(L24)!0.5!(Ld1)$) {\small (ii)};
    \draw[decoration={markings, mark=at position 0.5 with {\arrow{Latex}}},postaction={decorate}] (L24) to (Ld4);
    \node[right,xshift=2] at ($(L24)!0.5!(Ld4)$) {\small (iv)};
    \node[above] at (L24) {$\Lambda^0_{24}$};
    \node[below] at (Ld1) {$\Lambda^1_{24}  $};
    \node[below] at (Ld4) {$\Lambda^2_{24}  $};
    %
  \end{tikzpicture}
\end{center}
%
%

\bigskip

Next we find the further embeddings in the above invariant submodules.

We find that
  $ \Lambda^1_{24}  $ has the type (i) SV of weight $\L^2_{24}, $ 
while
  $ \Lambda^1_{24}  $ has no SV.

Thus, the  complete diagram for \textbf{A24} is:

\eqn{A24}
  \begin{tikzpicture}
    \node[circle,fill=black,scale=0.3] (L24) at (0,10) {};
    \node[circle,fill=black,scale=0.3] (Ld2) at (-2,8) {};
    \node[circle,fill=black,scale=0.3] (Ld4) at (0,6) {};
    \draw[decoration={markings, mark=at position 0.5 with {\arrow{Latex}}},postaction={decorate}] (L24) to (Ld2);
    \node[left,xshift=-2] at ($(L24)!0.5!(Ld2)$) {\small (ii)};
    \draw[decoration={markings, mark=at position 0.5 with {\arrow{Latex}}},postaction={decorate}] (L24) to (Ld4);
    \node[right,xshift=2] at ($(L24)!0.5!(Ld4)$) {\small (iv)};
    \draw[decoration={markings, mark=at position 0.5 with {\arrow{Latex}}},postaction={decorate}] (Ld2) to (Ld4);
    \node[left,xshift=-2] at ($(Ld2)!0.5!(Ld4)$) {\small (i)};
    \node[above] at (L24) {$\Lambda^0_{24}$};
    \node[left] at (Ld2) {$\Lambda^1_{24} $};
    \node[below] at (Ld4) {$\Lambda^2_{24} $};
  \end{tikzpicture}
 \ee


\subsection{Cases  \textbf{A1,A2,A4,A5}}

In \S \ref{SEC:EmbPat} we have found the initial embeddings of these cases which are given by the diagram:

\vspace{7mm}

\noindent
\textbf{A1, A2, A4, A5}

\begin{center}
  \begin{tikzpicture}
    \node[circle,fill=black,scale=0.3] (L1) at (0,0) {};
    \node[circle,fill=black,scale=0.3] (Ld1) at (3,0) {};
    \draw[decoration={markings, mark=at position 0.5 with {\arrow{Latex}}},postaction={decorate}] (L1) to (Ld1);
    \node[above] at ($(L1)!0.5!(Ld1)$) {\small (a)};
     \node[left] at (L1) {$\Lambda^0_a$};
     \node[right] at (Ld1) {$\Lambda^1_a$};
     \node[right,xshift=20mm] at (Ld1) {$ a = 1,2,4,5$};
  \end{tikzpicture}
\end{center}

\bigskip

Next we find the further embeddings in the above invariant submodules.

\subsubsection{Embedding diagram of $ V^{\Lambda^0_1} $ (case A1)}

We first recall  from \S \ref{SEC:EmbPat} the weights in the initial embedding diagram above
   for this case:
  \eqna{a1ini}
   \L^0_1 &=&    \Big( \Lambda_1, \Lambda_1+ \hf (p^1-1) \Big),
     \quad
     \Lambda_1 \neq \frac{1}{4}(5-p^1-p^4) \\
    \Lambda^1_1  &=& \Big( \Lambda_1 + \frac{p^1}{2},  \Lambda_1 -\hf \Big). \eena

We draw the embedding diagram by distinguishing the four cases.

\bigskip\noindent
Case 1) $ 5 \leq 4 \Lambda_1 $

\medskip
$ \L^1_1 $ is irreducible.

\bigskip\noindent
Case 2) $ 4 \Lambda_1 < 5 -2p^1 $

\medskip
We have the following embedding patterns:
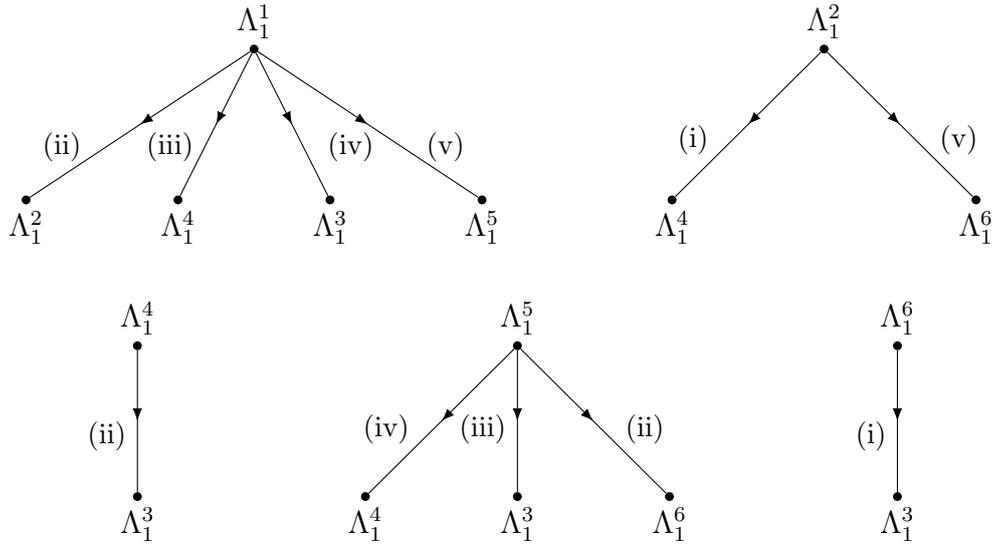
\begin{center}
  \begin{tikzpicture}
    \node[circle,fill=black,scale=0.3] (L234) at (0,10) {};
    \node[circle,fill=black,scale=0.3] (Ld2) at (-3,8) {};
    \node[circle,fill=black,scale=0.3] (Ld3) at (-1,8) {};
    \node[circle,fill=black,scale=0.3] (Ld4) at (1,8) {};
    \node[circle,fill=black,scale=0.3] (Ld5) at (3,8) {};
    \draw[decoration={markings, mark=at position 0.5 with {\arrow{Latex}}},postaction={decorate}] (L234) to (Ld2);
    \node[left,xshift=-5] at ($(L234)!0.65!(Ld2)$) {\small (ii)};
    \draw[decoration={markings, mark=at position 0.5 with {\arrow{Latex}}},postaction={decorate}] (L234) to (Ld3);
    \node[left] at ($(L234)!0.65!(Ld3)$) {\small (iii)};
    \draw[decoration={markings, mark=at position 0.5 with {\arrow{Latex}}},postaction={decorate}] (L234) to (Ld4);
    \node[right,xshift=5] at ($(L234)!0.65!(Ld4)$) {\small (iv)};
    \draw[decoration={markings, mark=at position 0.5 with {\arrow{Latex}}},postaction={decorate}] (L234) to (Ld5);
    \node[right,xshift=5] at ($(L234)!0.65!(Ld5)$) {\small (v)};
    \node[above] at (L234) {$\L^1_1$};
    \node[below] at (Ld2) {$\Lambda^{2}_1$};
    \node[below] at (Ld3) {$\Lambda^{4}_1$};
    \node[below] at (Ld4) {$\Lambda^{3}_1$};
    \node[below] at (Ld5) {$\Lambda^{5}_1$};
    \begin{scope}[xshift=75mm]
        \node[circle,fill=black,scale=0.3] (L14) at (0,10) {};
        \node[circle,fill=black,scale=0.3] (Ld1) at (-2,8) {};
        \node[circle,fill=black,scale=0.3] (Ld4) at (2,8) {};
        \draw[decoration={markings, mark=at position 0.5 with {\arrow{Latex}}},postaction={decorate}] (L14) to (Ld1);
        \node[left,xshift=-5] at ($(L14)!0.6!(Ld1)$) {\small (i)};
        \draw[decoration={markings, mark=at position 0.5 with {\arrow{Latex}}},postaction={decorate}] (L14) to (Ld4);
        \node[right,xshift=5] at ($(L14)!0.6!(Ld4)$) {\small (v)};
        \node[above] at (L14) {$\Lambda^{2}_1$};
        \node[below] at (Ld1) {$\Lambda^{4}_1$};
        \node[below] at (Ld4) {$\Lambda^{6}_1$};
    \end{scope}
  \end{tikzpicture}
\end{center}

\begin{center}
  \begin{tikzpicture}
    \node[circle,fill=black,scale=0.3] (L12) at (0,10) {};
    \node[circle,fill=black,scale=0.3] (L13) at (0,8) {};
     \draw[decoration={markings, mark=at position 0.5 with {\arrow{Latex}}},postaction={decorate}] (L12) to (L13);
     \node[left] at ($(L12)!0.6!(L13)$) {\small (ii)};
     \node[above] at (L12) {$\Lambda^{4}_1$};
     \node[below] at (L13) {$\Lambda^{3}_1$};
    \begin{scope}[xshift=50mm]
    \node[circle,fill=black,scale=0.3] (L14) at (0,10) {};
    \node[circle,fill=black,scale=0.3] (Ld1) at (-2,8) {};
    \node[circle,fill=black,scale=0.3] (Ld14) at (0,8) {};
    \node[circle,fill=black,scale=0.3] (Ld4) at (2,8) {};
    \draw[decoration={markings, mark=at position 0.5 with {\arrow{Latex}}},postaction={decorate}] (L14) to (Ld1);
    \node[left,xshift=-5] at ($(L14)!0.55!(Ld1)$) {\small (iv)};
    \draw[decoration={markings, mark=at position 0.5 with {\arrow{Latex}}},postaction={decorate}] (L14) to (Ld14);
    \node[left] at ($(L14)!0.55!(Ld14)$) {\small (iii)};
    \draw[decoration={markings, mark=at position 0.5 with {\arrow{Latex}}},postaction={decorate}] (L14) to (Ld4);
    \node[right,xshift=5] at ($(L14)!0.55!(Ld4)$) {\small (ii)};
    \node[above] at (L14) {$\Lambda^{5}_1$};
    \node[below] at (Ld1) {$\Lambda^{4}_1$};
    \node[below] at (Ld14) {$\Lambda^{3}_1$};
    \node[below] at (Ld4) {$\Lambda^{6}_1$};
    \end{scope}
    \begin{scope}[xshift=100mm]
    \node[circle,fill=black,scale=0.3] (L12) at (0,10) {};
    \node[circle,fill=black,scale=0.3] (L13) at (0,8) {};
     \draw[decoration={markings, mark=at position 0.5 with {\arrow{Latex}}},postaction={decorate}] (L12) to (L13);
     \node[left] at ($(L12)!0.6!(L13)$) {\small (i)};
     \node[above] at (L12) {$\Lambda^{6}_1$};
     \node[below] at (L13) {$\Lambda^{3}_1$};
     \end{scope}

 \end{tikzpicture}
\end{center}
where
 \eqna{a1case2}
      \Lambda^{2}_1 &= \Big( \Lambda_1+ \frac{p^1}{2}, -\Lambda_1+2 \Big),
      \\
      \Lambda^{4}_1 &= \Big( -\Lambda_1+ \frac{5}{2}, \Lambda_1+ \hf (p^1-1)\Big),
      \\
      \Lambda^{3}_1 &= \Big( -\Lambda_1+ \frac{5}{2}, -\Lambda_1-\frac{p^1}{2}+2 \Big),
      \\
      \Lambda^{5}_1 &= \Big( -\Lambda_1- \hf (p^1-5), \Lambda_1-\frac{1}{2} \Big),
      \\
      \Lambda^{6}_1 &=
      \Big(  -\Lambda_1 - \hf (p^1-5), -\Lambda_1+2 \Big).
 \eena
Thus the complete pattern for A1 Case 2 is given as follows:

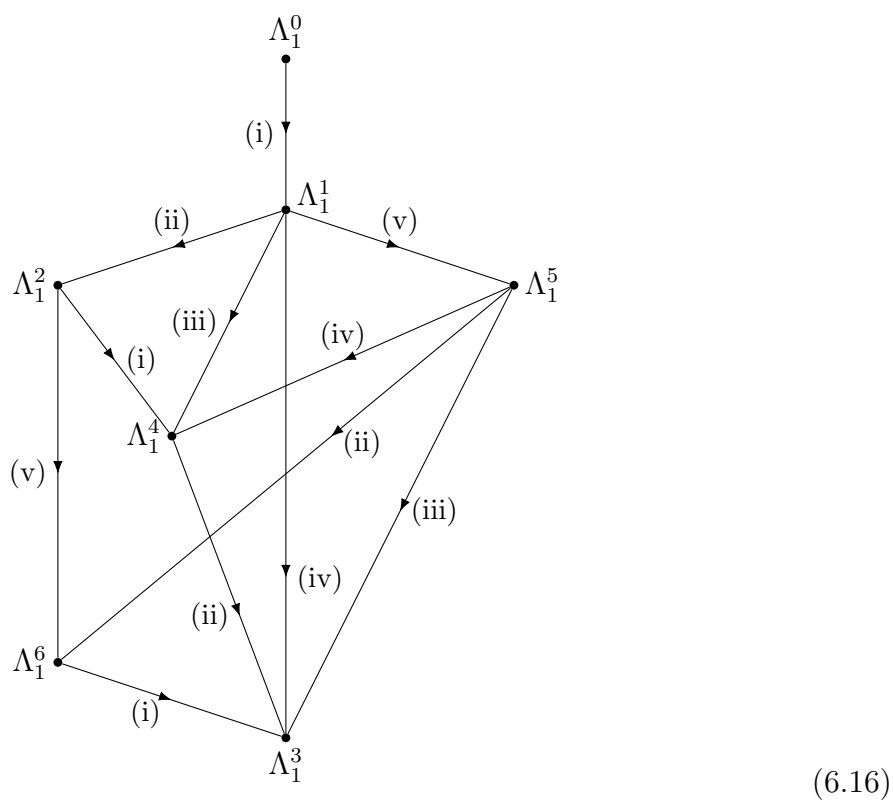
\eqn{A1C2}
  \begin{tikzpicture}
    \node[circle,fill=black,scale=0.3] (L0) at (0,10) {};
    \node[circle,fill=black,scale=0.3] (L1) at (0,8) {};
    \node[circle,fill=black,scale=0.3] (L12) at (-3,7) {};
    \node[circle,fill=black,scale=0.3] (L15) at (3,7) {};
    \node[circle,fill=black,scale=0.3] (L13) at (-1.5,5) {};
    \node[circle,fill=black,scale=0.3] (L125) at (-3,2) {};
    \node[circle,fill=black,scale=0.3] (L14) at (0,1) {};
    \draw[decoration={markings, mark=at position 0.5 with {\arrow{Latex}}},postaction={decorate}] (L0) to (L1);
    \node[left] at ($(L0)!0.5!(L1)$) {\small (i)};
    \draw[decoration={markings, mark=at position 0.5 with {\arrow{Latex}}},postaction={decorate}] (L1) to (L12);
    \node[above] at ($(L1)!0.5!(L12)$) {\small (ii)};
    \draw[decoration={markings, mark=at position 0.5 with {\arrow{Latex}}},postaction={decorate}] (L1) to (L13);
    \node[left] at ($(L1)!0.5!(L13)$) {\small (iii)};
    \draw[decoration={markings, mark=at position 0.7 with {\arrow{Latex}}},postaction={decorate}] (L1) to (L14);
    \node[right] at ($(L1)!0.7!(L14)$) {\small (iv)};
    \draw[decoration={markings, mark=at position 0.5 with {\arrow{Latex}}},postaction={decorate}] (L1) to (L15);
    \node[above] at ($(L1)!0.5!(L15)$) {\small (v)};
    \draw[decoration={markings, mark=at position 0.5 with {\arrow{Latex}}},postaction={decorate}] (L12) to (L125);
    \node[left] at ($(L12)!0.5!(L125)$) {\small (v)};
    \draw[decoration={markings, mark=at position 0.5 with {\arrow{Latex}}},postaction={decorate}] (L12) to (L13);
    \node[right] at ($(L12)!0.5!(L13)$) {\small (i)};
    \draw[decoration={markings, mark=at position 0.6 with {\arrow{Latex}}},postaction={decorate}] (L13) to (L14);
    \node[left] at ($(L13)!0.6!(L14)$) {\small (ii)};
    \draw[decoration={markings, mark=at position 0.4 with {\arrow{Latex}}},postaction={decorate}] (L15) to (L125);
    \node[right,yshift=-3] at ($(L15)!0.4!(L125)$) {\small (ii)};
    \draw[decoration={markings, mark=at position 0.5 with {\arrow{Latex}}},postaction={decorate}] (L15) to (L13);
    \node[above] at ($(L15)!0.5!(L13)$) {\small (iv)};
    \draw[decoration={markings, mark=at position 0.5 with {\arrow{Latex}}},postaction={decorate}] (L15) to (L14);
    \node[right] at ($(L15)!0.5!(L14)$) {\small (iii)};
    \draw[decoration={markings, mark=at position 0.5 with {\arrow{Latex}}},postaction={decorate}] (L125) to (L14);
    \node[left,yshift=-5] at ($(L125)!0.5!(L14)$) {\small (i)};
    \node[above] at (L0) {$\Lambda^0_1$};
    \node[right,yshift=5] at (L1) {$\L^1_1$};
    \node[left] at (L12) {$\Lambda^{2}_1$};
    \node[right] at (L15) {$\Lambda^{5}_1$};
    \node[left] at (L13) {$\Lambda^{4}_1$};
    \node[left] at (L125) {$\Lambda^{6}_1$};
    \node[below] at (L14) {$\Lambda^{3}_1$};

  \end{tikzpicture}
\ee

\bigskip\noindent
Case 3) $ 5-2p^1 < 4 \Lambda_1 < 5-p^1$

\medskip
In this case, the embedding pattern is simplified:
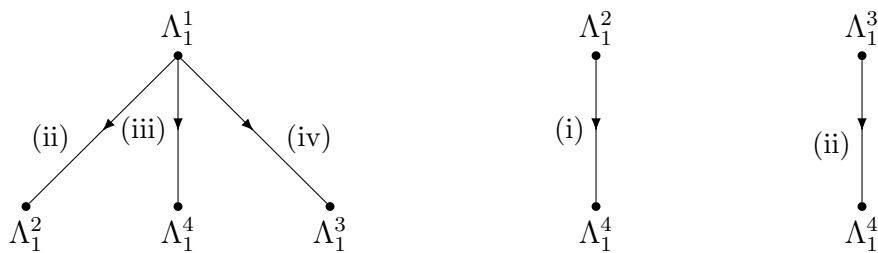
\begin{center}
  \begin{tikzpicture}
    \node[circle,fill=black,scale=0.3] (L14) at (0,10) {};
    \node[circle,fill=black,scale=0.3] (Ld1) at (-2,8) {};
    \node[circle,fill=black,scale=0.3] (Ld14) at (0,8) {};
    \node[circle,fill=black,scale=0.3] (Ld4) at (2,8) {};
    \draw[decoration={markings, mark=at position 0.5 with {\arrow{Latex}}},postaction={decorate}] (L14) to (Ld1);
    \node[left,xshift=-5] at ($(L14)!0.55!(Ld1)$) {\small (ii)};
    \draw[decoration={markings, mark=at position 0.5 with {\arrow{Latex}}},postaction={decorate}] (L14) to (Ld14);
    \node[left] at ($(L14)!0.5!(Ld14)$) {\small (iii)};
    \draw[decoration={markings, mark=at position 0.5 with {\arrow{Latex}}},postaction={decorate}] (L14) to (Ld4);
    \node[right,xshift=5] at ($(L14)!0.55!(Ld4)$) {\small (iv)};
    \node[above] at (L14) {$\Lambda^{1}_1$};
    \node[below] at (Ld1) {$\Lambda^{2}_1$};
    \node[below] at (Ld14) {$\Lambda^{4}_1$};
    \node[below] at (Ld4) {$\Lambda^{3}_1$};
   \begin{scope}[xshift=55mm]
    \node[circle,fill=black,scale=0.3] (L12) at (0,10) {};
    \node[circle,fill=black,scale=0.3] (L13) at (0,8) {};
     \draw[decoration={markings, mark=at position 0.5 with {\arrow{Latex}}},postaction={decorate}] (L12) to (L13);
     \node[left] at ($(L12)!0.5!(L13)$) {\small (i)};
     \node[above] at (L12) {$\Lambda^{2}_1$};
     \node[below] at (L13) {$\Lambda^{4}_1$};
    \end{scope}

    \begin{scope}[xshift=90mm]
    \node[circle,fill=black,scale=0.3] (L12) at (0,10) {};
    \node[circle,fill=black,scale=0.3] (L13) at (0,8) {};
     \draw[decoration={markings, mark=at position 0.5 with {\arrow{Latex}}},postaction={decorate}] (L12) to (L13);
     \node[left] at ($(L12)!0.6!(L13)$) {\small (ii)};
     \node[above] at (L12) {$\Lambda^{3}_1$};
     \node[below] at (L13) {$\Lambda^{4}_1$};
     \end{scope}

 \end{tikzpicture}
\end{center}
Then the complete embedding pattern A1 Case 3 is given by:

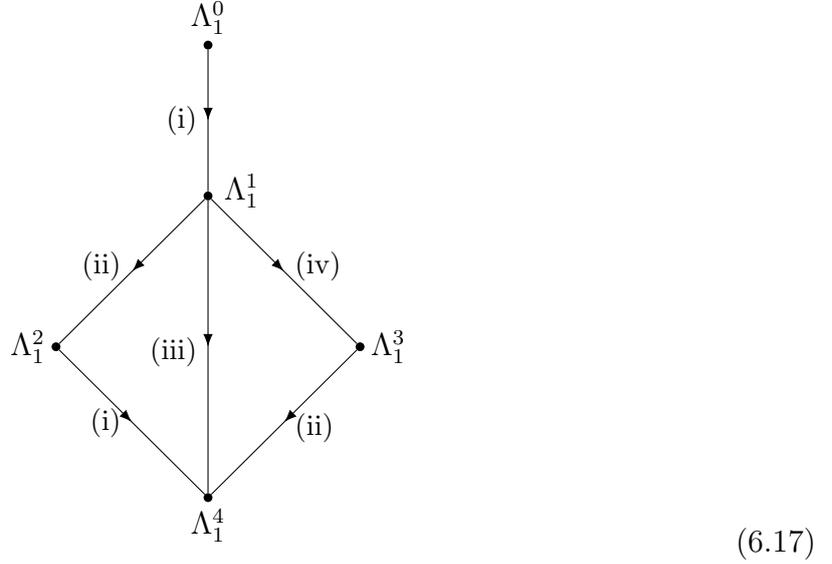
\eqn{A1C3}
  \begin{tikzpicture}
    \node[circle,fill=black,scale=0.3] (L0) at (0,12) {};
    \node[circle,fill=black,scale=0.3] (L1) at (0,10) {};
    \node[circle,fill=black,scale=0.3] (L12) at (-2,8) {};
    \node[circle,fill=black,scale=0.3] (L14) at (2,8) {};
    \node[circle,fill=black,scale=0.3] (L13) at (0,6) {};
    \draw[decoration={markings, mark=at position 0.5 with {\arrow{Latex}}},postaction={decorate}] (L0) to (L1);
    \node[left] at ($(L0)!0.5!(L1)$) {\small (i)};
    \draw[decoration={markings, mark=at position 0.5 with {\arrow{Latex}}},postaction={decorate}] (L1) to (L12);
    \node[left,yshift=2] at ($(L1)!0.5!(L12)$) {\small (ii)};
    \draw[decoration={markings, mark=at position 0.5 with {\arrow{Latex}}},postaction={decorate}] (L1) to (L14);
    \node[right,yshift=2] at ($(L1)!0.5!(L14)$) {\small (iv)};
    \draw[decoration={markings, mark=at position 0.5 with {\arrow{Latex}}},postaction={decorate}] (L1) to (L13);
    \node[left,yshift=-2] at ($(L1)!0.5!(L13)$) {\small (iii)};
    \draw[decoration={markings, mark=at position 0.5 with {\arrow{Latex}}},postaction={decorate}] (L12) to (L13);
    \node[left] at ($(L12)!0.5!(L13)$) {\small (i)};
    \draw[decoration={markings, mark=at position 0.5 with {\arrow{Latex}}},postaction={decorate}] (L14) to (L13);
    \node[right,yshift=-2] at ($(L14)!0.5!(L13)$) {\small (ii)};
    \node[above] at (L0) {$\Lambda^{0}_1$};
    \node[right,xshift=2,yshift=2] at (L1) {$\Lambda^{1}_1$};
    \node[left] at (L12) {$\Lambda^{2}_1$};
    \node[right] at (L14) {$\Lambda^{3}_1$};
    \node[below] at (L13) {$\Lambda^{4}_1$};
 \end{tikzpicture}
\ee

\bigskip\noindent
Case 4) $ 5-p^1 < 4 \Lambda_1 < 5$

\medskip
In A1 Case 4 the embedding pattern is very simple:

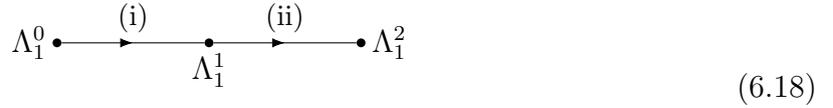
\eqn{A1C4}
  \begin{tikzpicture}
    \node[circle,fill=black,scale=0.3] (L0) at (0,0) {};
    \node[circle,fill=black,scale=0.3] (L1) at (2,0) {};
    \node[circle,fill=black,scale=0.3] (L12) at (4,0) {};
    \draw[decoration={markings, mark=at position 0.5 with {\arrow{Latex}}},postaction={decorate}] (L0) to (L1);
    \node[above] at ($(L0)!0.5!(L1)$) {\small (i)};
    \draw[decoration={markings, mark=at position 0.5 with {\arrow{Latex}}},postaction={decorate}] (L1) to (L12);
    \node[above] at ($(L1)!0.5!(L12)$) {\small (ii)};

   \node[left] at (L0) {$\Lambda^{0}_1$};
   \node[below] at (L1) {$\Lambda^{1}_1$};
   \node[right] at (L12) {$\Lambda^{2}_1$};
  \end{tikzpicture}
\ee
Note that the diagram of Case 4) may be considered as part of the diagram of case 3) taking into account only the
VMs ~$\L^0_1$, ~$ \L^1_1$, ~$\Lambda^{2}_1$ and the two embeddings connecting them, also noting in for case 4) the VM
~$\Lambda^{2}_1$ is irreducible.


\subsubsection{Embedding diagram of $ V^{\Lambda^0_2} $ (case A2)}

We first recall  from  \eqref{L2} the weights in the initial embedding diagram above
   for this case:
  \eqna{a2ini}
     \Lambda^0_2 &=&  \Big( \Lambda_1, \frac{3}{4}-\frac{p^2}{2} \Big),
      \\
    && -2\Lambda_1 + \frac{5}{2} + p^2 \notin \bbn \nn\\
    \Lambda^1_2  &= &\Big( \Lambda_1, \frac{3}{4}+\frac{p^2}{2} \Big).
    \eena

   We draw the embedding diagrams by distinguishing the four cases.

\bigskip\noindent
Case 1) $ 5+2p^2 \leq 4 \Lambda_1 $

\medskip
$\Lambda^1_2 = \Lambda'^{\mathrm{(ii)}}$ is irreducible.

\bigskip\noindent
Case 2) $  4 \Lambda_1 < 5-2p^2 $

\medskip
We have the following patterns:
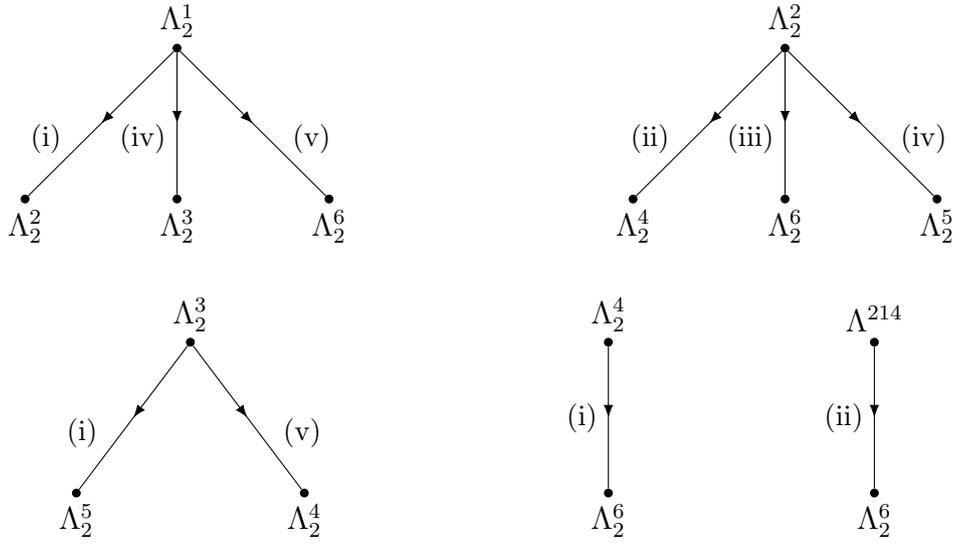
\begin{center}
  \begin{tikzpicture}
    \node[circle,fill=black,scale=0.3] (L234) at (0,10) {};
    \node[circle,fill=black,scale=0.3] (Ld2) at (-2,8) {};
    \node[circle,fill=black,scale=0.3] (Ld3) at (0,8) {};
    \node[circle,fill=black,scale=0.3] (Ld4) at (2,8) {};
    \draw[decoration={markings, mark=at position 0.5 with {\arrow{Latex}}},postaction={decorate}] (L234) to (Ld2);
    \node[left,xshift=-5] at ($(L234)!0.6!(Ld2)$) {\small (i)};
    \draw[decoration={markings, mark=at position 0.5 with {\arrow{Latex}}},postaction={decorate}] (L234) to (Ld3);
    \node[left] at ($(L234)!0.6!(Ld3)$) {\small (iv)};
    \draw[decoration={markings, mark=at position 0.5 with {\arrow{Latex}}},postaction={decorate}] (L234) to (Ld4);
    \node[right,xshift=5] at ($(L234)!0.6!(Ld4)$) {\small (v)};
    \node[above] at (L234) {$\Lambda^1_2$};
    \node[below] at (Ld2) {$\Lambda^{2}_2$};
    \node[below] at (Ld3) {$\Lambda^{3}_2$};
    \node[below] at (Ld4) {$\Lambda^{6}_2$};
    \begin{scope}[xshift=80mm]
    \node[circle,fill=black,scale=0.3] (L234) at (0,10) {};
    \node[circle,fill=black,scale=0.3] (Ld2) at (-2,8) {};
    \node[circle,fill=black,scale=0.3] (Ld3) at (0,8) {};
    \node[circle,fill=black,scale=0.3] (Ld4) at (2,8) {};
    \draw[decoration={markings, mark=at position 0.5 with {\arrow{Latex}}},postaction={decorate}] (L234) to (Ld2);
    \node[left,xshift=-5] at ($(L234)!0.6!(Ld2)$) {\small (ii)};
    \draw[decoration={markings, mark=at position 0.5 with {\arrow{Latex}}},postaction={decorate}] (L234) to (Ld3);
    \node[left] at ($(L234)!0.6!(Ld3)$) {\small (iii)};
    \draw[decoration={markings, mark=at position 0.5 with {\arrow{Latex}}},postaction={decorate}] (L234) to (Ld4);
    \node[right,xshift=5] at ($(L234)!0.6!(Ld4)$) {\small (iv)};
    \node[above] at (L234) {$\Lambda^{2}_2$};
    \node[below] at (Ld2) {$\Lambda^{4}_2$};
    \node[below] at (Ld3) {$\Lambda^{6}_2$};
    \node[below] at (Ld4) {$\Lambda^{5}_2$};
    \end{scope}

   \end{tikzpicture}
\end{center}

\begin{center}
  \begin{tikzpicture}
    \node[circle,fill=black,scale=0.3] (L234) at (0,10) {};
    \node[circle,fill=black,scale=0.3] (Ld2) at (-1.5,8) {};
    \node[circle,fill=black,scale=0.3] (Ld4) at (1.5,8) {};
    \draw[decoration={markings, mark=at position 0.5 with {\arrow{Latex}}},postaction={decorate}] (L234) to (Ld2);
    \node[left,xshift=-5] at ($(L234)!0.6!(Ld2)$) {\small (i)};
    \draw[decoration={markings, mark=at position 0.5 with {\arrow{Latex}}},postaction={decorate}] (L234) to (Ld4);
    \node[right,xshift=5] at ($(L234)!0.6!(Ld4)$) {\small (v)};
    \node[above] at (L234) {$\Lambda^{3}_2$};
    \node[below] at (Ld2) {$\Lambda^{5}_2$};
    \node[below] at (Ld4) {$\Lambda^{4}_2$};
    \begin{scope}[xshift=55mm]
      \node[circle,fill=black,scale=0.3] (L2) at (0,10) {};
      \node[circle,fill=black,scale=0.3] (L23) at (0,8) {};
      \draw[decoration={markings, mark=at position 0.5 with {\arrow{Latex}}},postaction={decorate}] (L2) to (L23);
      \node[left] at ($(L2)!0.5!(L23)$) {\small (i)};
      \node[above] at (L2) {$\Lambda^{4}_2$};
      \node[below] at (L23) {$\Lambda^{6}_2$};
    \end{scope}
    \begin{scope}[xshift=90mm]
      \node[circle,fill=black,scale=0.3] (L2) at (0,10) {};
      \node[circle,fill=black,scale=0.3] (L23) at (0,8) {};
      \draw[decoration={markings, mark=at position 0.5 with {\arrow{Latex}}},postaction={decorate}] (L2) to (L23);
      \node[left] at ($(L2)!0.5!(L23)$) {\small (ii)};
      \node[above] at (L2) {$\Lambda^{214}$};
      \node[below] at (L23) {$\Lambda^{6}_2$};
    \end{scope}
   \end{tikzpicture}
\end{center}
\begin{align}
     \Lambda^{2}_2   &=
         \Big( \frac{5}{4}+\frac{p^2}{2}, \Lambda_1-\hf \Big),
     \\
     \Lambda^{3}_2  &=
         \Big( \frac{5}{4}-\frac{p^2}{2}, -\Lambda_1+2 \Big),
     \\
     \Lambda^{6}_2  &=
         \Big( -\Lambda_1+ \frac{5}{2}, \frac{3}{4}+\frac{p^2}{2} \Big),
     \\
     \Lambda^{4}_2 &= \Big( \frac{5}{4}+\frac{p^2}{2}, -\Lambda_1+2 \Big),
     \\
     \Lambda^{5}_2 &= \Big( -\Lambda_1+\frac{5}{2}, \frac{3}{4}-\frac{p^2}{2} \Big).
\end{align}
Thus the complete embedding pattern of A2 Case 2 is given by:

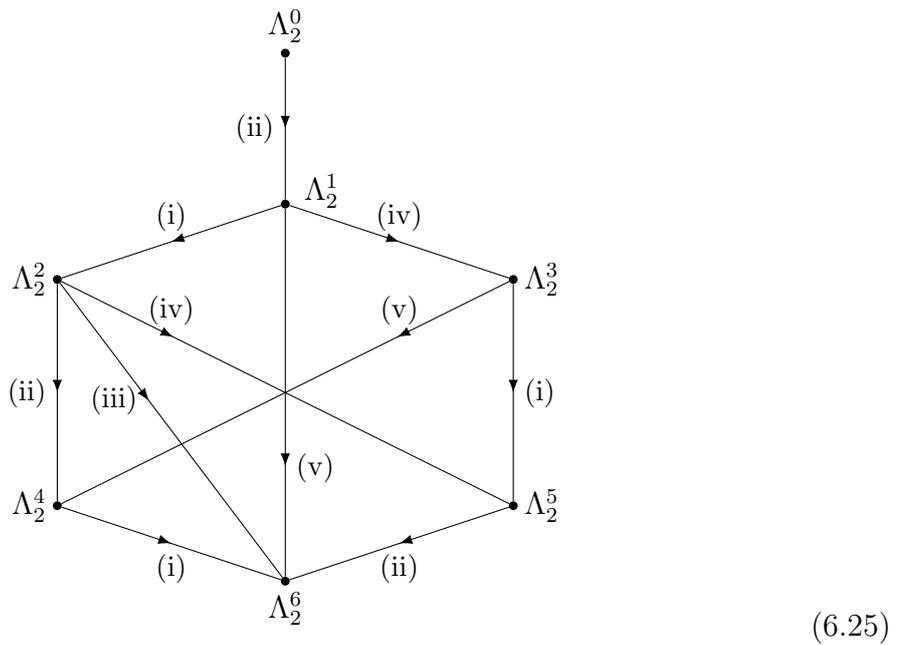
\eqn{A2C2}
  \begin{tikzpicture}
    \node[circle,fill=black,scale=0.3] (L0) at (0,10) {};
    \node[circle,fill=black,scale=0.3] (L2) at (0,8) {};
    \node[circle,fill=black,scale=0.3] (L21) at (-3,7) {};
    \node[circle,fill=black,scale=0.3] (L24) at (3,7) {};
    \node[circle,fill=black,scale=0.3] (L212) at (-3,4) {};
    \node[circle,fill=black,scale=0.3] (L214) at (3,4) {};
    \node[circle,fill=black,scale=0.3] (L25) at (0,3) {};
    \draw[decoration={markings, mark=at position 0.5 with {\arrow{Latex}}},postaction={decorate}] (L0) to (L2);
    \node[left] at ($(L0)!0.5!(L2)$) {\small (ii)};
    \draw[decoration={markings, mark=at position 0.5 with {\arrow{Latex}}},postaction={decorate}] (L2) to (L21);
    \node[above] at ($(L2)!0.5!(L21)$) {\small (i)};
    \draw[decoration={markings, mark=at position 0.5 with {\arrow{Latex}}},postaction={decorate}] (L2) to (L24);
    \node[above] at ($(L2)!0.5!(L24)$) {\small (iv)};
    \draw[decoration={markings, mark=at position 0.7 with {\arrow{Latex}}},postaction={decorate}] (L2) to (L25);
    \node[right] at ($(L2)!0.7!(L25)$) {\small (v)};
    \draw[decoration={markings, mark=at position 0.5 with {\arrow{Latex}}},postaction={decorate}] (L21) to (L212);
    \node[left] at ($(L21)!0.5!(L212)$) {\small (ii)};
    \draw[decoration={markings, mark=at position 0.4 with {\arrow{Latex}}},postaction={decorate}] (L21) to (L25);
    \node[left] at ($(L21)!0.4!(L25)$) {\small (iii)};
    \draw[decoration={markings, mark=at position 0.25 with {\arrow{Latex}}},postaction={decorate}] (L21) to (L214);
    \node[above] at ($(L21)!0.25!(L214)$) {\small (iv)};
    \draw[decoration={markings, mark=at position 0.25 with {\arrow{Latex}}},postaction={decorate}] (L24) to (L212);
    \node[above] at ($(L24)!0.25!(L212)$) {\small (v)};
    \draw[decoration={markings, mark=at position 0.5 with {\arrow{Latex}}},postaction={decorate}] (L24) to (L214);
    \node[right] at ($(L24)!0.5!(L214)$) {\small (i)};
    \draw[decoration={markings, mark=at position 0.5 with {\arrow{Latex}}},postaction={decorate}] (L212) to (L25);
    \node[below] at ($(L212)!0.5!(L25)$) {\small (i)};
    \draw[decoration={markings, mark=at position 0.5 with {\arrow{Latex}}},postaction={decorate}] (L214) to (L25);
    \node[below] at ($(L214)!0.5!(L25)$) {\small (ii)};
    \node[above] at (L0) {$\Lambda^{0}_2$};
    \node[right,xshift=3,yshift=5] at (L2) {$\Lambda^1_{2}$};
    \node[left] at (L21) {$\Lambda^{2}_2$};
    \node[left] at (L212) {$\Lambda^{4}_2$};
    \node[right] at (L24) {$\Lambda^{3}_2$};
    \node[right] at (L214) {$\Lambda^{5}_2$};
    \node[below] at (L25) {$\Lambda^{6}_2$};
  \end{tikzpicture}
\ee

\bigskip\noindent
Case 3) $  5-2p^2 < 4 \Lambda_1 < 5 $

\medskip
Here for A2 Case 3 the complete embedding pattern is simple:

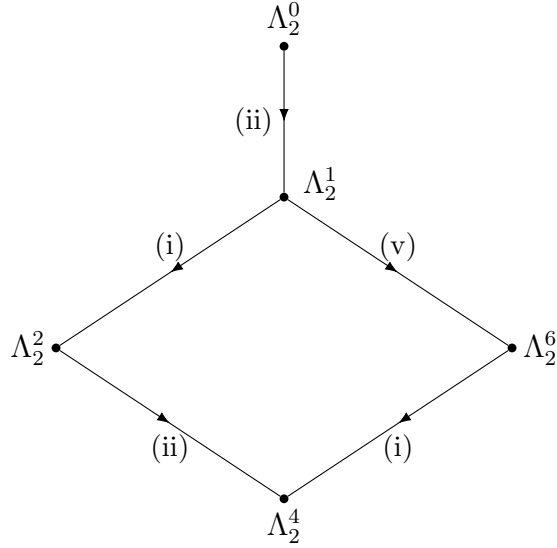
\eqn{A2C3}
  \begin{tikzpicture}
    \node[circle,fill=black,scale=0.3] (L0) at (0,10) {};
    \node[circle,fill=black,scale=0.3] (L2) at (0,8) {};
    \node[circle,fill=black,scale=0.3] (L21) at (-3,6) {};
    \node[circle,fill=black,scale=0.3] (L25) at (3,6) {};
    \node[circle,fill=black,scale=0.3] (L212) at (0,4) {};
    \draw[decoration={markings, mark=at position 0.5 with {\arrow{Latex}}},postaction={decorate}] (L0) to (L2);
    \node[left] at ($(L0)!0.5!(L2)$) {\small (ii)};
    \draw[decoration={markings, mark=at position 0.5 with {\arrow{Latex}}},postaction={decorate}] (L2) to (L21);
    \node[above] at ($(L2)!0.5!(L21)$) {\small (i)};
    \draw[decoration={markings, mark=at position 0.5 with {\arrow{Latex}}},postaction={decorate}] (L2) to (L25);
    \node[above] at ($(L2)!0.5!(L25)$) {\small (v)};
    \draw[decoration={markings, mark=at position 0.5 with {\arrow{Latex}}},postaction={decorate}] (L21) to (L212);
    \node[below] at ($(L21)!0.5!(L212)$) {\small (ii)};
    \draw[decoration={markings, mark=at position 0.5 with {\arrow{Latex}}},postaction={decorate}] (L25) to (L212);
    \node[below] at ($(L25)!0.5!(L212)$) {\small (i)};
    \node[above] at (L0) {$\Lambda^{0}_2$};
    \node[right,xshift=3,yshift=5] at (L2) {$\Lambda^1_{2}$};
    \node[left] at (L21) {$\Lambda^{2}_2$};
    \node[below] at (L212) {$\Lambda^{4}_2$};
    \node[right] at (L25) {$\Lambda^{6}_2$};
  \end{tikzpicture}
\ee

\bigskip\noindent
Case 4) $  5 < 4 \Lambda_1 < 5+2p^2 $

\medskip

For A2 Case 4 the complete embedding pattern is simple:

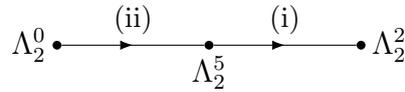
\eqn{A2C4}
  \begin{tikzpicture}
    \node[circle,fill=black,scale=0.3] (L0) at (0,0) {};
    \node[circle,fill=black,scale=0.3] (L2) at (2,0) {};
    \node[circle,fill=black,scale=0.3] (L21) at (4,0) {};
    \draw[decoration={markings, mark=at position 0.5 with {\arrow{Latex}}},postaction={decorate}] (L0) to (L2);
    \node[above] at ($(L0)!0.5!(L2)$) {\small (ii)};
    \draw[decoration={markings, mark=at position 0.5 with {\arrow{Latex}}},postaction={decorate}] (L2) to (L21);
    \node[above] at ($(L2)!0.5!(L21)$) {\small (i)};
    \node[left] at (L0) {$\Lambda^{0}_2$};
    \node[below] at (L2) {$\Lambda^5_{2}$};
    \node[right] at (L21) {$\Lambda^{2}_2$};
  \end{tikzpicture}
\ee

\bigskip
Like \textbf{A1}, one may combine Cases 3) and 4) into  single diagram.

   \bigskip


\subsubsection{Embedding diagram of $ V^{\Lambda^0_4} $ (case A4)}

 We first recall  from  \eqref{L4} the weights in the initial embedding diagram above
   for this case:

   \eqna{a4ini}
     \L^0_4   &=& \Big( \Lambda_1, -\Lambda_1 + 2 -\frac{p^4}{2}  \Big),
     \\
&&     -4\Lambda_1 +5 -p^4 \notin\bbn , \nn\\
    \Lambda^1_4  &=& \Big( \Lambda_1+ \frac{p^4}{2}, -\Lambda_1 + 2 \Big).  \eena

We find that
  $ V^{\Lambda^{1}_4} $ 
has the following SVs:
\begin{itemize}
  \item $ 5 < 4\Lambda_1 \ \Rightarrow \ $ type (ii) SV of weight
    \begin{equation}
     \Lambda^{2}_4  =
      \Big( \Lambda_1+\frac{p^4}{2}, \Lambda_1-\hf \Big).
    \end{equation}
  \item $ 4\Lambda_1 < 5-2p^4 \ \Rightarrow \ $ type (i) (v) SVs of weights:
     \eqna{a44}
     \Lambda^{3}_4   &=&
      \Big( -\Lambda_1+\frac{5}{2}, \Lambda_1 + \hf (p^4-1) \Big),
      \\
    \Lambda^{4}_4   &=&
      \Big( -\Lambda_1-\hf(p^4-5), -\Lambda_1+2 \Big).
     \eena
  \item $ 5-2p^4 < 4\Lambda_1 < 5-p^4 \ \Rightarrow \ $ type (i) SV of weight $ \Lambda^{3}_4  .$
\end{itemize}

 Further embeddings for case A4:

 We find that  for ~$5< 4\Lambda_1$~  the VM  $\Lambda^{2}_4$ has no SV.

 Next we find that if ~$4\Lambda_1 < 5-p^4$~ the VM ~$\Lambda^{3}_4$ has the type (ii) SV of weight
\eqn{a444}
  \Lambda^{5}_4  =
   \Big( -\Lambda_1+\frac{5}{2} , -\Lambda_1-\frac{p^4}{2}+2 \Big).
\end{equation}

 Next we  find that the VM $  \Lambda^{4}_4  $ has the type (i) SV of weight \eqref{a444}.

 Finally, we find that if ~$4\Lambda_1 < 5-2p^4$~ the  VM $ \Lambda^{2}_4   $  has no SV.

\bigskip

 Thus, we find the complete diagram for \textbf{A4}:

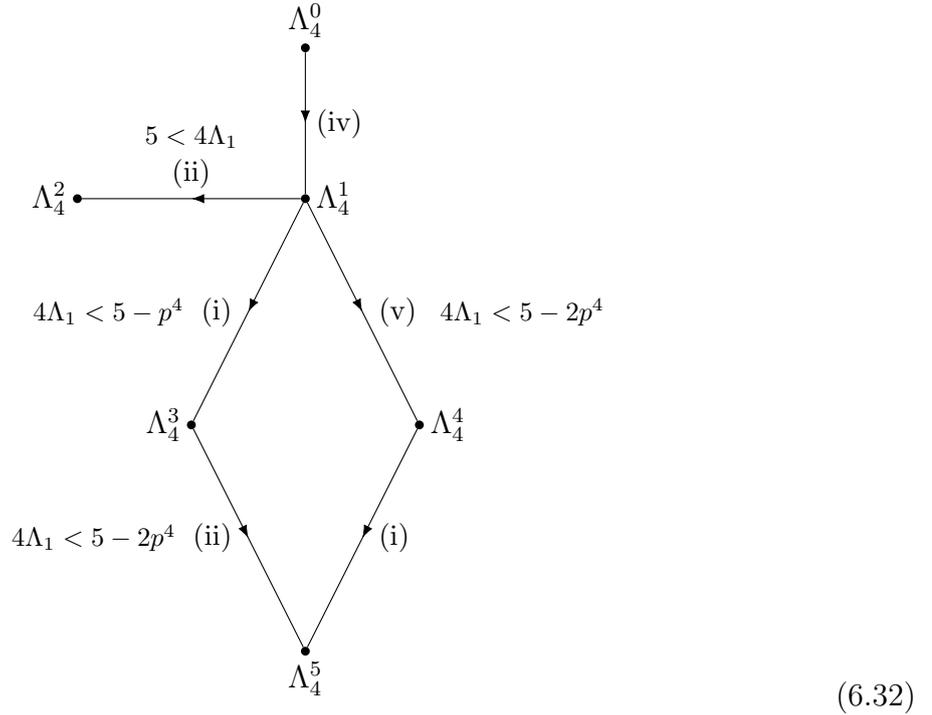
\eqn{A5}
  \begin{tikzpicture}
    \node[circle,fill=black,scale=0.3] (L4) at (0,10) {};
    \node[circle,fill=black,scale=0.3] (Ld4) at (0,8) {};
    \node[circle,fill=black,scale=0.3] (Ldd2) at (-3,8) {};
    \node[circle,fill=black,scale=0.3] (Ldd1) at (-1.5,5) {};
    \node[circle,fill=black,scale=0.3] (Ldd5) at (1.5,5) {};
    \node[circle,fill=black,scale=0.3] (Lddd2) at (0,2) {};
    \draw[decoration={markings, mark=at position 0.5 with {\arrow{Latex}}},postaction={decorate}] (L4) to (Ld4);
    \node[right] at ($(L4)!0.5!(Ld4)$) {\small (iv)};
    \draw[decoration={markings, mark=at position 0.5 with {\arrow{Latex}}},postaction={decorate}] (Ld4) to (Ldd2);
    \node[above] at ($(Ld4)!0.5!(Ldd2)$) {\small (ii)};
    \draw[decoration={markings, mark=at position 0.5 with {\arrow{Latex}}},postaction={decorate}] (Ld4) to (Ldd1);
    \node[left,xshift=-2] at ($(Ld4)!0.5!(Ldd1)$) {\small (i)};
    \draw[decoration={markings, mark=at position 0.5 with {\arrow{Latex}}},postaction={decorate}] (Ld4) to (Ldd5);
    \node[right,xshift=2] at ($(Ld4)!0.5!(Ldd5)$) {\small (v)};
    \draw[decoration={markings, mark=at position 0.5 with {\arrow{Latex}}},postaction={decorate}] (Ldd1) to (Lddd2);
    \node[left,xshift=-2] at ($(Ldd1)!0.5!(Lddd2)$) {\small (ii)};
    \draw[decoration={markings, mark=at position 0.5 with {\arrow{Latex}}},postaction={decorate}] (Ldd5) to (Lddd2);
    \node[right,xshift=2] at ($(Ldd5)!0.5!(Lddd2)$) {\small (i)};
    \node[above] at (L4) {$\Lambda^0_4$};
    \node[right] at (Ld4) {$\Lambda^1_4$};
    \node[left] at (Ldd2) {$ \Lambda^{2}_4  $};
    \node[left] at (Ldd1) {$ \Lambda^{3}_4  $};
    \node[right] at (Ldd5) {$ \Lambda^{4}_4  $};
    \node[below] at (Lddd2) {$\Lambda^{5}_4  $};
    \node[above,yshift=15] at ($(Ld4)!0.5!(Ldd2)$)  {\footnotesize $ 5 < 4\Lambda_1$};
    \node[left,xshift=-20] at ($(Ld4)!0.5!(Ldd1)$) {\footnotesize $ 4\Lambda_1 <5-p^4 $};
    \node[right,xshift=25] at ($(Ld4)!0.5!(Ldd5)$) {\footnotesize $ 4\Lambda_1 <5-2p^4 $};
    \node[left,xshift=-23] at ($(Ldd1)!0.5!(Lddd2)$) {\footnotesize $ 4\Lambda_1 <5-2p^4 $};
  \end{tikzpicture}
\ee

\bigskip


\subsubsection{Embedding diagram of $ V^{\Lambda^0_5} $ (case A5)}

 We first recall  from  \eqref{L5} the weights in the initial embedding diagram above
   for this case:
   \eqna{a5ini}
 \Lambda^0_5 &=& \Big(\frac{5}{4}-\frac{p^5}{2}, \Lambda_2  \Big),   \\
  &&  \Lambda_2 \neq \frac{3}{4}-\frac{r}{2}, \ r \in \mathbb{Z} \nn\\
     \Lambda^1_5  &=& \Big(   \frac{5}{4}+\frac{p^5}{2}, \Lambda_2   \Big) \eena

We find
 that $ V^{\Lambda^1_5} $ has the following SVs:
\begin{itemize}
  \item $ 3+2p^5 < 4\Lambda_2 \ \Rightarrow \ $ type (i) SV of weight
    \begin{equation}
     V^{\Lambda^{2}_5}  = \Big( \Lambda_2+\hf, \frac{3}{4}+\frac{p^5}{2} \Big).
    \end{equation}
  \item $ 4\Lambda_2 < 3-2p^5 \ \Rightarrow \ $ type (ii) (iii) (iv) SVs of weights
      \eqna{a55}
   V^{\Lambda^{3}_5}    &=& \Big( \frac{5}{4}+\frac{p^5}{2},-\Lambda_2+\frac{3}{2} \Big),
      \\
   V^{\Lambda^{5}_5}    &=& \Big( -\Lambda_2+2, \frac{3}{4}+\frac{p^5}{2} \Big),
      \\
   V^{\Lambda^{4}_5}    &=& \Big( -\Lambda_2+2, \frac{3}{4}-\frac{p^5}{2} \Big).
     \eena
  \item $ 3-2p^5 < 4\Lambda_2 < 3 \ \Rightarrow \ $ type (ii) SV of weight $\Lambda^{3}_5  . $
\end{itemize}


 Next we find that if ~$3+2p^5 < 4\Lambda_2$~ the   VM $ \Lambda^{2}_5  $  has no SV.

Next we find that if ~$4\Lambda_2 < 3$~ the VM  $ \Lambda^{3}_5  $
  has the type (i) SV of weight  \rf{a55}{b}.

Next we find that if ~$4\Lambda_2 < 3-2p^5$~ the VM $ \Lambda^{5}_5  $
    has no SV.

Next we find that if~ $4\Lambda_2 < 3-2p^5$~ the VM $ \Lambda^{4}_5  $
has the type (ii) SV of weight  \rf{a55}{b}.

Combining all results for case A5 we find that it has the following diagram:

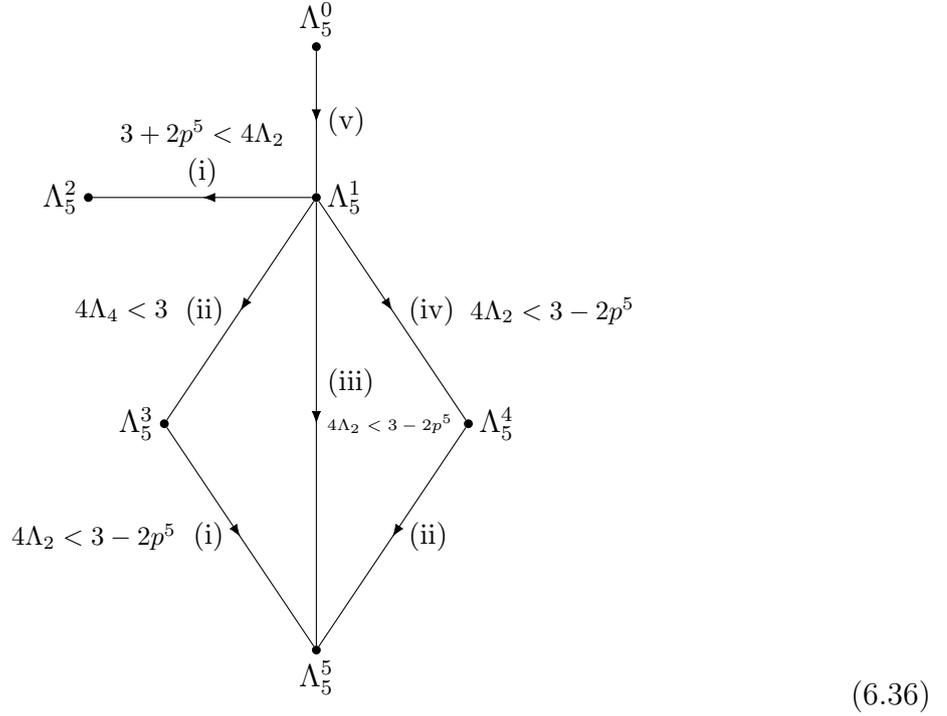
\eqn{A5dg}
  \begin{tikzpicture}
    \node[circle,fill=black,scale=0.3] (L5) at (0,10) {};
    \node[circle,fill=black,scale=0.3] (Ld5) at (0,8) {};
    \node[circle,fill=black,scale=0.3] (Ldd1) at (-3,8) {};
    \node[circle,fill=black,scale=0.3] (Ldd2) at (-2,5) {};
    \node[circle,fill=black,scale=0.3] (Ldd3) at (0,2) {};
    \node[circle,fill=black,scale=0.3] (Ldd4) at (2,5) {};
    \draw[decoration={markings, mark=at position 0.5 with {\arrow{Latex}}},postaction={decorate}] (L5) to (Ld5);
    \node[right] at ($(L5)!0.5!(Ld5)$) {\small (v)};
    \draw[decoration={markings, mark=at position 0.5 with {\arrow{Latex}}},postaction={decorate}] (Ld5) to (Ldd1);
    \node[above] at ($(Ld5)!0.5!(Ldd1)$) {\small (i)};
    \draw[decoration={markings, mark=at position 0.5 with {\arrow{Latex}}},postaction={decorate}] (Ld5) to (Ldd2);
    \node[left,xshift=-2] at ($(Ld5)!0.5!(Ldd2)$) {\small (ii)};
    \draw[decoration={markings, mark=at position 0.5 with {\arrow{Latex}}},postaction={decorate}] (Ld5) to (Ldd3);
    \node[right] at ($(Ld5)!0.41!(Ldd3)$) {\small (iii)};
    \draw[decoration={markings, mark=at position 0.5 with {\arrow{Latex}}},postaction={decorate}] (Ld5) to (Ldd4);
    \node[right,xshift=2] at ($(Ld5)!0.5!(Ldd4)$) {\small (iv)};
    \draw[decoration={markings, mark=at position 0.5 with {\arrow{Latex}}},postaction={decorate}] (Ldd4) to (Ldd3);
    \node[right,xshift=2] at ($(Ldd4)!0.5!(Ldd3)$) {\small (ii)};
    \draw[decoration={markings, mark=at position 0.5 with {\arrow{Latex}}},postaction={decorate}] (Ldd2) to (Ldd3);
    \node[left,xshift=-2] at ($(Ldd2)!0.5!(Ldd3)$) {\small (i)};
    \node[above] at (L5) {$\Lambda^0_5$};
    \node[right] at (Ld5) {$\Lambda^1_5$};
    \node[left] at (Ldd1) {$\Lambda^{2}_5  $};
    \node[left] at (Ldd2) {$\Lambda^{3}_5  $};
    \node[below] at (Ldd3) {$\Lambda^{5}_5 $};
    \node[right] at (Ldd4) {$\Lambda^{4}_5 $};
    \node[above,yshift=15] at ($(Ld5)!0.5!(Ldd1)$) {\footnotesize $ 3+2p^5 < 4\Lambda_2 $};
    \node[left,xshift=-23] at ($(Ld5)!0.5!(Ldd2)$) {\footnotesize $ 4\Lambda_4 < 3 $};
    \node[right,xshift=25] at ($(Ld5)!0.5!(Ldd4)$) {\footnotesize $ 4 \Lambda_2 < 3-2p^5$};
    \node[right] at (0,5) {\tiny $ 4 \Lambda_2 < 3-2p^5$};
    \node[left,xshift=-20] at ($(Ldd2)!0.5!(Ldd3)$) {\footnotesize $ 4 \Lambda_2 < 3-2p^5$};
  \end{tikzpicture}
\ee

\bigskip

\section*{Acknowledgments}

\nt
VKD acknowledges partial support from Bulgarian NSF Grant DN-18/1. 
The authors are grateful to the referees for valuable comments. 

 \bigskip

%
%
%

\end{document}